\documentclass[a4paper,12pt,twoside]{article}
\usepackage[pdftex,colorlinks,bookmarksnumbered,bookmarksopen,citecolor=red,urlcolor=red]{hyperref}
\usepackage{graphicx}\usepackage{a4wide}
\usepackage{amssymb,amsmath,amsthm}
\usepackage{subfigure,color,colordvi,graphicx,xcolor}
\usepackage[only,llbracket,rrbracket,interleave]{stmaryrd}
\usepackage[superscript,sort,compress]{cite}

\newcommand{\bm}[1]{\text{\boldmath $#1$\unboldmath}}

\newcommand{\vect}[1]{\mathbf{#1}}
\newcommand{\mat}[1]{\mathbf{#1}}

\newcommand{\grad}{\bm{\nabla}}

\newcommand{\RR}{\mathbb{R}}

\newcommand{\VhHat}{\ensuremath{\mathcal{\hat{V}}^h}}
\newcommand{\Vh}{\ensuremath{\mathcal{V}^h}}

\newcommand{\eltwo}{\ensuremath{\mathcal{L}_2}}

\newcommand{\ndof}  {\ensuremath{\texttt{n}_{\texttt{dof}}}}

\newcommand{\nsd}  {\ensuremath{\texttt{n}_{\texttt{sd}}}}
\newcommand{\numel}{\ensuremath{\texttt{n}_{\texttt{el}}}}
\newcommand{\numfa}{\ensuremath{\texttt{n}_{\texttt{fa}}}}

\newcommand{\hv}{\hat{v}}
\newcommand{\hu}{\hat{u}}
\newcommand{\bu}{\bm{u}}

\newcommand{\bhw}{\widehat{\bw}}
\newcommand{\bhu}{\widehat{\bu}}
\newcommand{\bq}{\bm{q}}

\newcommand{\bn}{\bm{n}}
\newcommand{\bx}{\bm{x}}
\newcommand{\bL}{\bm{L}}

\newcommand{\qe}{\vect{q}_e}
\newcommand{\ue}{\text{u}_e}
\newcommand{\uHi}{\hat{\text{u}}_i}
\newcommand{\uHj}{\hat{\text{u}}_j}

\newcommand{\Le}{\mat{L}_e}
\newcommand{\bue}{\vect{u}_e}
\newcommand{\pe}{\text{p}_e}
\newcommand{\buHi}{\hat{\vect{u}}_i}
\newcommand{\buHj}{\hat{\vect{u}}_j}

\newcommand{\nDeg}{\ensuremath{k}}
\newcommand{\Pk}{\ensuremath{\mathcal{P}^{\nDeg}}}

\newcommand{\Insd}{\mat{I}_{\nsd}}

\newcommand{\jump}[1]{\llbracket #1\rrbracket}

\newenvironment{keywords}{\begin{quote}\emph{\textbf{Keywords:}}}{\end{quote}}

\theoremstyle{definition}

\newtheorem{remark}{Remark}

\newcommand{\bw}{\bm{w}}

\newcommand{\Aset}{\mathcal{A}_e}
\newcommand{\Bset}{\mathcal{B}_e}
\newcommand{\Dset}{\mathcal{D}_e}
\newcommand{\Nset}{\mathcal{N}_e}
\newcommand{\Mset}{\mathcal{M}_e}
\newcommand{\Iset}{\mathcal{I}_e}

\graphicspath{{figsPDF/}{figsPNG/}}

\begin{document}
\title{A face-centred finite volume method for second-order elliptic problems}

\author{Ruben Sevilla\\[-1ex]
             \small Zienkiewicz Centre for Computational Engineering, \\[-1ex]
             \small College of Engineering, Swansea University, Wales, UK \\[1em]
             Matteo Giacomini, and Antonio Huerta\\[-1ex]
             \small Laboratori de C\`alcul Num\`eric (LaC\`aN), \\[-1ex]
             \small ETS de Ingenieros de Caminos, Canales y Puertos, \\[-1ex]
             \small Universitat Polit\`ecnica de Catalunya, Barcelona, Spain}
\date{\today}
\maketitle

\begin{abstract}
This work proposes a novel finite volume paradigm, the face-centred finite volume (FCFV) method. Contrary to the popular vertex (VCFV) and cell (CCFV) centred finite volume methods, the novel FCFV defines the solution on the mesh faces (edges in 2D) to construct locally-conservative numerical schemes. 
The idea of the FCFV method stems from a hybridisable discontinuous Galerkin (HDG) formulation with constant degree of approximation, thus inheriting the convergence properties of the classical HDG. The resulting FCFV features a global problem in terms of a piecewise constant function defined on the faces of the mesh. The solution and its gradient in each element are then recovered by solving a set of independent element-by-element problems.
The mathematical formulation of FCFV for Poisson and Stokes equation is derived and numerical evidence of optimal convergence in 2D and 3D is provided. Numerical examples are presented to illustrate the accuracy, efficiency and robustness of the proposed methodology. The results show that, contrary to other FV methods, the accuracy of the FCFV method is not sensitive to mesh distortion and stretching. In addition, the FCFV method shows its better performance, accuracy and robustness using simplicial elements, facilitating its application to problems involving complex geometries in 3D.
\end{abstract}


\section{Introduction}
\label{sc:Intro}

Starting from its first appearance in the 1960s, the finite volume method (FVM) has experienced a growing success, especially within the computational fluid dynamics (CFD) community. Stemming from the fundamental work of Godunov~\cite{MR0119433}, Varga~\cite{MR1753713} and Preissmann~\cite{preissmann1961propagation}, the FVM made its official appearance in~\cite{McDonald71,RizziInouye73}, where the authors considered its application to hyperbolic systems of conservation laws.  Nowadays, the FVM is the most widely spread methodology implemented in open-source, commercial and industrial CFD solvers. 

Over the years, several variations of the original FVM have been proposed. For a complete introduction, the interested reader is referred to classical textbooks~\cite{MR1925043,MR2731357} and to the review papers~\cite{MR2417929,BO2004,EYMARD2000713,droniou2014review}. The two most popular approaches are the so-called cell-centred finite volume (CCFV) method and the vertex-centred finite volume (VCFV) method. 

The CCFV approach defines the solution at the centre of the mesh elements (i.e. cells) such that their values represent cell averages of the unknown quantities. Several techniques to accurately compute the gradient of the solution at the element faces, based on node averaging or least squares, have been proposed and compared~\cite{diskin2010comparison,diskin2011comparison}. In all cases, a reconstruction of the gradient of the solution is required to guarantee second-order convergence of the solution error. This is crucial to guarantee a first-order convergence of the solution gradient, which is required to accurately compute engineering quantities of interest (e.g. lift and drag). The accuracy of the reconstruction is heavily dependent on the quality of the mesh and some approaches fail to provide a second-order scheme on highly stretched and deformed grids.

The VCFV strategy defines the solution at the mesh nodes. A control volume is constructed around each node by using the centroid of the mesh elements and mid-edge points (and face centroids in three dimensions). The control volumes form a non-overlapping set of subdomains that cover the whole domain and form the so-called dual mesh. The resulting approximation is locally piecewise constant on each dual element where the values of the unknowns represent control volume averages. Similarly to CCFV scheme, the VCFV method requires the reconstruction of the gradient of the solution at each dual face. A first order accurate reconstruction scheme is required to provide a second-order VCFV method~\cite{diskin2010comparison,diskin2011comparison}.

In parallel to the development of FV schemes, a great effort was dedicated during the 1970s to the application of finite element methods to CFD problems~\cite{zienkiewicz2000finite}. The difficulties encountered due to the convection dominated nature of many fluid flow problems prompted the development of the so-called stabilised finite element techniques~\cite{donea2003finite} and discontinuous Galerkin methods~\cite{ReedHill:1973}.

More recently, a great effort has been dedicated to reinterpret finite volume schemes within a continuous and discontinuous finite element framework. In \cite{MR2417929}, Morton and Sonar motivate their exposition of finite volume schemes as Petrov-Galerkin finite element methods owing to the flexibility the latter approaches show in terms of approximation using unstructured meshes and the solid theoretical framework developed for their analysis. In \cite{MR2765501}, Vohral\'\i k exploits similar ideas to develop a unified theory of \emph{a posteriori} error estimators valid for both finite volume and finite element approximations. Within this context, the CCFV scheme may be interpreted as a discontinuous Galerkin method with piecewise constant degree of approximation within each element~\cite{EYMARD2000713,MR1842161}. In a similar fashion, a VCFV scheme on simplicial meshes may be interpreted as a conforming piecewise linear continuous finite element method~\cite{MR899703,selmin1993node,idelsohn1994finite}.

In this paper, an alternative to the discussed finite volume strategies is proposed by defining the unknowns over the faces of the mesh. 
As for CCFV and VCFV, the resulting face-centred finite volume (FCFV) method may be interpreted as a lowest-order finite element method. More precisely, FCFV is derived from the recently proposed hybridisable discontinuous Galerkin (HDG) method by Cockburn and co-workers \cite{cockburn2004characterization,Jay-CG:05,Jay-CG:05-GAMM,Jay-CGL:09} by imposing a constant degree of approximation. As such, the method requires the solution of a global system of equations equal to the total number of element faces. The solution and its gradient in each element are then recovered by solving a set of independent element-by-element problems.

The proposed FCFV method provides first-order accuracy on both the solution and its gradient without the need to perform a reconstruction of the gradients to accurately compute the fluxes at the element or control volume boundary. Therefore, its accuracy is not compromised in the presence of highly stretched or distorted elements. In addition, due to the definition of the unknowns on the element faces, the global system of equations that must be solved, provides a less degree of coupling of the information when compared to other finite volume schemes. The application to scalar and vector second-order elliptic problems is considered, namely the Poisson and the Stokes problems respectively. For the solution of Stokes flow problems, the FVFC method does not require the solution of a Poisson problem for computing the pressure, as required by segregated schemes such as the semi-implicit method for pressure-linked equations (SIMPLE) algorithm~\cite{patankar1980numerical}. In addition, contrary to other mixed finite element methods, with the FCVC it is possible to use the same space of approximation for both velocity and pressure, circumventing the so-called Ladyzhenskaya-Babu{\v s}ka-Brezzi  (LBB) condition~\cite{donea2003finite}.

The rest of this paper is organised as follows. Section~\ref{sc:Poisson} presents the proposed FCFV method for the solution of the Poisson equation. The extension to Stokes flow problems is described in Section~\ref{sc:Stokes}. Section \ref{sc:computational} discusses some computational aspects of the FCFV rationale and recalls its theoretical convergence properties. An exhaustive set of numerical studies is presented in Section~\ref{sc:studies}. These studies include mesh convergence tests, a comparison in terms of the computational cost and the influence of the stabilisation parameter, the  mesh distortion and the element stretching. The studies consider both the Poisson and Stokes equations, using different element types and in two and three dimensional domains. In Section~\ref{sc:examples} large three dimensional examples are considered to show the potential of the proposed methodology. Finally, section~\ref{sc:Conclusion} summarises the conclusions of the work that has been presented.

\section{FCFV for the Poisson equation}
\label{sc:Poisson}

\subsection{Problem statement and mixed formulation}
\label{sc:PoissonStatement}

Let $\Omega\in\mathbb{R}^{\nsd}$ be an open bounded domain with boundary $\partial\Omega = \overline{\Gamma}_D \cup \overline{\Gamma}_N$, $\overline{\Gamma}_D \cap \overline{\Gamma}_N = \emptyset$ and $\nsd$ the number of spatial dimensions. The strong form for the second-order elliptic problem can be written as
\begin{equation} \label{eq:Poisson}
\left\{\begin{aligned}
-\grad\cdot\grad u &= s       &&\text{in $\Omega$,}\\
u &= u_D  &&\text{on $\Gamma_D$,}\\
\bn\cdot\grad u &= t        &&\text{on $\Gamma_N$,}\\
\end{aligned}\right.
\end{equation}
where $s\in\eltwo(\Omega)$ is a source term, $\bn$ is the outward unit normal vector to $\partial\Omega$ and $u_D$ and $t$ respectively are the Dirichlet and Neumann data imposed on the external boundary. Other boundary conditions may also be considered but, for the sake of simplicity  (and without any loss of generality), solely the Dirichlet-Neumann case will be detailed.

Let us assume that $\Omega$ is partitioned in $\numel$ disjoint subdomains $\Omega_e$
\begin{equation}
\overline{\Omega} =  \bigcup_{e=1}^{\numel} \overline{\Omega}_e, \quad 
\Omega_e \cap \Omega_l = \emptyset \text{ for } e\neq l ,
\end{equation}
with boundaries $\partial\Omega_e$, which define an internal interface $\Gamma$
\begin{equation}\label{eq:Gamma}
\Gamma := \Big[ \bigcup_{e=1}^{\numel} \partial\Omega_e \Big]\setminus\partial\Omega
\end{equation}
Moreover, it is also convenient to write the boundary of each element as the union of the individual element faces (edges in two dimensions), namely
\begin{equation}
\partial\Omega_e :=  \bigcup_{j=1}^{\numfa^e} \Gamma_{e,j},
\end{equation}
where $\numfa^e$ denotes the number of faces of the element $\Omega_e$.

Following the definition in \cite{AdM-MFH:08}, the \emph{jump} $\jump{\cdot}$ operator is introduced. That is, along each portion of the interface $\Gamma$ it sums the values from the left and right of say, $\Omega_e$ and $\Omega_l$, namely 
\begin{equation}
\jump{\odot} = \odot_e + \odot_l .
\end{equation}
It is important to observe that this definition always requires the normal vector $\bn$ in the argument and always produces functions in the same space as the argument.

The second-order elliptic problem \eqref{eq:Poisson} can thus be written in mixed form in the \emph{broken} computational domain as a system of first-order equations, namely
\begin{equation} \label{eq:PoissonBrokenFirstOrder}
\left\{\begin{aligned}
\bq+\grad u &= \bm{0} &&\text{in $\Omega_e$, and for $e=1,\dotsc ,\numel$,}\\	
\grad\cdot\bq &= s          &&\text{in $\Omega_e$, and for $e=1,\dotsc ,\numel$,}\\
u &= u_D     &&\text{on $\Gamma_D$,}\\
\bn\cdot\bq &= -t         &&\text{on $\Gamma_N$,}\\
\jump{u\bn} &=\bm{0}  &&\text{on $\Gamma$,}\\
\jump{\bn\cdot \bq} &= 0  &&\text{on $\Gamma$,}\\
\end{aligned} \right.
\end{equation}
where the two last equations correspond to the imposition of the continuity of respectively the primal variable $u$ and the normal fluxes along the internal interface $\Gamma$.

\subsection{Strong form of the local and global problems}
\label{sc:PoissonStrong}

In this subsection, the classical formulation of the hybridisable discontinuous Galerkin method is recalled. The HDG method for second-order elliptic problems has been studied in a series of papers by Cockburn and co-workers\cite{Jay-CGL:09,cockburn2009hybridizable,Nguyen-NPC:09} and relies on rewriting Equation~\eqref{eq:PoissonBrokenFirstOrder} as two equivalent problems. First, the local (element-by-element) problem with Dirichlet boundary conditions is defined, namely 
\begin{equation} \label{eq:Dlocal-strong}
\left\{\begin{aligned}
\bq_e + \grad u_e &=\bm{0}  &&\text{in $\Omega_e$, }\\
\grad\cdot\bq_e &= s          &&\text{in $\Omega_e$,}\\		
u_e &= u_D     &&\text{on $\partial\Omega_e\cap\Gamma_D$,}\\
u_e &=\hu  &&\text{on $\partial\Omega_e\setminus\Gamma_D$,}
\end{aligned} \right.
\end{equation}
for $e=1,\dotsc ,\numel$.
In each element $\Omega_e$ this problem produces an element-by-element solution $\bq_e$ and $u_e$ as a function of the unknown $\hu\in\eltwo(\Gamma\cup\Gamma_N)$. Note that these problems can be solved independently element-by-element.

Second, a global problem is defined to determine $\hu$. It corresponds to the imposition of the Neumann boundary condition and the so-called \emph{transmission conditions}, see \cite{Jay-CGL:09}. 
\begin{equation} \label{eq:Dtransmission}
\left\{\begin{aligned}
\jump{u\bn} &=\bm{0}  &&\text{on $\Gamma$,}\\
\jump{\bn\cdot \bq} &= 0  &&\text{on $\Gamma$,}\\
\bn\cdot\bq &= -t         &&\text{on $\Gamma_N$.}\\
\end{aligned} \right.
\end{equation}
These transmission conditions were introduced in \eqref{eq:PoissonBrokenFirstOrder} to ensure inter-element continuity when the broken computational domain formulation was presented. Note that the first equation in \eqref{eq:Dtransmission} imposes continuity of $u$ across $\Gamma$. Owing to the Dirichlet boundary condition $u=\hu$ on $\Gamma$ as imposed by the local problems \eqref{eq:Dlocal-strong} and the uniqueness of the hybrid variable $\hu$, 
the continuity of the primal variable, $\jump{\hu\bn} =\bm{0}$, is automatically verified. Hence, the global problem reduces to the second and third equation in \eqref{eq:Dtransmission}.

\subsection{Weak form of the local and global problems}
\label{sc:PoissonWeak}

First, following the notation in \cite{RS-SH:16}, the discrete functional spaces are introduced:
\begin{subequations}\label{eq:HDG-Poisson-Spaces}
\begin{align} 
\Vh(\Omega) & := \{ v \in \eltwo(\Omega) : v \vert_{\Omega_e}\in \Pk(\Omega_e) \;\forall\Omega_e \, , \, e=1,\dotsc ,\numel \} , \label{eq:spaceScalarElem} \\
\VhHat(S) & := \{ \hv \in \eltwo(S) : \hv\vert_{\Gamma_i}\in \Pk(\Gamma_i)
\;\forall\Gamma_i\subset S\subseteq\Gamma\cup\partial\Omega \}, \label{eq:spaceScalarFace}
\end{align}
\end{subequations}
where $\mathcal{P}^{k}(\Omega_e)$ and $\mathcal{P}^{k}(\Gamma_i)$ stand for the spaces of polynomial functions of complete degree at most $k$ in $\Omega_e$ and on $\Gamma_i$ respectively. 
Moreover, recall the notation for the classical internal products of scalar functions in $\eltwo(\Omega_e)$ and $\eltwo(\Gamma_i)$ 
\begin{equation} \label{eq:innerScalar}
(p,q)_{\Omega_e} := \int_{\Omega_e} p q \ d\Omega , \quad \langle \hat{p}, \hat{q} \rangle_{\partial\Omega_e} := \sum_{\Gamma_i \subset \partial\Omega_e} \int_{\Gamma_i} \hat{p} \hat{q} \ d\Gamma 
\end{equation}
and the internal product of vector valued functions in $[\eltwo(\Omega)]^{\nsd}$
\begin{equation} \label{eq:innerVector}
	(\bm{p},\bm{q})_{\Omega_e} := \int_{\Omega_e} \bm{p} \cdot\bm{q} \ d\Omega 
\end{equation}

The discrete weak formulation of the previously introduced local problems is obtained by multiplying the problems by a test function in an appropriate discrete functional space and integrating by parts.
For $e=1,\dotsc ,\numel$, seek $(\bq_e^h ,u_e^h)\in [\Vh(\Omega_e)]^{\nsd}\times\Vh(\Omega_e)$ such that for all $(\bw ,v)\in [\Vh(\Omega_e)]^{\nsd}\times\Vh(\Omega_e)$ it holds
\begin{gather*} 
- (\bw,\bq_e^h)_{\Omega_e} 
+ (\grad \cdot \bw, u_e^h)_{\Omega_e} 
= \langle \bn_e\cdot\bw, u_D \rangle_{\partial\Omega_e \cap \Gamma_D}
+ \langle \bn_e\cdot\bw, \hu^h \rangle_{\partial\Omega_e\setminus\Gamma_D} , 
\\	
- (\grad v, \bq_e^h)_{\Omega_e} 
+\langle v,\bn_e\cdot\widehat{\bq}_e^h \rangle_{\partial\Omega_e}  = (v,s)_{\Omega_e} ,
\end{gather*}
The traces of the numerical fluxes $\widehat{\bq}_e^h$ have to be properly defined in order to guarantee the stability of the method \cite{Jay-CGL:09}. More precisely, they are defined element-by-element (i.e.\ for $e=1,\dotsc ,\numel$) as  
\begin{equation} \label{eq:EBENumFlux}
\bn_e\cdot\widehat{\bq}_e^h := \begin{cases}
\bn_e\cdot\bq_e^h + \tau_e (u_e^h - u_D    ) & \text{on $\partial\Omega_e\cap\Gamma_D$,} \\
\bn_e\cdot\bq_e^h + \tau_e (u_e^h - \hu^h) & \text{elsewhere,}  
\end{cases}
\end{equation}
with $\tau_e$ being a stabilization parameter which may assume different values on each face of the boundary $\partial\Omega_e$ and whose selection has an important effect on the stability, accuracy and convergence properties of the resulting HDG method. The influence of the stabilization parameter has been studied extensively by Cockburn and co-workers, see for instance~\cite{Jay-CGL:09,Cockburn-CDG:08}. \\
By exploiting the definition of the numerical fluxes given by \eqref{eq:EBENumFlux} and integrating by parts again the left-hand side of the second equation in order to retrieve a symmetric formulation, the discrete weak problem becomes:
for $e=1,\dotsc ,\numel$, seek $(\bq_e^h ,u_e^h)\in [\Vh(\Omega_e)]^{\nsd}\times\Vh(\Omega_e)$ that for all $(\bw ,v)\in [\Vh(\Omega_e)]^{\nsd}\times\Vh(\Omega_e)$ satisfies 
\begin{subequations}\label{eq:HDG-Poisson-Dlocal}
	\begin{gather} 
	- (\bw,\bq_e^h)_{\Omega_e} 
	+ (\grad\cdot\bw , u_e^h)_{\Omega_e}
	= \langle \bn_e\cdot\bw, u_D \rangle_{\partial\Omega_e \cap \Gamma_D}
	+ \langle \bn_e\cdot\bw, \hu^h \rangle_{\partial\Omega_e\setminus\Gamma_D} ,
	\label{eq:Dweak-continuous2}
	\\
	(v, \grad \cdot \bq_e^h )_{\Omega_e}
	+ \langle v,\tau_e\, u_e^h \rangle_{\partial\Omega_e}  
	= (v,s)_{\Omega_e}
	+ \langle v,\tau_e\, u_D \rangle_{\partial\Omega_e\cap\Gamma_D}  
	+ \langle v,\tau_e\,\hu^h \rangle_{\partial\Omega_e\setminus\Gamma_D}  .
	\label{eq:Dweak-continuous1}		
	\end{gather}
\end{subequations}

In a similar fashion, the following discrete formulation is derived for the global problem: seek $\hu^h\in\VhHat(\Gamma\cup\Gamma_N)$ such that for all $\hv\in\VhHat(\Gamma\cup\Gamma_N)$ it holds
\begin{equation}
\sum_{e=1}^{\numel}
\langle \hv,\bn_e\cdot\widehat{\bq}_e^h \rangle_{\partial\Omega_e\setminus\partial\Omega}
+
\sum_{e=1}^{\numel} 
\langle \hv,\bn_e\cdot\widehat{\bq}_e^h + t  \rangle_{\partial\Omega_e\cap\Gamma_N}
= 0 ,
\end{equation}
or, equivalently, 
\begin{equation}\label{eq:HDG-Poisson-Dglobal}
\sum_{e=1}^{\numel}\Bigl\{
\langle \hv,\bn_e\cdot\bq_e^h \rangle_{\partial\Omega_e\setminus\Gamma_D}
+ \langle \hv,\tau_e\, u_e^h \rangle_{\partial\Omega_e\setminus\Gamma_D}
- \langle \hv,\tau_e\,\hu^h \rangle_{\partial\Omega_e\setminus\Gamma_D}\Bigr\}
= -\sum_{e=1}^{\numel} \langle \hv, t \rangle_{\partial\Omega_e\cap\Gamma_N}.
\end{equation}

Henceforth, to simplify the notation the superindex $^h$ expressing the discrete approximations will be dropped, unless needed in order to follow the development.

\subsection{FCFV discretisation}
\label{sc:PoissonFV}

For the sake of readability, introduce the following notation for the sets of faces: $\Aset := \{1, \ldots, \numfa^e \}$ is the set of indices for all the faces of element $\Omega_e$; $\Dset := \{j \in \Aset \; | \; \Gamma_{e,j} \cap \Gamma_D \neq \emptyset \}$ is the set of indices for all the faces of element $\Omega_e$ on the Dirichlet boundary $\Gamma_D$; $\Nset := \{j \in \Aset \; | \; \Gamma_{e,j} \cap \Gamma_N \neq \emptyset \}$ is the set of indices for all the faces of element $\Omega_e$ on the Neumann boundary $\Gamma_N$; $\Bset := \Aset \setminus \Dset = \{j \in \Aset \; | \; \Gamma_{e,j} \cap \Gamma_D = \emptyset \}$ is the set of indices for all the faces of element $\Omega_e$ not on the Dirichlet boundary $\Gamma_D$.

The discretisation of the local problems given by Equation~\eqref{eq:HDG-Poisson-Dlocal} with a degree of approximation $k=0$ in each element for both $\bq_e$ and $u_e$ and also a degree of approximation $k=0$ in each face/edge for $\hu$ leads to 
\begin{subequations}\label{eq:HDG-Poisson-DlocalK0}
	\begin{equation}
	- |\Omega_e| \qe =  \sum_{j\in\Dset} |\Gamma_{e,j}| \bn_j u_{D,j}  + \sum_{j\in\Bset} |\Gamma_{e,j}| \bn_j \uHj  ,
	\end{equation}
	\begin{equation}
	\sum_{j\in\Aset} |\Gamma_{e,j}| \tau_j \ue = |\Omega_e| s_e + \sum_{j\in\Dset} |\Gamma_{e,j}| \tau_j u_{D,j} + \sum_{j\in\Bset} |\Gamma_{e,j}| \tau_j \uHj
	\end{equation}
\end{subequations}
where $\tau_j$ is the value of the stabilization parameter on the $j$-th face of the element.

\begin{remark}
The local problem of Equation~\eqref{eq:HDG-Poisson-DlocalK0} assumes that the integrals in the weak formulation are computed with a numerical quadrature with a single integration point in each element and each face/edge. It is worth noting that this is exact for some of the integrals but in other cases it introduces an error of order $\mathcal{O}(h)$, where $h$ is the characteristic element size. The two situations where the integral with one integration point is not exact are when the data (source term, Dirichlet and Neumann data) that is not constant per element or face/edge and when the outward unit normal to the face changes within the face (i.e. when elements with non-planar quadrilateral faces are considered in three dimensions).
\end{remark}

In general, the HDG local problem requires the solution of a small system of equations to obtain $\bq$ and $u$ in terms of $\hu$. However, for the particular choice of a constant interpolation in each element, Equations~\eqref{eq:HDG-Poisson-DlocalK0} are uncoupled and provide the following explicit expressions of $\bq_e$ and $u_e$ in terms of $\hu$:
\begin{subequations}\label{eq:HDG-Poisson-DlocalK0explicit}
	\begin{equation}
	\qe =  -|\Omega_e|^{-1} \vect{z}_e  -|\Omega_e|^{-1} \sum_{j\in\Bset} |\Gamma_{e,j}| \bn_j \uHj  
	\label{eq:Poisson-locQ}
	\end{equation}
	\begin{equation}
	\ue = \alpha_e^{-1} \beta_e + \alpha_e^{-1} \sum_{j\in\Bset} |\Gamma_{e,j}| \tau_j \uHj,
	\label{eq:Poisson-locU}
	\end{equation}
\end{subequations}
where
\begin{equation} \label{eq:poissonPrecomp}
\alpha_e := \sum_{j\in\Aset} |\Gamma_{e,j}| \tau_j 
, \quad
\beta_e := |\Omega_e| s_e  + \sum_{j\in\Dset} |\Gamma_{e,j}| \tau_j u_{D,j} 
, \quad
\vect{z}_e := \sum_{j\in\Dset} |\Gamma_{e,j}| \bn_j u_{D,j} .
\end{equation}

Similarly, the discretisation of the global problem given by Equation \eqref{eq:HDG-Poisson-Dglobal} with a degree of approximation $k=0$ for $\hu$ leads to
\begin{equation}\label{eq:HDG-Poisson-DglobalK0}
\sum_{e=1}^{\numel}\Bigl\{
|\Gamma_{e,i}| \bn_i \cdot \qe + |\Gamma_{e,i}| \tau_i \ue - |\Gamma_{e,i}| \tau_i \uHi \Bigr\}
= -\sum_{e=1}^{\numel} \Bigl\{|\Gamma_{e,i}| t_i \, \chi_{\Nset}(i) \Bigr\} \quad \text{for $i \in \Bset$} ,
\end{equation}
where $\chi_{\Nset}$ is the indicator function of the set $\Nset$, i.e.
\begin{equation}
\chi_{\Nset}(i) = 
\bigg\{
\begin{array}{ll}
1 & \text{ if } \ i\in\Nset \\
0 & \text{ otherwise}
\end{array}.
\end{equation}
By plugging \eqref{eq:Poisson-locQ} and \eqref{eq:Poisson-locU} into \eqref{eq:HDG-Poisson-DglobalK0}, the following system of equations containing only $\hu$ as an unknown is obtained:
\begin{equation}
\mat{\widehat{K}}\vect{\hu}=\vect{\hat{f}},
\end{equation}
where the global matrix $\mat{\widehat{K}}$ and right hand side vector $\vect{\hat{f}}$ are computed by assembling the contributions given by
\begin{subequations}\label{HDG-Poisson-globalSystem}
	\begin{align}
	&
	\mat{\widehat{K}}^e_{i,j} :=  |\Gamma_{e,i}| \left( \alpha_e^{-1} |\Gamma_{e,j}| \tau_i \tau_j - |\Omega_e|^{-1} |\Gamma_{e,j}| \bn_i \cdot \bn_j - \tau_i \delta_{ij} \right) 
	, \\
	&
	\vect{\widehat{f}}^e_i :=  |\Gamma_{e,i}| \left( |\Omega_e|^{-1}  \bn_i \cdot \vect{z}_e - t_i \, \chi_{\Nset}(i) - \alpha_e^{-1} \beta_e \tau_i  \right),
	\end{align}
\end{subequations}
for $i,j \in \Bset$, being $\delta_{ij}$ the classical Kronecker delta.

\section{FCFV for the Stokes equation}
\label{sc:Stokes}

\subsection{Problem statement and mixed formulation}
\label{sc:StokesStatement}

Following the above rationale, the formulation of the face-centred finite volume method for the approximation of Stokes flows is derived next. The strong form  of the velocity-pressure formulation of the Stokes equation can be written as
\begin{equation} \label{eq:Stokes}
\left\{\begin{aligned}
- \grad\cdot(\nu \grad \bu -  p \Insd) &= \bm{s}       &&\text{in $\Omega$,}\\
\grad\cdot\bu &= 0  &&\text{in $\Omega$,}\\
\bu &= \bu_D  &&\text{on $\Gamma_D$,}\\
\bn \cdot \bigl(\nu \grad \bu  -  p \Insd \bigr) &= \bm{t}  &&\text{on $\Gamma_N$,}\\                                          
\end{aligned}\right.
\end{equation}
where the couple $(\bu,p)$ represents the velocity and pressure field, $\nu  > 0$ is the viscosity coefficient and $\bm{s}$, $\bu_D$ and $\bm{t}$ respectively are the volumetric source term, the Dirichlet boundary datum to impose the value of the velocity on $\Gamma_D$ and the pseudo-traction applied on the Neumann boundary $\Gamma_N$.

Assuming that $\Omega$ is partitioned in $\numel$ disjoint subdomains and splitting the second-order momentum conservation equation in two first-order equations, Equation~\eqref{eq:Stokes} can be written in the \emph{broken} computational domain as
\begin{equation} \label{eq:StokesBrokenMixed}
\left\{\begin{aligned}
\bL + \sqrt{\nu} \grad\bu & = \bm{0} &&\text{in $\Omega_e$, and for $e=1,\dotsc ,\numel$,}\\
\grad\cdot\bigl(\sqrt{\nu} \bL + p\Insd \bigr) &= \bm{s}         &&\text{in $\Omega_e$, and for $e=1,\dotsc ,\numel$,}\\
\grad\cdot\bu &= 0  &&\text{in $\Omega_e$, and for $e=1,\dotsc ,\numel$,}\\
\bu &= \bu_D     &&\text{on $\Gamma_D$,}\\
\bn \cdot \bigl(\sqrt{\nu} \bL + p \Insd\bigr)  &= - \bm{t}  &&\text{on $\Gamma_N$,}\\
\jump{\bu \otimes \bn} &=\bm{0}  &&\text{on $\Gamma$,}\\
\jump{\bn \cdot \bigl(\sqrt{\nu} \bL + p \Insd\bigr)} &= \bm{0}  &&\text{on $\Gamma$.}\\
\end{aligned} \right.
\end{equation}
where, as for the Poisson equation, the last two equations enforce the continuity of respectively the primal variable and the normal trace of the flux across the interface $\Gamma$.

\subsection{Strong form of the local and global problems}
\label{sc:StokesStrong}

The HDG formulation for the Stokes equation has been developed in a series of publications by Cockburn and co-workers~\cite{MR2485446,Nguyen-NPC:10,Nguyen-CNP:10,Nguyen-NPC:11}.
As previously discussed, the hybridisable discontinuous Galerkin method relies on writing Equation~\eqref{eq:StokesBrokenMixed} as a set of $\numel$ local problems defined element-by-element and featuring purely Dirichlet boundary conditions and a global problem to compute the hybrid variable defined on the mesh skeleton. First, for $e=1,\dotsc,\numel$ a solution $(\bL_e,\bu_e,p_e)$ is sought as a function of the unknown hybrid variable $\bhu$, namely
\begin{equation} \label{eq:StokesLocalStrong}
\left\{\begin{aligned}
\bL_e + \sqrt{\nu} \grad\bu_e & = \bm{0} &&\text{in $\Omega_e$,}\\
\grad\cdot\bigl(\sqrt{\nu} \bL_e + p_e\Insd \bigr) &= \bm{s}         &&\text{in $\Omega_e$,}\\
\grad\cdot\bu_e &= 0  &&\text{in $\Omega_e$}\\
\bu_e &= \bu_D     &&\text{on $\partial\Omega_e \cap \Gamma_D$,}\\
\bu_e &= \bhu     &&\text{on $\partial\Omega_e \setminus \Gamma_D$.}\\                       
\end{aligned} \right.
\end{equation}
It is worth noting that Equation~\eqref{eq:StokesLocalStrong} features a problem with only Dirichlet boundary conditions, hence the pressure is determined up to a constant. An additional constraint (e.g.\ setting the mean value of the pressure on the element boundary) is added to avoid the indeterminacy, namely
\begin{equation}\label{eq:pressureConstraint}
\frac{1}{|\partial\Omega_e|} \langle p_e, 1 \rangle_{\partial\Omega_e} = \rho_e,
\end{equation}
where $\rho_e$ denotes the mean pressure on the boundary of element $\Omega_e$.

In addition, the free divergence condition in Equation~\eqref{eq:StokesLocalStrong} induces the compatibility condition
\begin{equation}\label{eq:divergenceFreeConstraint}
\langle \bhu \cdot \bn_e , 1 \rangle_{\partial \Omega_e \setminus \Gamma_D} + \langle \bu_D \cdot \bn_e, 1 \rangle_{\partial \Omega_e\cap \Gamma_D} = 0.
\end{equation}

A global problem is defined to determine the trace of the velocity on the mesh skeleton, $\bhu$, and the mean pressure in each element, $\rho_e$, namely
\begin{equation} \label{eq:StokesGlobalStrong}
\left\{\begin{aligned}
\jump{\bu \otimes \bn} &=\bm{0}  &&\text{on $\Gamma$,}\\
\jump{\bn \cdot \bigl(\sqrt{\nu} \bL + p \Insd\bigr)} &= \bm{0}  &&\text{on $\Gamma$.}\\
\bn \cdot \bigl(\sqrt{\nu} \bL + p \Insd\bigr)  &= - \bm{t}  &&\text{on $\Gamma_N$.}\\  		
\end{aligned} \right.
\end{equation}
As previously discussed, the first condition in Equation~\eqref{eq:StokesGlobalStrong} is automatically satisfied due to the unique definition of the hybrid variable $\bhu$ on each face and the Dirichlet boundary condition $\bu_e = \bhu$ imposed in the local problems.

\subsection{Weak form of the local and global problems}
\label{sc:StokesWeak}

In addition to the internal products introduced in Equations~\eqref{eq:innerScalar} and \eqref{eq:innerVector}, the following internal products are defined
\begin{equation} \label{eq:innerTensorVector}
	(\bm{P},\bm{Q})_{\Omega_e} := \int_{\Omega_e} \bm{P} : \bm{Q} \ d\Omega , \quad \langle \hat{\bm{p}}, \hat{\bm{q}} \rangle_{\partial\Omega_e} := \sum_{\Gamma_i \subset \partial\Omega_e} \int_{\Gamma_i} \hat{\bm{p}} \cdot \hat{\bm{q}} \ d\Gamma 
\end{equation}
for tensor valued functions in $[\eltwo(\Omega)]^{\nsd \times \nsd}$ and vector valued functions in $[\eltwo(\Gamma_i)]^{\nsd}$ respectively.

The discrete weak formulation of the local problems reads as follows: for $e=1,\dotsc \numel$, find $(\bL^h_e,\bu^h_e,p^h_e)\in[\Vh(\Omega_e)]^{\nsd \times \nsd}\times[\Vh(\Omega_e)]^{\nsd}\times\Vh(\Omega_e)$ such that 
\begin{subequations}\label{HDG-Stokes-Dlocal}
	\begin{align}
	&
	- (\bm{G},\bL^h_e)_{\Omega_e} + (\grad\cdot\bm{G}, \sqrt{\nu} \bu^h_e)_{\Omega_e} =   \langle \bn_e \cdot \bm{G}  , \sqrt{\nu}\bu_D\rangle_{\partial\Omega_e\cap\Gamma_D} + \langle \bn_e \cdot\bm{G} , \sqrt{\nu}\bhu^h \rangle_{\partial\Omega_e\setminus\Gamma_D} , \label{eq:HDG-Stokes-Dlocal-L}
	\\
	&
	- (\grad\bw, \sqrt{\nu} \bL^h_e)_{\Omega_e} - (\grad\cdot\bw,p^h_e)_{\Omega_e} +\langle \bw, \bn_e \cdot \bigl(\widehat{\sqrt{\nu} \bL^h_e + p^h_e \Insd} \bigr)  \rangle_{\partial\Omega_e}  =  (\bw,\bm{s})_{\Omega_e}  , \label{eq:HDG-Stokes-Dlocal-Mom}
	\\
	&
	(\grad q, \bu^h_e)_{\Omega_e} = \langle q, \bu_D \cdot \bn_e \rangle_{\partial\Omega_e\cap\Gamma_D} + \langle q, \bhu^h \cdot \bn_e \rangle_{\partial\Omega_e\setminus\Gamma_D} , \label{eq:HDG-Stokes-Dlocal-U} 
	\\
	&
	\frac{1}{|\partial\Omega_e|}\langle p^h_e, 1 \rangle_{\partial\Omega_e} = \rho_e , \label{eq:HDG-Stokes-Dlocal-P}
	\end{align}
\end{subequations}
for all $(\bm{G},\bw,q)\in[\Vh(\Omega_e)]^{\nsd \times \nsd}\times[\Vh(\Omega_e)]^{\nsd}\times\Vh(\Omega_e)$. 

Integrating by parts Equation~\eqref{eq:HDG-Stokes-Dlocal-Mom} and introducing the definition of the trace of the numerical normal flux
\begin{equation} \label{eq:traceStokes}
\bn_e \cdot \bigl(\widehat{\sqrt{\nu} \bL^h_e + p^h_e \Insd} \bigr) := 
\begin{cases}
\bn_e \cdot \bigl( \sqrt{\nu} \bL^h_e + p^h_e \Insd \bigr) + \tau_e (\bu^h_e - \bu_D) & \text{on $\partial\Omega_e\cap\Gamma_D$,} \\
\bn_e \cdot \bigl( \sqrt{\nu} \bL^h_e + p^h_e \Insd \bigr) + \tau_e (\bu^h_e - \bhu^h) & \text{elsewhere.}  
\end{cases}
\end{equation}
leads to the following local problem:
\begin{subequations}\label{HDG-Stokes-DlocalSymm}
  \begin{gather}
	- (\bm{G},\bL^h_e)_{\Omega_e} + (\grad\cdot\bm{G}, \sqrt{\nu} \bu^h_e)_{\Omega_e} =  \langle \bn_e \cdot \bm{G}, \sqrt{\nu} \bu_D\rangle_{\partial\Omega_e\cap\Gamma_D} + \langle \bn_e \cdot \bm{G}, \sqrt{\nu} \bhu^h \rangle_{\partial\Omega_e\setminus\Gamma_D} , \label{eq:HDG-Stokes-Dlocal-LSymm} 
	\\
	\begin{split}
	(\bw, \sqrt{\nu} \grad \cdot \bL^h_e)_{\Omega_e} +  \langle \bw, \tau_e \bu^h_e  \rangle_{\partial\Omega_e} &+  (\bw, \grad p^h_e)_{\Omega_e} \\ 
	&=  (\bw,\bm{s})_{\Omega_e}  + \langle \bw, \tau_e \bu_D \rangle_{\partial\Omega_e \cap \Gamma_D} + \langle \bw, \tau_e \bhu^h \rangle_{\partial\Omega_e \setminus \Gamma_D} , \label{eq:HDG-Stokes-Dlocal-MomSymm}
	\end{split}
	\\
	(\grad q, \bu^h_e)_{\Omega_e} = \langle q, \bu_D \cdot \bn_e \rangle_{\partial\Omega_e\cap\Gamma_D} + \langle q, \bhu^h \cdot \bn_e \rangle_{\partial\Omega_e\setminus\Gamma_D} , \label{eq:HDG-Stokes-Dlocal-USymm} 
	\\
	\frac{1}{|\partial\Omega_e|} \langle  p^h_e, 1 \rangle_{\partial\Omega_e} = \rho_e , \label{eq:HDG-Stokes-Dlocal-PSymm}
 \end{gather}
\end{subequations}

In a similar fashion, the following global problem accounting for the transmission conditions and the Neumann boundary condition given in Equation~\eqref{eq:StokesGlobalStrong} is: find $\bhu^h\in[\VhHat(\Gamma\cup\Gamma_N)]^{\nsd}$ and $\rho \in \mathbb{R}^{\numel}$ 
\begin{subequations}\label{HDG-Stokes-Dglobal}
  \begin{gather}
  \begin{split}
	\sum_{e=1}^{\numel}\Bigl\{
	\langle \bhw, \bn_e \cdot \sqrt{\nu} \bL^h_e &\rangle_{\partial\Omega_e\setminus\Gamma_D}
	+ \langle \bhw,p^h_e \bn_e \rangle_{\partial\Omega_e\setminus\Gamma_D}  \\
	&+ \langle \bhw,\tau_e\, \bu^h_e \rangle_{\partial\Omega_e\setminus\Gamma_D} 
	 - \langle \bhw,\tau_e\,\bhu^h \rangle_{\partial\Omega_e\setminus\Gamma_D}\Bigr\}
	= -\sum_{e=1}^{\numel} \langle \bhw, \bm{t} \rangle_{\partial\Omega_e\cap\Gamma_N},
  \end{split}
  \\
	\langle \bhu^h \cdot \bn_e, 1 \rangle_{\partial \Omega_e \setminus \Gamma_D} = - \langle \bu_D \cdot \bn_e , 1 \rangle_{\partial \Omega_e\cap \Gamma_D} \quad \text{ for } e=1,\dotsc,\numel,
	\label{eq:weakCompatibilityStokes}
  \end{gather}
\end{subequations}
for all $\bhw\in[\VhHat(\Gamma\cup\Gamma_N)]^{\nsd}$.

Henceforth, to simplify the notation, the superindex $^h$ expressing the discrete approximations will be dropped, unless needed in order to follow the development.

\subsection{FCFV discretisation}
\label{sc:StokesFV}

The discretisation of the local problems given by Equation~\eqref{HDG-Stokes-DlocalSymm} with a degree of approximation $k=0$ in each element for both $\bL_e$, $\bu_e$ and $p_e$ and also a degree of approximation $k=0$ in each face/edge for $\bhu$ leads to 
\begin{subequations}\label{HDG-Stokes-localK0}
	\begin{align}
	&
	-| \Omega_e | \Le =  \sqrt{\nu} \sum_{j \in \Dset} | \Gamma_{e,j} | \bn_j \otimes \bu_{D,j}   + \sqrt{\nu} \sum_{j \in \Bset} | \Gamma_{e,j} | \bn_j \otimes \buHj  , \label{eq:HDG-Stokes-Dlocal-LK0} 
	\\
	&
	\sum_{j \in \Aset} | \Gamma_{e,j} | \tau_j \bue =  | \Omega_e | \bm{s}_e  + \sum_{j \in \Dset} | \Gamma_{e,j} | \tau_j \bu_{D,j} + \sum_{j \in \Bset} | \Gamma_{e,j} | \tau_j \buHj , \label{eq:HDG-Stokes-Dlocal-MomK0} 
	\\
	&
	0 = \sum_{j \in \Dset} | \Gamma_{e,j} | \bu_{D,j} \cdot \bn_j + \sum_{j \in \Bset} | \Gamma_{e,j} | \buHj \cdot \bn_j , \label{eq:HDG-Stokes-Dlocal-UK0} 
	\\
	&
	\pe = \rho_e, \label{eq:HDG-Stokes-Dlocal-PK0}
	\end{align}
\end{subequations}
for $e=1,\dotsc,\numel$.

It is important to note that Equation~\eqref{eq:HDG-Stokes-Dlocal-UK0} coincides with the discretised version of the global compatibility condition of Equation~\eqref{eq:weakCompatibilityStokes}, thus it may be neglected in the local computations and be imposed solely in the global problem. The three remaining equations are uncoupled and provide the following expressions of $\bL_e$, $\bu_e$ and $p_e$ as functions of the global unknowns $\bhu$ and $\rho_e$:
\begin{subequations}\label{HDG-Stokes-localVars}
	\begin{align}
	&
	\Le = - \sqrt{\nu} | \Omega_e |^{-1} \mat{Z}_e  - \sqrt{\nu} | \Omega_e |^{-1} \sum_{j \in \Bset} | \Gamma_{e,j} | \bn_j \otimes  \buHj,
	\label{eq:Stokes-locL}
	\\
	&
	\bue =  \alpha_e^{-1}\bm{\beta}_e + \alpha_e^{-1} \sum_{j \in \Bset} | \Gamma_{e,j} | \tau_j \buHj ,
	\label{eq:Stokes-locU}
	\\
	&
	\pe = \rho_e,
	\label{eq:Stokes-locP}
	\end{align}
\end{subequations}
where the following quantities only depend upon the data of the problem and may be precomputed:
\begin{equation} \label{eq:stokesPrecomp}
\alpha_e  := \sum_{j \in \Aset} | \Gamma_{e,j} | \tau_j , 
\qquad
\bm{\beta}_e  :=  | \Omega_e | \bm{s}_e  + \sum_{j \in \Dset} | \Gamma_{e,j} | \tau_j \bu_{D,j}, 
\qquad
\mat{Z}_e  := \sum_{j \in \Dset} | \Gamma_{e,j} |\bn_j \otimes  \bu_{D,j}.
\end{equation}

In a similar fashion, the global problem of Equation~\eqref{HDG-Stokes-Dglobal} particularised for a constant degree of approximation $k=0$ leads to
\begin{subequations}\label{HDG-Stokes-globalK0}
  \begin{gather}
  \begin{split} 
	\sum_{e=1}^{\numel}\Bigl\{
	\sqrt{\nu} | \Gamma_{e,i} | \bn_i \cdot \Le  + | \Gamma_{e,i} | \pe \bn_i &+ | \Gamma_{e,i} | \tau_i \bue -  | \Gamma_{e,i} | \tau_i \buHi  \Bigr\} \\
	&= -\sum_{e=1}^{\numel} | \Gamma_{e,i} | \bm{t}_i \, \chi_{\Nset}(i) \quad \text{ for } i \in \Bset , \label{HDG-Stokes-globalK0-U}
  \end{split}
  \\
	\sum_{j \in \Bset} | \Gamma_{e,j} | \buHj \cdot \bn_j = - \sum_{j \in \Dset} | \Gamma_{e,j} | \bu_{D,j} \cdot \bn_j \quad \text{ for } e=1,\dotsc,\numel . \label{HDG-Stokes-globalK0-Rho}
  \end{gather}
\end{subequations}

By plugging Equation~\eqref{eq:Stokes-locL}, \eqref{eq:Stokes-locU} and \eqref{eq:Stokes-locP} into Equation~\eqref{HDG-Stokes-globalK0}, the global problem can be written in terms of the global unknowns $\bhu$ and $\rho$. That is, being $\vect{\hat{u}}$ the vector containing the value of the hybrid variable on the faces on $\Gamma \cup \Gamma_N$ and $\bm{\rho}$ the vector containing the values of the mean pressure on each element $\Omega_e$, the following linear system is obtained:
\begin{equation} \label{eq:globalSystemStokes}
\begin{bmatrix}
\mat{\widehat{K}}_{\hu \hu} & \mat{\widehat{K}}_{\hu \rho} \\
\mat{\widehat{K}}_{\hu \rho}^T & \mat{0}_{\numel}
\end{bmatrix}
\begin{Bmatrix}
\vect{\hat{u}}  \\
\bm{\rho}
\end{Bmatrix}
=
\begin{Bmatrix}
\vect{\hat{f}}_{\hu} \\
\vect{\hat{f}}_{\rho}
\end{Bmatrix} ,
\end{equation}
where the blocks composing the matrices and the vectors of the previous linear system  are computed by assembling the contributions given by
\begin{subequations}\label{HDG-Stokes-globalSystem}
	\begin{align}
	&
	(\mat{\widehat{K}}_{\hu \hu})^e_{i,j} :=  |\Gamma_{e,i}| \left( \alpha_e^{-1} \tau_i \tau_j |\Gamma_{e,j}| - \nu | \Omega_e |^{-1} |\Gamma_{e,j}| \bn_i \cdot \bn_j  - \tau_i \delta_{ij} \right) \mat{I}_{\nsd}
	, \\
	&
	(\mat{\widehat{K}}_{\hu \rho})^e_i := |\Gamma_{e,i}| \bn_i , \\
	&
	(\vect{\hat{f}}_{\hu})^e_i :=  |\Gamma_{e,i}| \left( \nu | \Omega_e |^{-1} \bn_i \cdot \mat{Z}_e  - \bm{t}_i \, \chi_{\Nset}(i)  - \alpha_e^{-1} \tau_i \bm{\beta}_e  \right) , \\
	&
	(\hat{\text{f}}_{\rho})^e := -\sum_{j \in \Dset} | \Gamma_{e,j} | \bu_{D,j} \cdot \bn_j,
	\end{align}
\end{subequations}
for $i,j \in \Bset$.
It is important to emphasise that $(\mat{\widehat{K}}_{\hu \hu})^e_{i,j}$ denotes a matrix,  $(\mat{\widehat{K}}_{\hu \rho})^e_i$ and $(\vect{\hat{f}}_{\hu})^e_i$ are vectors and $(\hat{\text{f}}_{\rho})^e$ is a scalar.

\begin{remark}
When Dirichlet boundary conditions are imposed in the whole boundary (i.e. $\Gamma_D=\partial \Omega$ and $\Gamma_N = \emptyset$) an additional constraint  must be imposed to avoid
the indeterminacy of the pressure. It is common~\cite{MR2772094,Nguyen-NPC:10} to impose zero mean pressure on the domain, namely
\begin{equation}
\sum_{e=1}^{\numel} |\Omega_e| \rho_e = 0.
\end{equation}
In this case, the global system of Equation~\eqref{HDG-Stokes-globalSystem} must be modified to account for the extra constraint on the pressure.
\end{remark}

\section{Computational aspects}
\label{sc:computational}

In a series of papers by Cockburn and co-workers~\cite{Jay-CGL:09,MR2772094,MR2904582,MR3044180,MR3340089,MR3626531}, the optimal rate of convergence of HDG has been proved and numerically verified for a wide class of problems. More precisely, for Poisson and Stokes equations using constant degree of approximation, both the primal variables ($u$ in the Poisson equation and velocity $\bu$ and pressure $p$ in Stokes equation) and the dual variables representing the fluxes ($\bm{q} = - \grad u$ in Poisson equation and $\bm{L} = - \sqrt{\nu} \grad\bu$ in Stokes equation) converge with first-order accuracy.
The FCFV method inherits the convergence properties of HDG, that is, it experiences first-order convergence for both the primal variables and their fluxes.
The result for the primal variables is the same as other finite volume techniques  (e.g.\ cell-centred and vertex-centred finite volumes). The most distinctive feature of the proposed FCFV method is that first-order convergence is also achieved for the gradient of the solution without any reconstruction.

For the solution of Stokes flow problems, the proposed FCFV method does not require the solution of a Poisson problem for computing the pressure, as required by segregated schemes such as the semi-implicit method for pressure-linked equations (SIMPLE) algorithm~\cite{patankar1980numerical}. In addition, contrary to other mixed finite element methods, with the FCFV it is possible to use the same space of approximation for both velocity and pressure, circumventing the so-called LBB condition~\cite{donea2003finite}.

A remarkable property of the proposed FCFV method, inherited from the HDG method, is the uncoupled nature of the variables appearing in the local problems. More precisely, an analytical closed form of the primal and dual variables, $(\bm{q}_e,u_e)$ and $(\bm{L}_e,\bu_e,p_e)$ for Poisson and Stokes respectively is given in terms of the global variables, $\hu$ and  $(\bhu,\rho)$ for Poisson and Stokes respectively.

A drawback of the FCFV method with respect to other finite volume strategies is represented by the higher number of degrees of freedom. This issue is due to the higher number of faces with respect to the elements and the vertices for a given cardinality of the mesh (cf. tables \ref{tab:meshCardinality2D} and \ref{tab:meshCardinality3D} for two and three dimensional meshes respectively).
Nevertheless, a detailed comparison of the computational costs induced by CCFV, VCFV and FCFV strategies to compute a solution for a given precision should be performed in order to state any final conclusion on the advantages and disadvantages of the method concerning its computational cost. Note for instance that CCFV and FCFV, in contrast to VCFV, have, even for unstructured meshes, a fixed connectivity, which has major influences in the computability costs.
\begin{table}[hbt]
	\centering
	\begin{tabular}[hbt]{| l || c | c | c |}
		\hline
		Type & Vertices & Cells & Edges \\
		\hline & & & 
		\\ [-1em] \hline
		Triangles & $n$ & $2n$ & $3n$ \\
		\hline
		Quadrilaterals & $n$ & $n$ & $2n$ \\
		\hline
	\end{tabular}
	\caption{Number of vertices, cells and edges for meshes in two dimensions.}
	\label{tab:meshCardinality2D}
\end{table}
\begin{table}[hbt]
	\centering
	\begin{tabular}[hbt]{| l || c | c |  c |}
		\hline
		Type & Vertices & Cells &  Faces \\
		\hline & & &
		\\ [-1em] \hline
		Tetrahedrons & $n$ & $5n$ &  $10n$ \\
		\hline
		Hexahedrons & $n$ & $n$ & $3n$ \\
		\hline
		Prisms & $n$ & $2n$  & $5n/3$ \\
		\hline
		Pyramids & $n$ & $8n/5$ & $4n$ \\
		\hline
	\end{tabular}
	\caption{Number of vertices, cells and faces for meshes in three dimensions.}
	\label{tab:meshCardinality3D}
\end{table}

Despite the increased number of degrees of freedom compared to the VCFV, the proposed FCFV method requires a significantly low number of operations to construct the global system of equations. For solving the Poisson problem, the computation of the elemental matrix and right hand side vector in Equation~\eqref{HDG-Poisson-globalSystem} only requires a total of $4\nsd+12$ operations. For the Stokes problem, the computation of the elemental matrices and right hand side vectors in Equation~\eqref{HDG-Stokes-globalSystem} only requires a total of $2\nsd^2 + (|\Dset| + 6)\nsd + 2|\Dset| + 9$ operations. The operation counts does not include the operations required to compute the terms in Equations~\eqref{eq:poissonPrecomp} and \eqref{eq:stokesPrecomp} as these scalars and vectors can usually be computed once and stored without incurring in a significant memory consumption. For instance, the memory required for storing the terms in Equation~\eqref{eq:poissonPrecomp} is equal to 8(\nsd+2) MB per million elements. Similarly, the memory required for storing the terms in Equation~\eqref{eq:stokesPrecomp} is equal to ($\nsd^2 + \nsd +1$) MB per million elements.

It is also worth noting that the expressions appearing in Equations~\eqref{HDG-Poisson-globalSystem} and \eqref{HDG-Stokes-globalSystem} 
may be further simplified under some assumptions on the nature of the mesh. That is, some of the terms involving the outward unit normals to the element faces vanish and result in simpler expressions exploiting the orthogonality properties of Cartesian grids.

\section{Numerical studies}
\label{sc:studies}

This Section presents a set of numerical studies to verify the optimal convergence properties of the proposed approach, to compare the performance for different element types and to study the influence of numerical parameters and mesh metrics on the accuracy of the approximation. Two and three dimensional examples are considered for both the Poisson and the Stokes equations.

\subsection{Optimal convergence of the FCFV scheme for Poisson equation}
\label{sc:PoissonVerification}

The model problem of Equation~\eqref{eq:Poisson} is solved in two dimensions to test the optimal convergence of the proposed FCFV solver. The computational domain is $\Omega=[0,1]^2$. The source term is selected so that the analytical solution is 
\begin{equation}
u(x,y) = \exp \big( \alpha \sin(ax+cy) + \beta\cos(bx+dy) \big),
\end{equation}
with $\alpha=0.1$, $\beta=0.3$, $a=5.1$, $b=4.3$, $c=-6.2$ and $d=3.4$. Neumann boundary conditions, corresponding to the analytical normal flux, are imposed in $\Gamma_N = \{(x,y) \in \mathbb{R}^2 \; | \; y=0\}$ and Dirichlet boundary conditions, corresponding to the analytical solution, are imposed in $\Gamma_D = \partial \Omega \setminus \Gamma_N$. 

\begin{figure}[!tb]
	\centering
	\subfigure[Mesh 1]{\includegraphics[width=0.24\textwidth]{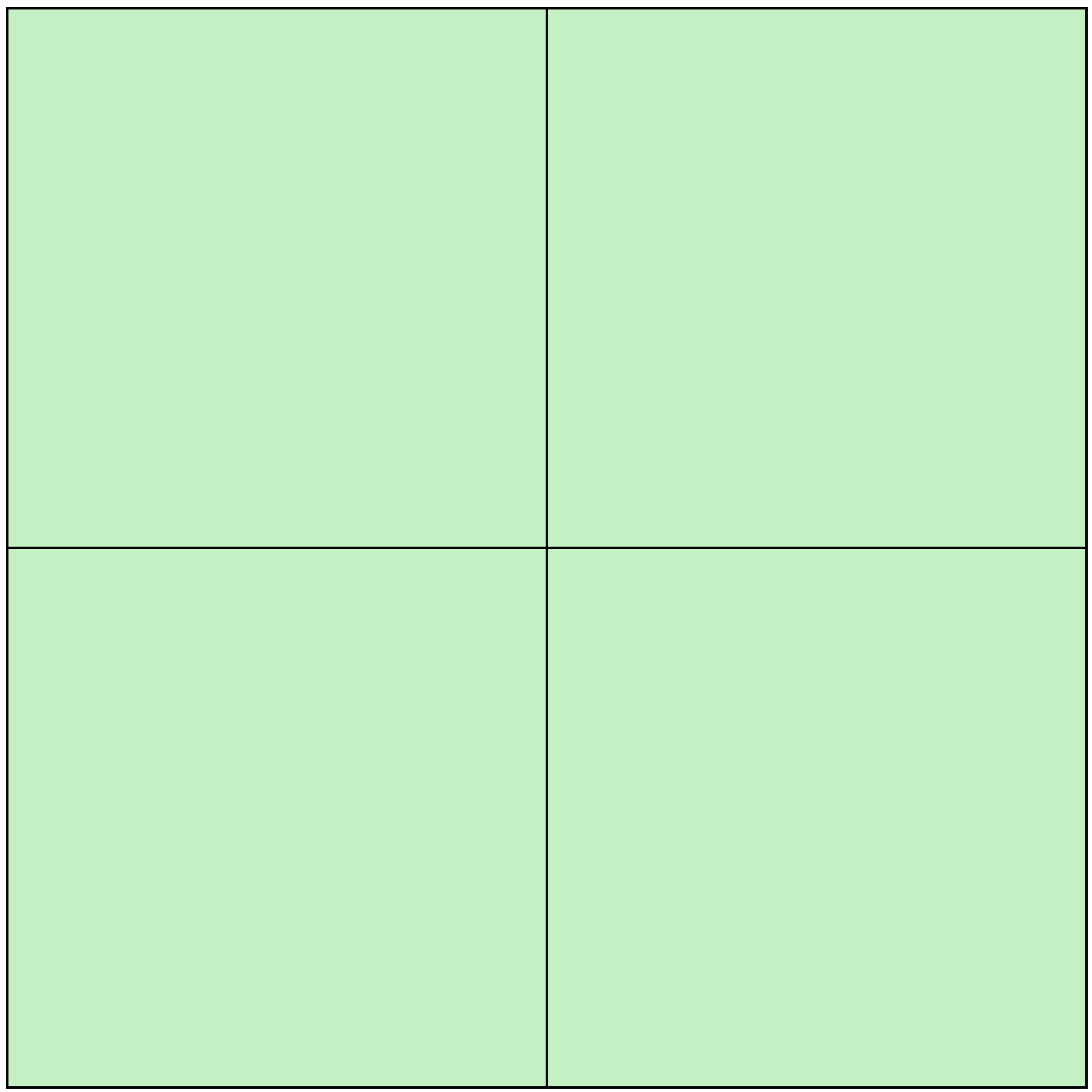}}
	\subfigure[Mesh 2]{\includegraphics[width=0.24\textwidth]{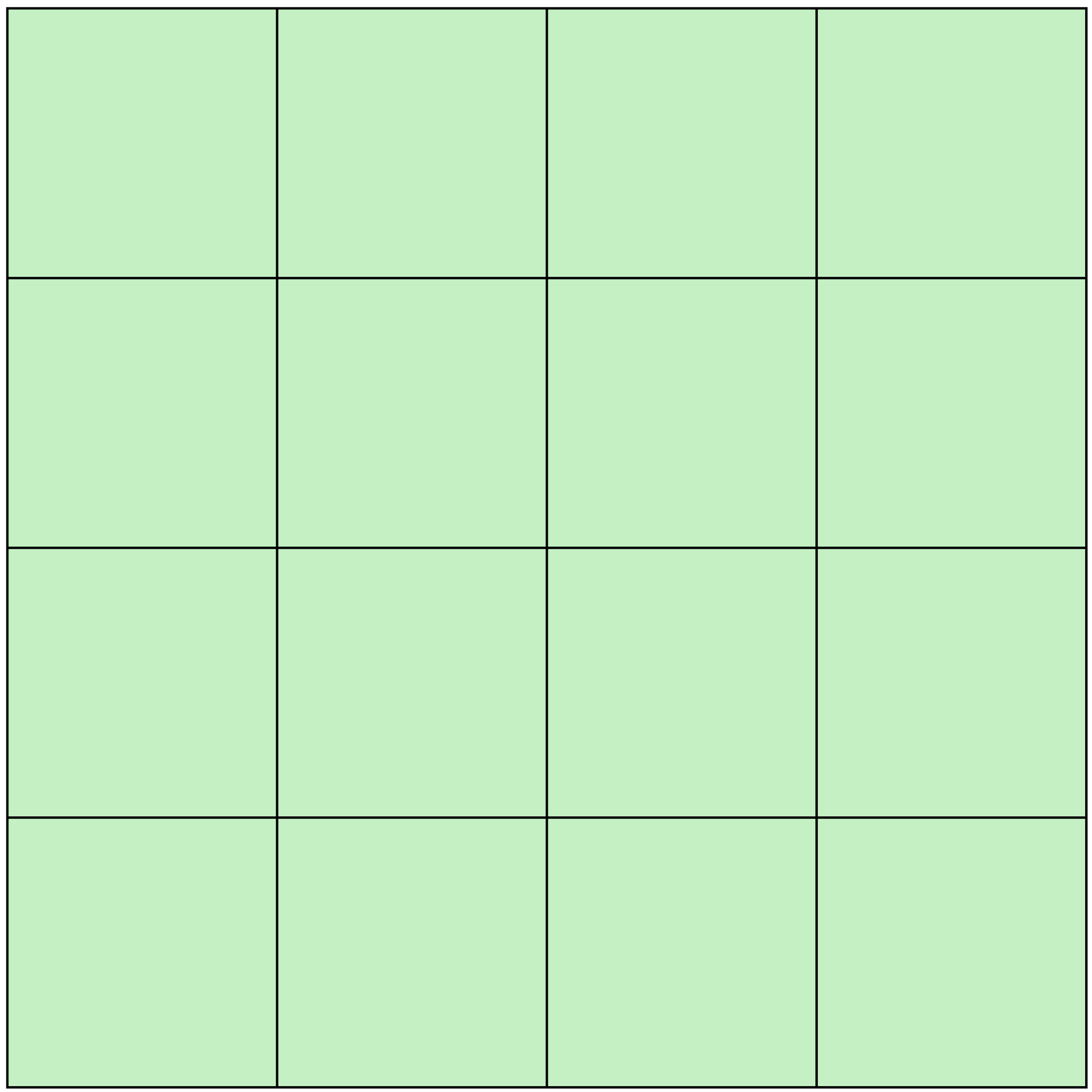}}
	\subfigure[Mesh 3]{\includegraphics[width=0.24\textwidth]{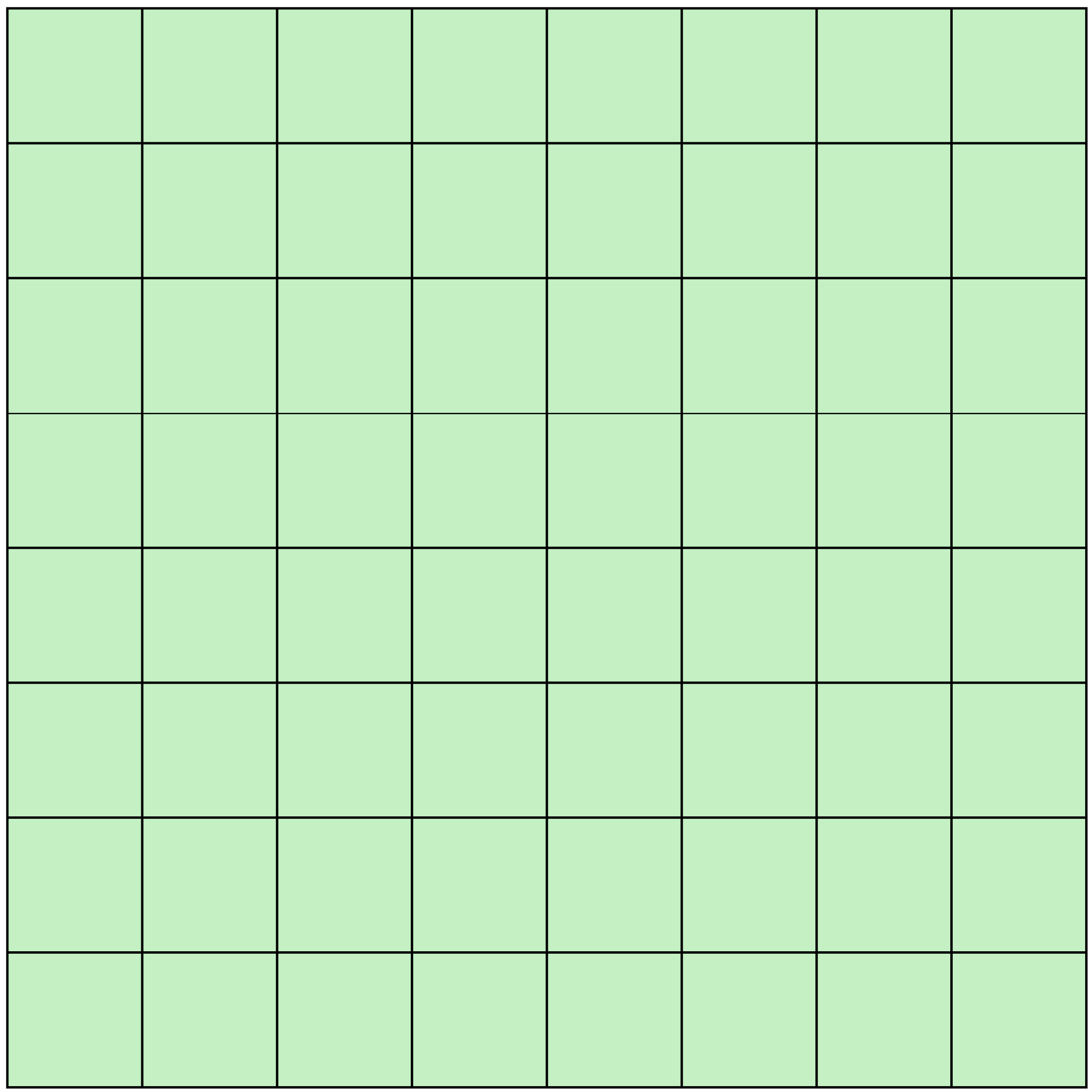}}
	\subfigure[Mesh 4]{\includegraphics[width=0.24\textwidth]{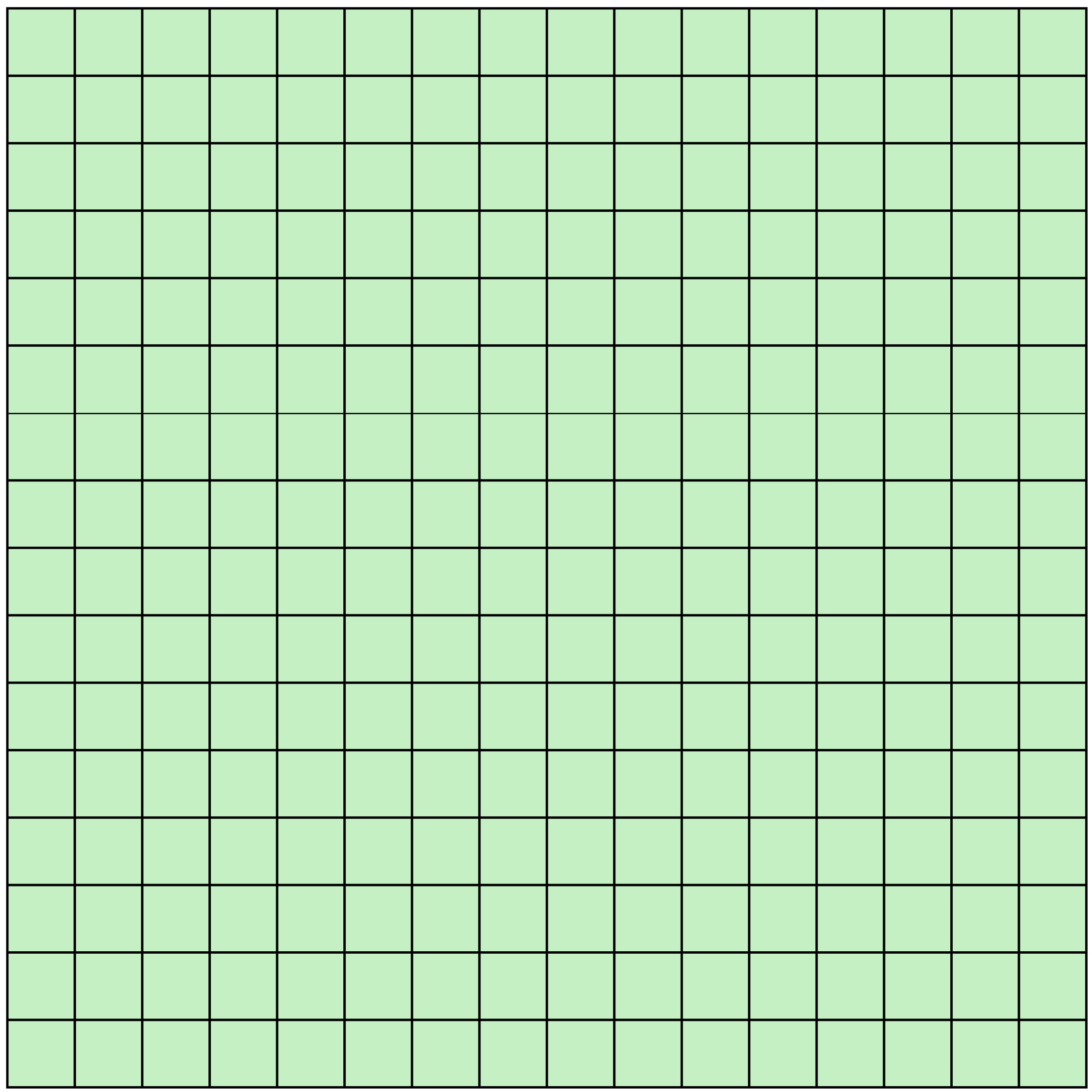}}
	\caption{Quadrilateral meshes of $\Omega=[0,1]^2$ for the mesh convergence analysis in 2D.}
	\label{fig:poisson2D_QUAmeshes}
\end{figure}

\begin{figure}[!tb]
	\centering
	\subfigure[Mesh 1]{\includegraphics[width=0.24\textwidth]{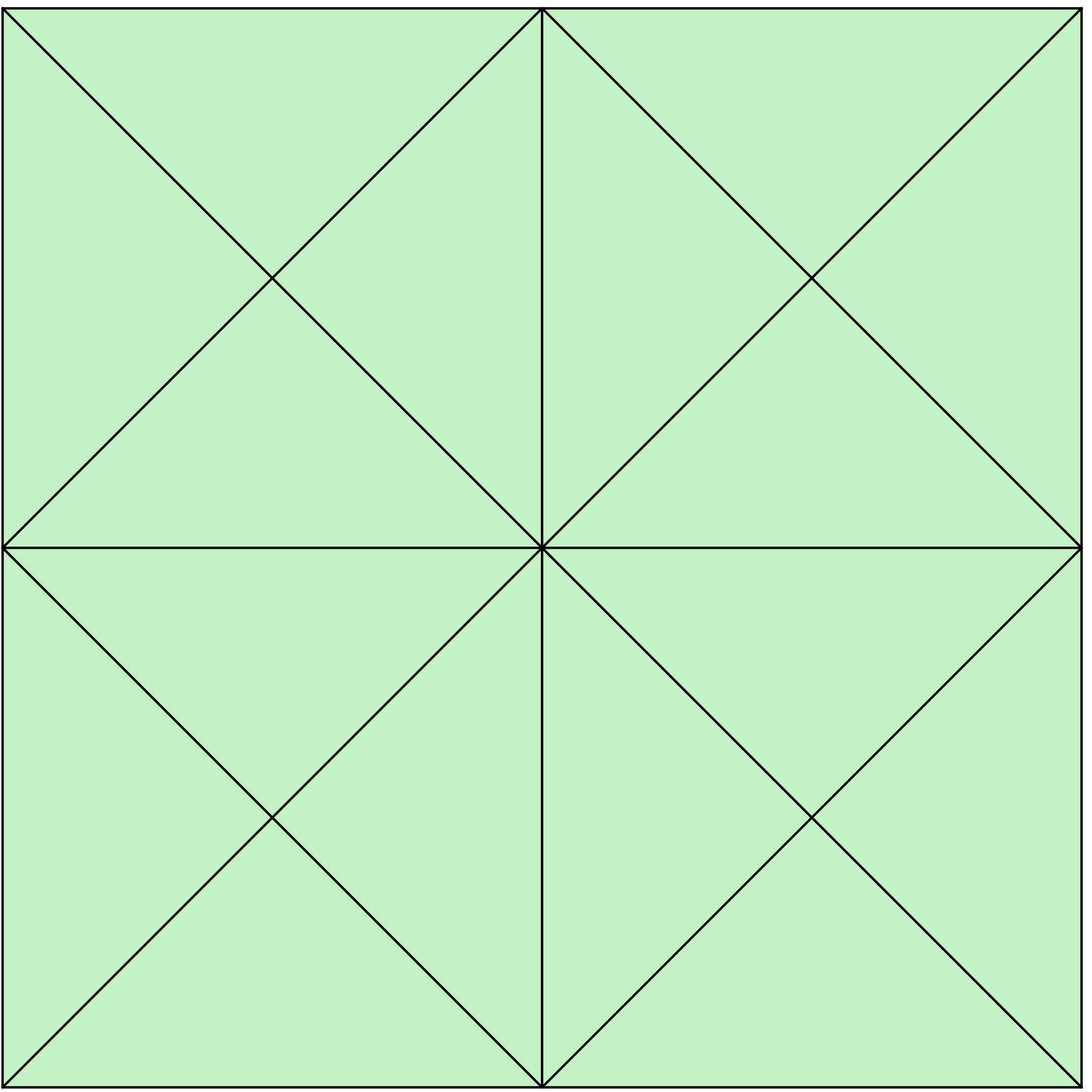}}
	\subfigure[Mesh 2]{\includegraphics[width=0.24\textwidth]{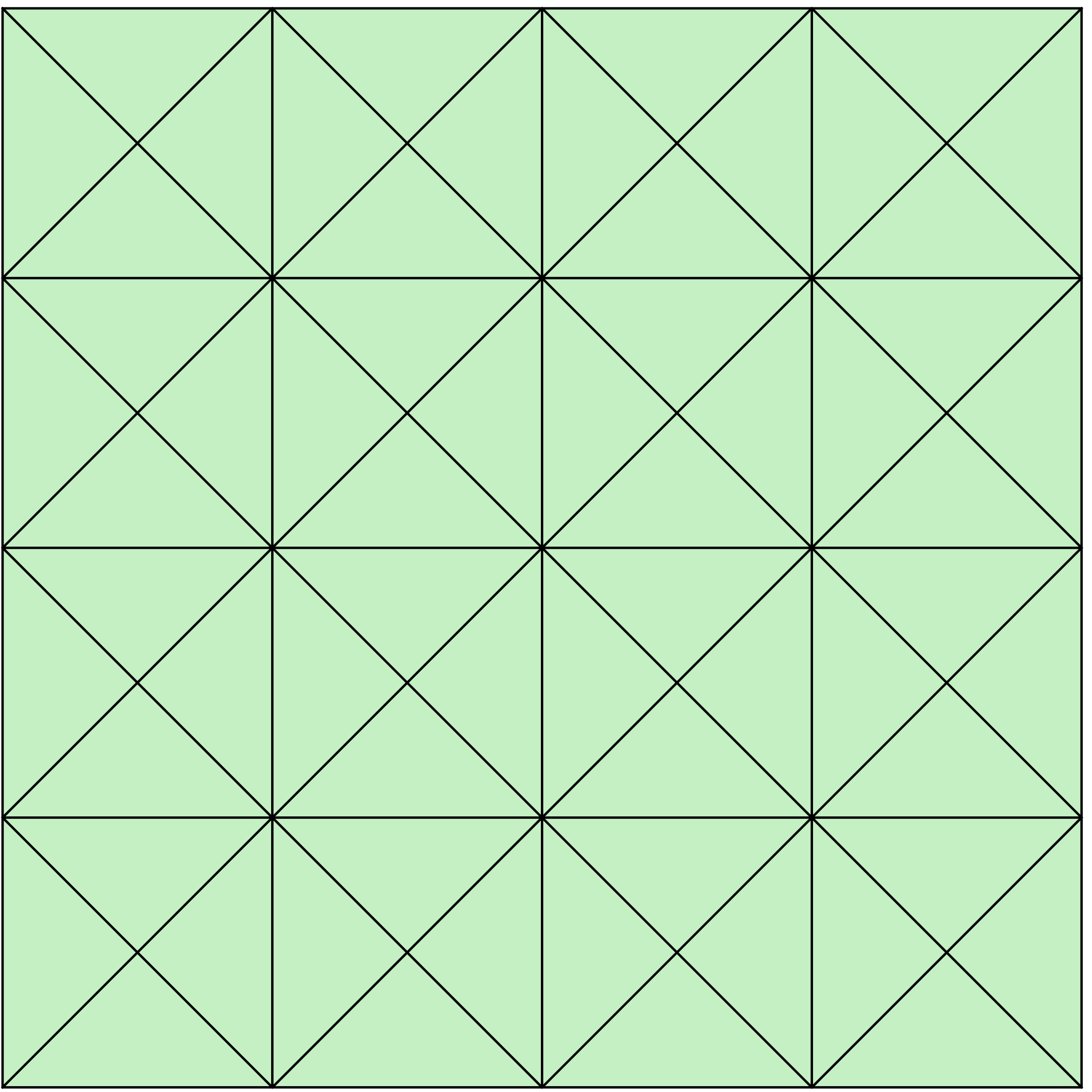}}
	\subfigure[Mesh 3]{\includegraphics[width=0.24\textwidth]{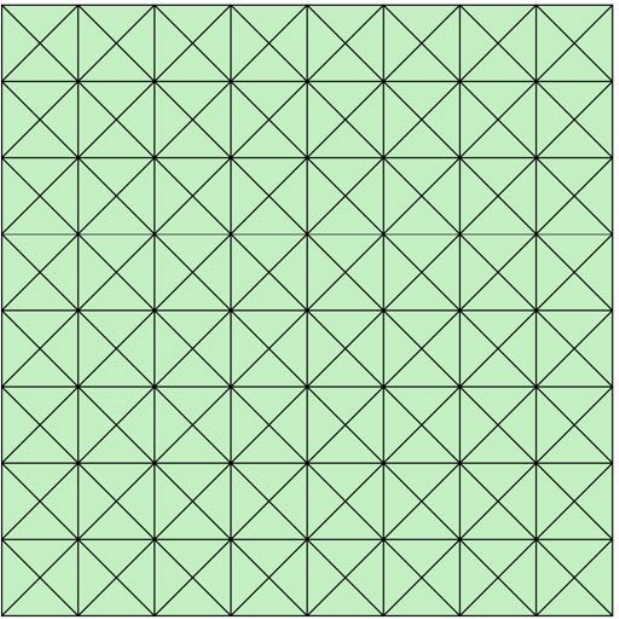}}
	\subfigure[Mesh 4]{\includegraphics[width=0.24\textwidth]{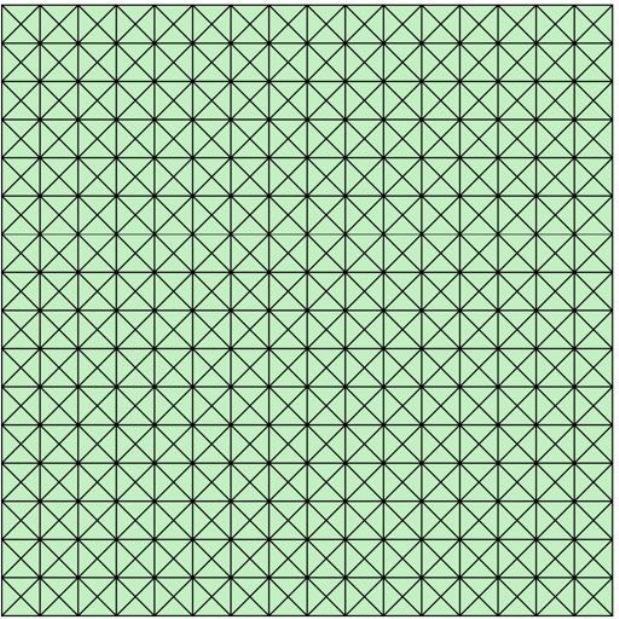}}
	\caption{Triangular meshes of $\Omega=[0,1]^2$ for the mesh convergence analysis in 2D.}
	\label{fig:poisson2D_TRImeshes}
\end{figure}

Quadrilateral and triangular uniform meshes are considered to perform an $h$-convergence study. The first four quadrilateral and triangular meshes are shown in Figures~\ref{fig:poisson2D_QUAmeshes} and $\ref{fig:poisson2D_TRImeshes}$ respectively. 

\begin{figure}[!tb]
	\centering
	\subfigure[Mesh 3]{\includegraphics[width=0.24\textwidth]{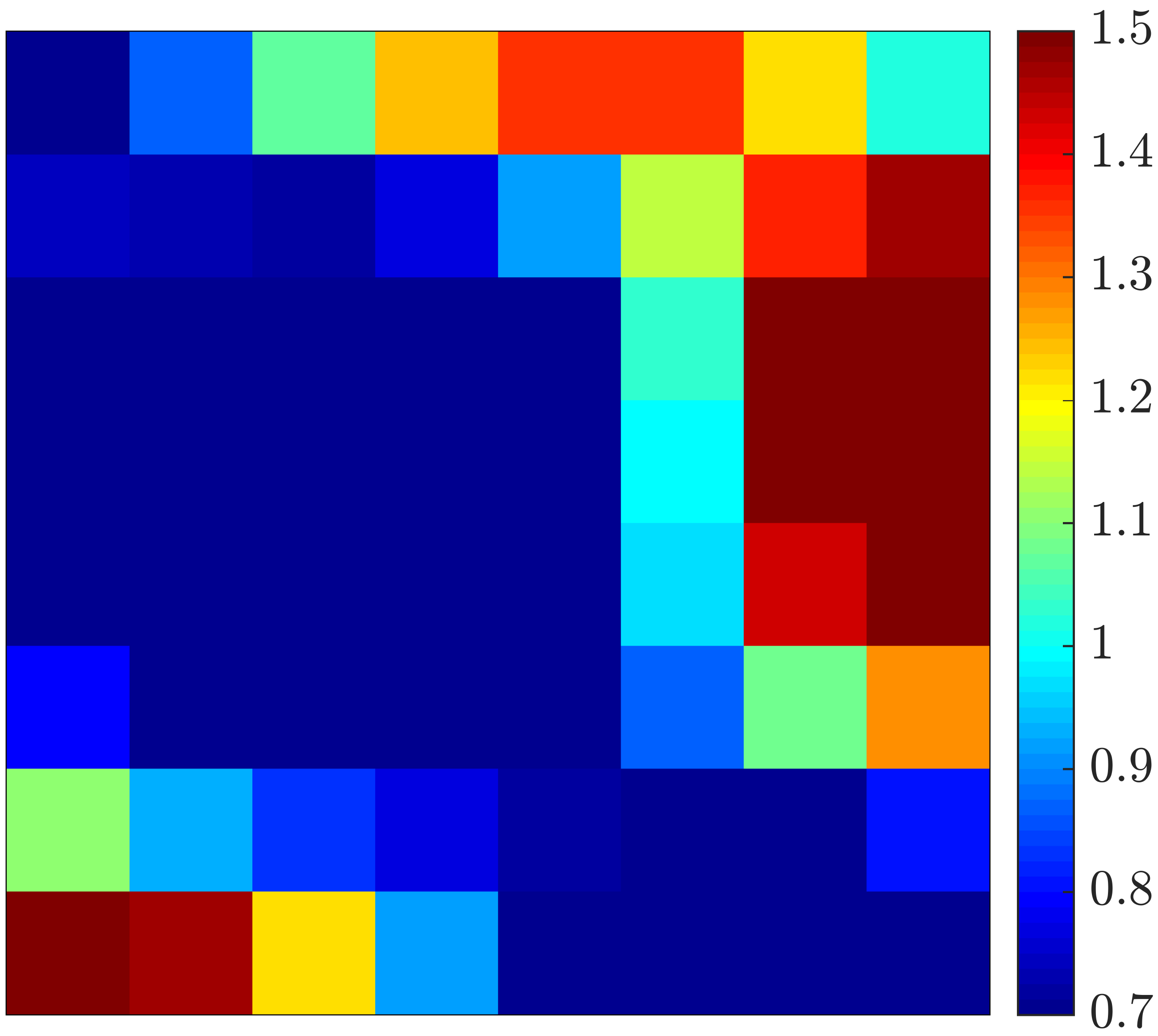}}
	\subfigure[Mesh 5]{\includegraphics[width=0.24\textwidth]{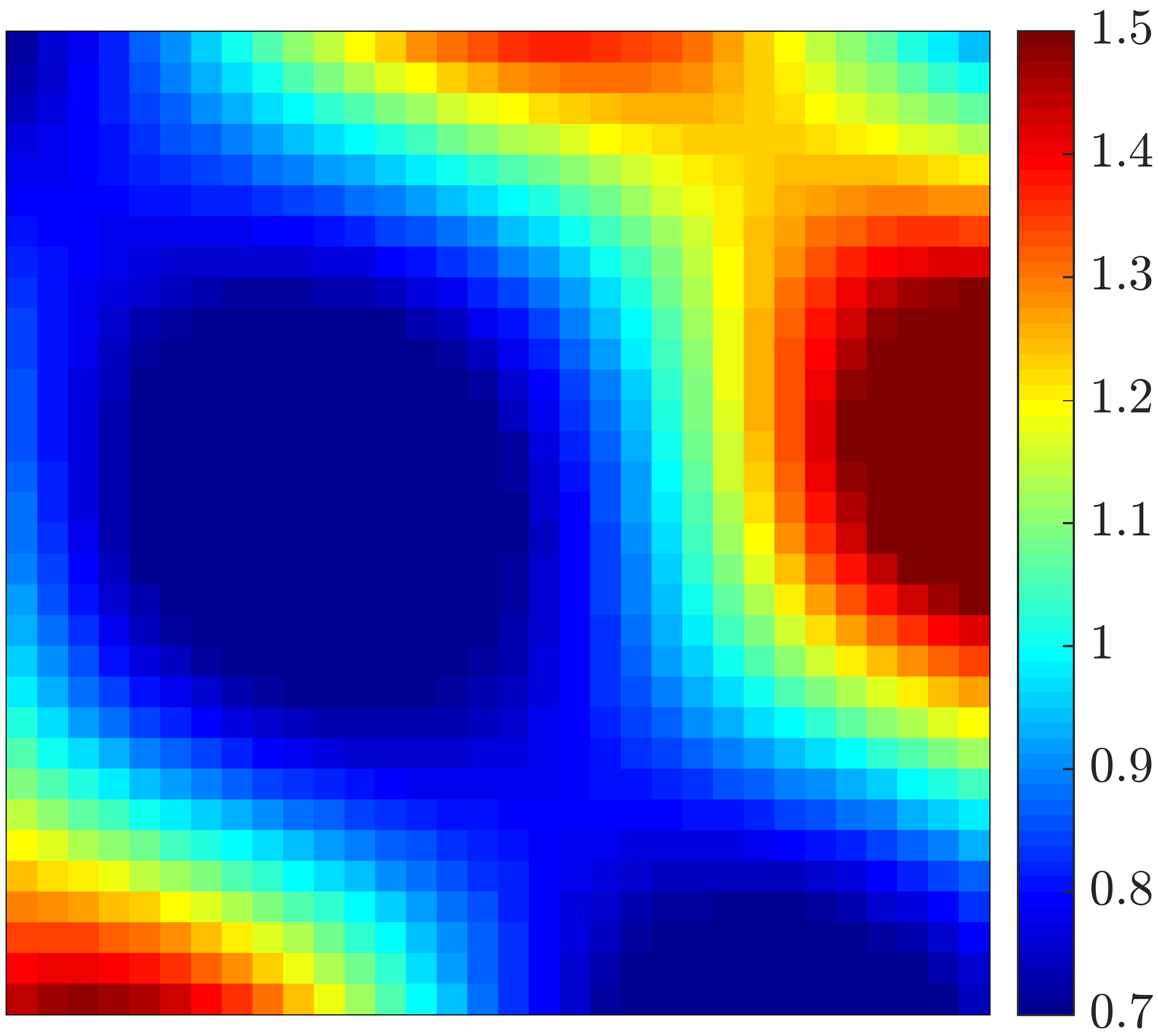}}
	\subfigure[Mesh 7]{\includegraphics[width=0.24\textwidth]{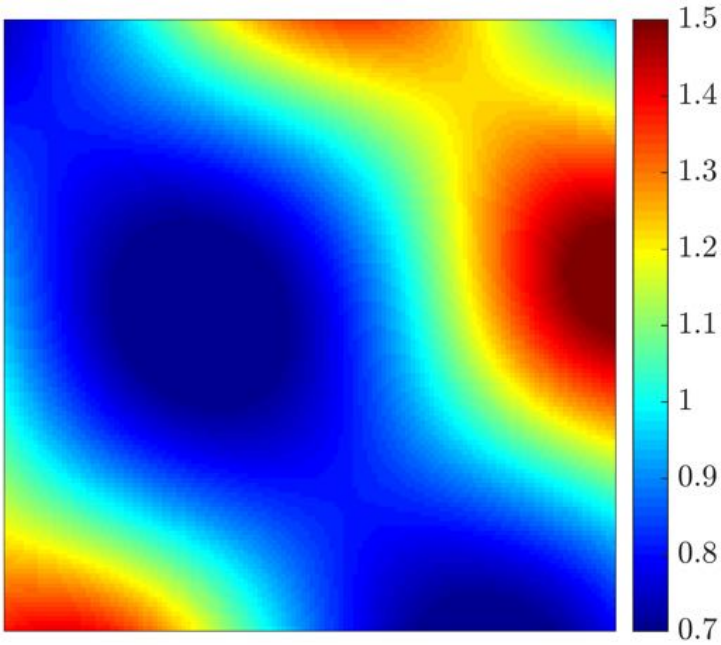}}
	\subfigure[Mesh 9]{\includegraphics[width=0.24\textwidth]{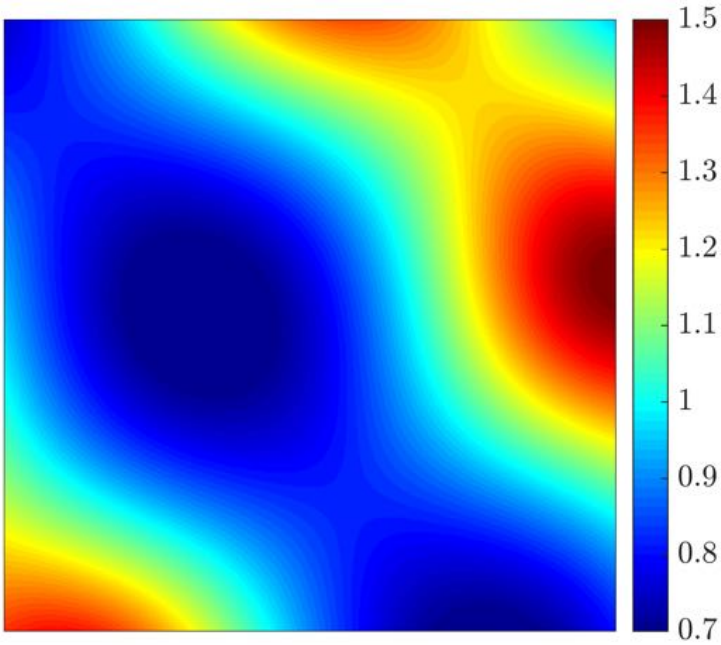}}
	\caption{Solution of the 2D Poisson problem using quadrilateral meshes.}
	\label{fig:poisson2D_QUAsolutions}
\end{figure}

\begin{figure}[!tb]
	\centering
	\subfigure[Mesh 3]{\includegraphics[width=0.24\textwidth]{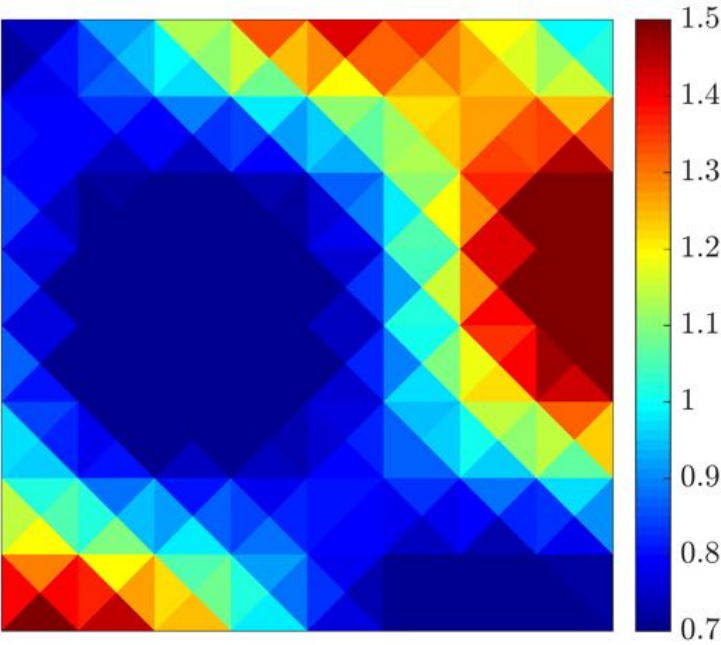}}
	\subfigure[Mesh 5]{\includegraphics[width=0.24\textwidth]{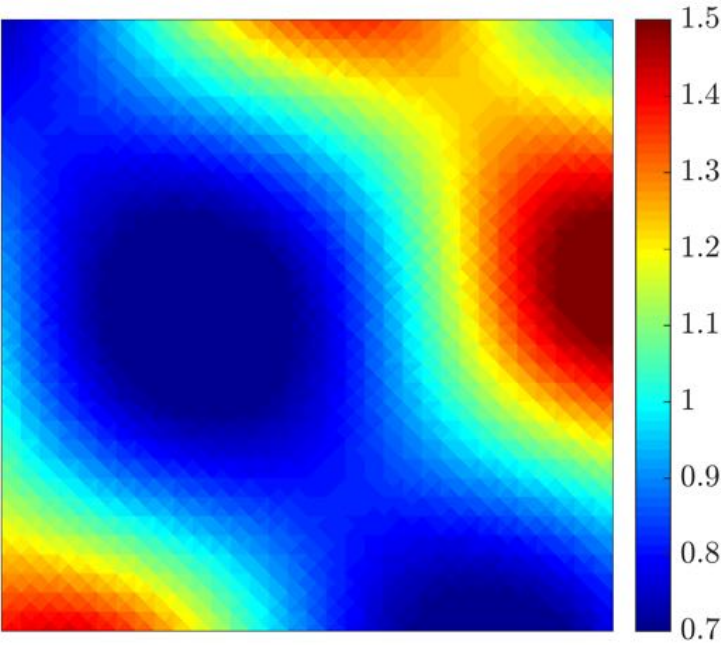}}
	\subfigure[Mesh 7]{\includegraphics[width=0.24\textwidth]{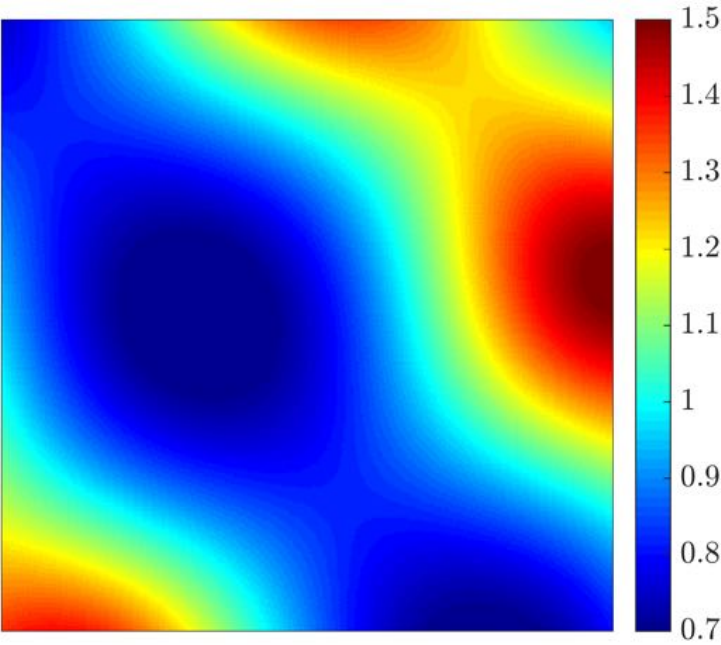}}
	\subfigure[Mesh 9]{\includegraphics[width=0.24\textwidth]{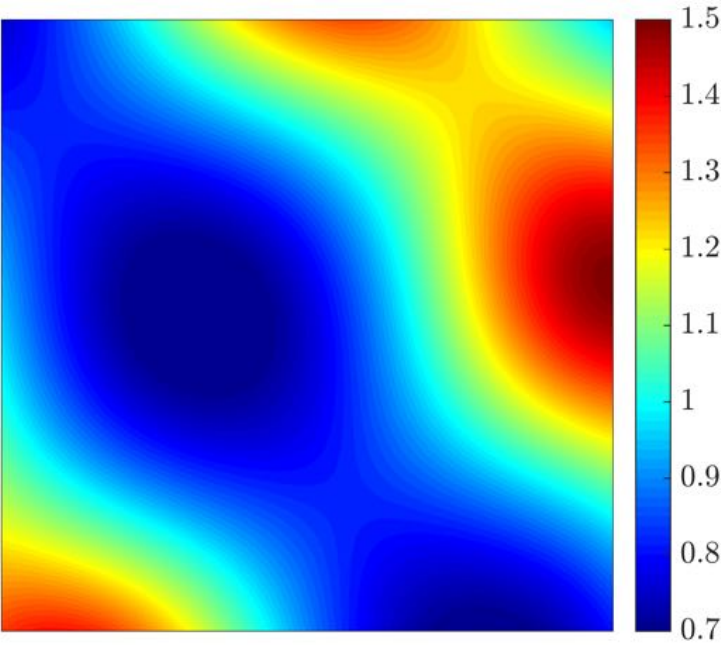}}
	\caption{Solution of the 2D Poisson problem using triangular meshes.}
	\label{fig:poisson2D_TRIsolutions}
\end{figure}

The numerical solution obtained with the proposed FCFV scheme for selected quadrilateral and triangular meshes is depicted in Figures~\ref{fig:poisson2D_QUAsolutions} and \ref{fig:poisson2D_TRIsolutions} respectively. The results clearly illustrate the constant degree of approximation used within each element and the increased accuracy obtained as the mesh is refined. 

\begin{figure}[!tb]
	\centering
	\subfigure[$u$]{\includegraphics[width=0.4\textwidth]{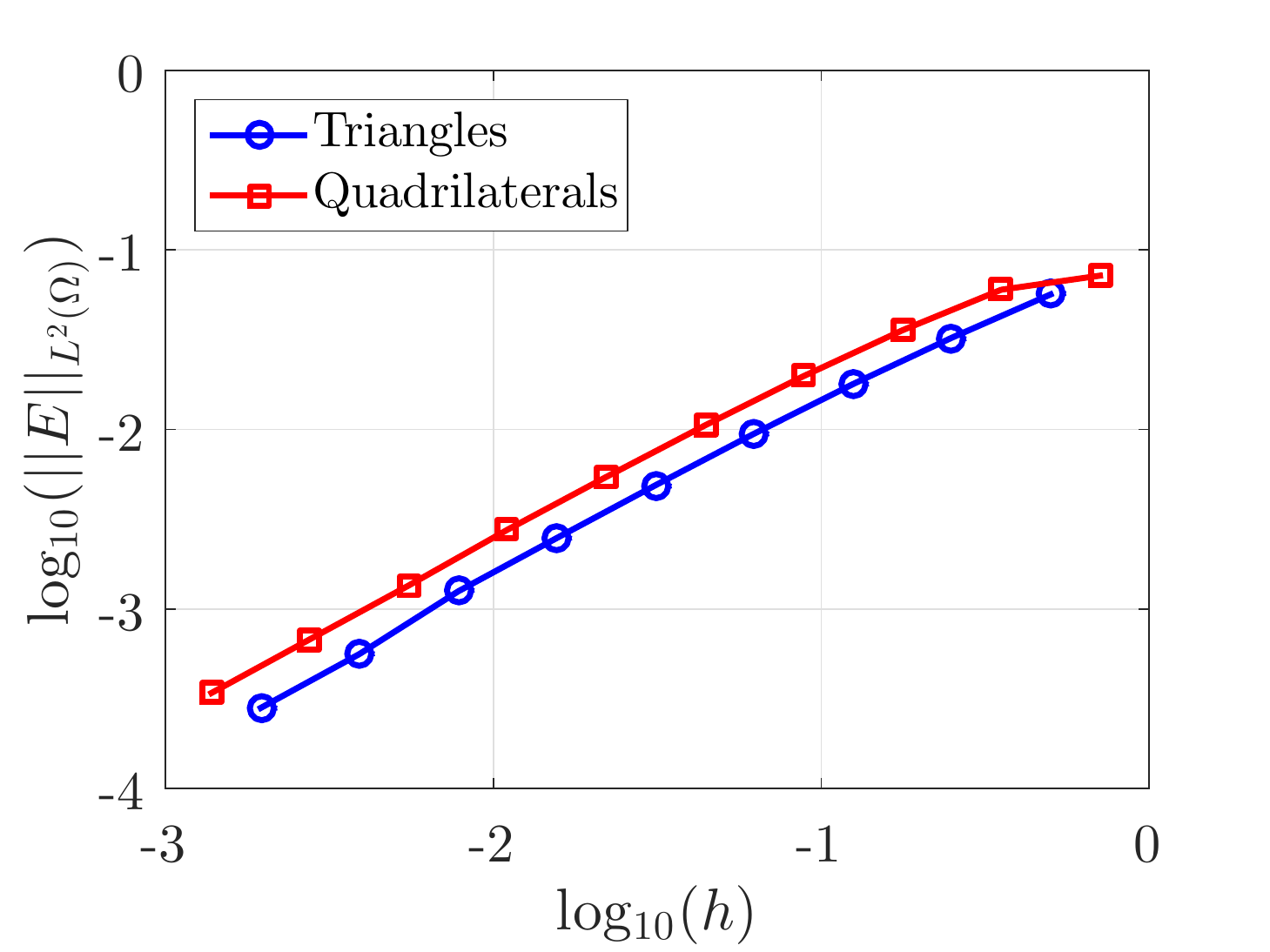}}
	\subfigure[$\bq$]{\includegraphics[width=0.4\textwidth]{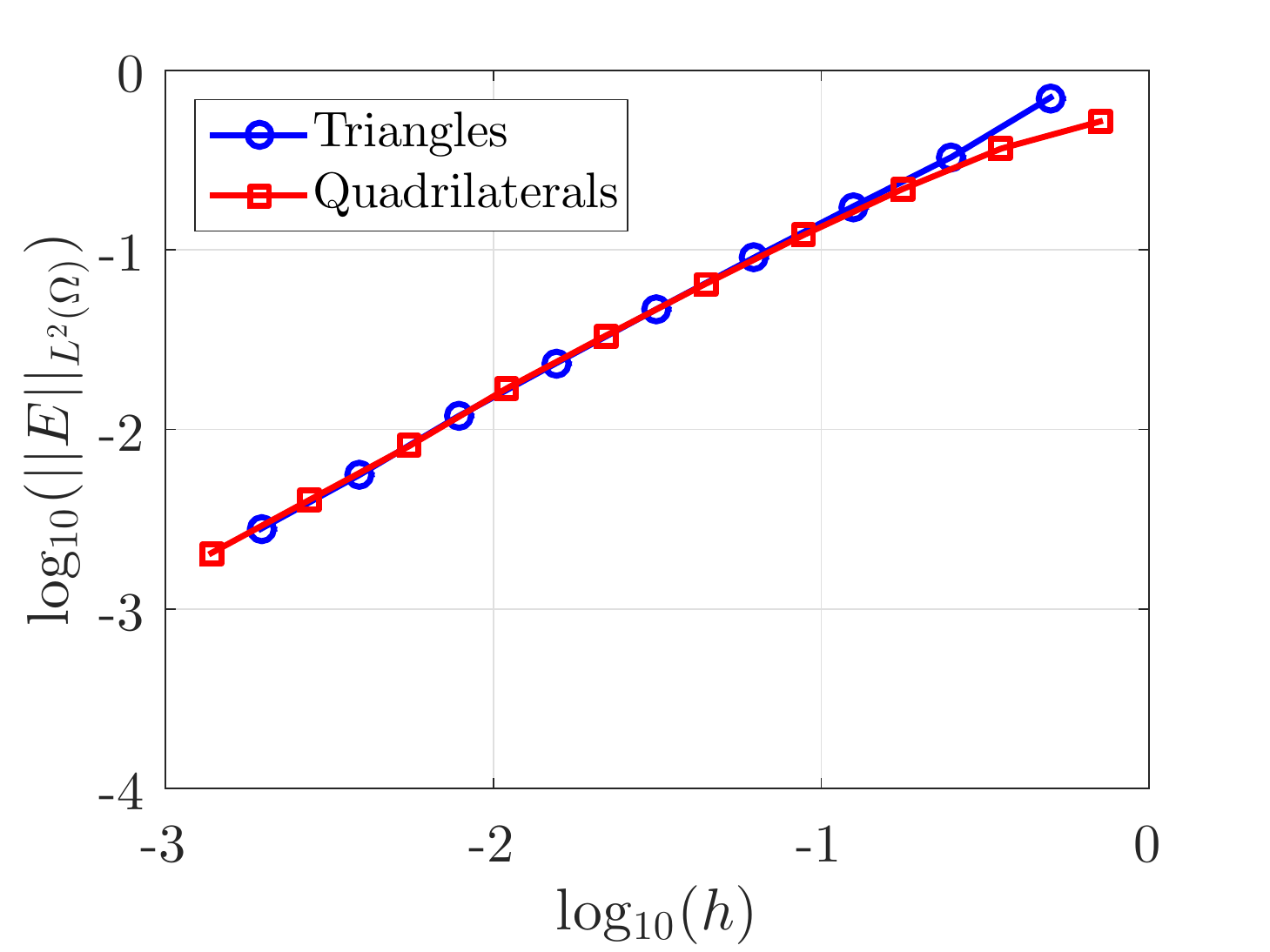}}
	\caption{Mesh convergence of the error of the solution and its gradient in the $\eltwo(\Omega)$ norm for the 2D Poisson problem.}
	\label{fig:poisson2D_hConv}
\end{figure}

The convergence of the error of the primal variable $u$, measured in the $\eltwo(\Omega)$ norm, as a function of the characteristic element size $h$ is depicted in Figure~\ref{fig:poisson2D_hConv} (a) for both triangular and quadrilateral elements. Similarly, the convergence of the error of the dual variable $\bq$, measured in the $\eltwo(\Omega)$ norm, as a function of the characteristic element size $h$ is depicted in Figure~\ref{fig:poisson2D_hConv} (b) for both triangular and quadrilateral elements. 
In all the examples the characteristic element size is defined as the maximum diameter of all elements,
\begin{equation}
h = \max_{e} \left\{ \text{diam}(\Omega_e) \right\}.
\end{equation}
For the regular meshes considered here, $h$ corresponds to the diagonal of a quadrilateral element or the largest edge of a triangle.

The results confirm the expected linear rate of convergence for both variables and by using triangular and quadrilateral elements.

Next, a three dimensional test case is considered. The computational domain is $\Omega=[0,1]^3$ and the source term is selected so that the analytical solution is 
\begin{equation}
u(x,y) = \exp \big( \alpha \sin(ax+cy+ez) + \beta\cos(bx+dy+fz) \big),
\end{equation}
with $\alpha=0.1$, $\beta=0.3$, $a=5.1$, $b=4.3$, $c=-6.2$, $d=3.4$, $e=1.8$ and $f=1.7$. Neumann boundary conditions, corresponding to the analytical normal flux, are imposed in $\Gamma_N = \{(x,y,z) \in \mathbb{R}^3 \; | \; z=0\}$ and Dirichlet boundary conditions, corresponding to the analytical solution, are imposed in $\Gamma_D = \partial \Omega \setminus \Gamma_N$.

\begin{figure}[!tb]
	\centering
	\subfigure[Hexahedrons]{\includegraphics[width=0.24\textwidth]{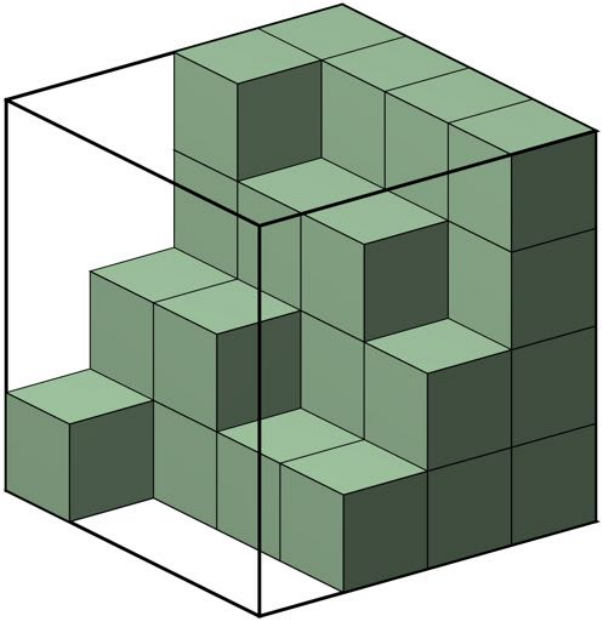}}
	\subfigure[Tetrahedrons]{\includegraphics[width=0.24\textwidth]{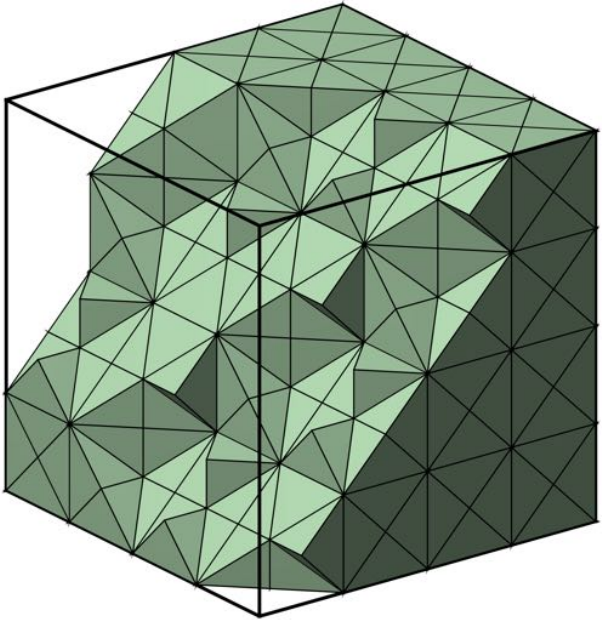}}		
	\subfigure[Prisms]{\includegraphics[width=0.24\textwidth]{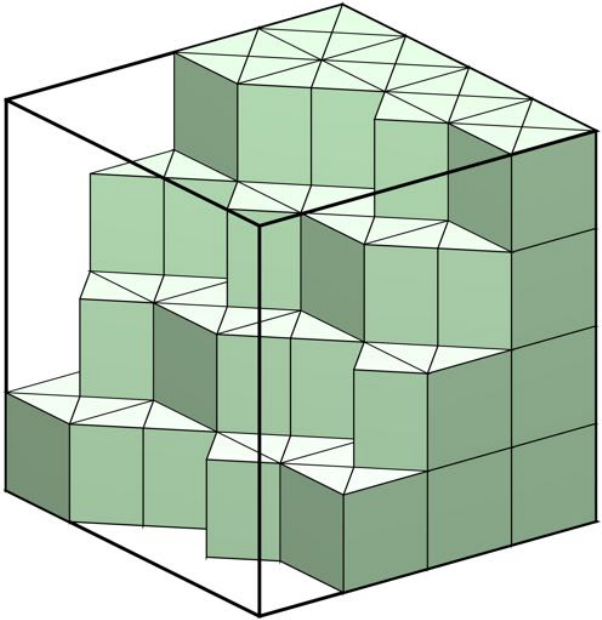}}
	\subfigure[Pyramids]{\includegraphics[width=0.24\textwidth]{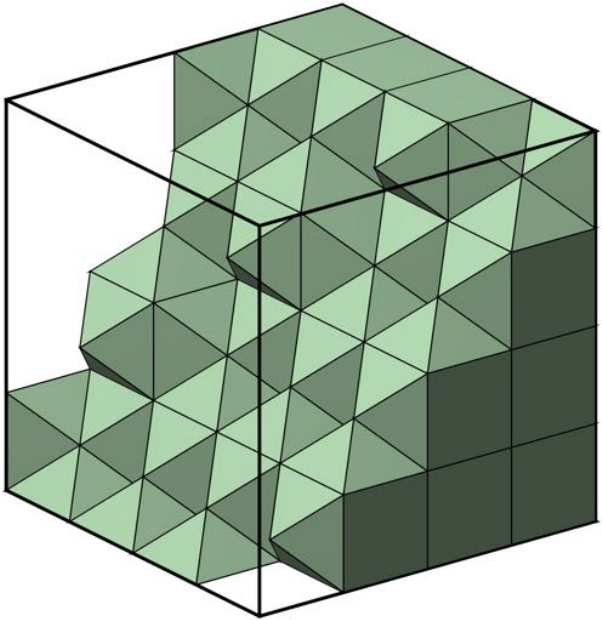}}
	\caption{Third level of mesh refinement for the meshes of $\Omega=[0,1]^3$ employed to test the optimal convergence in 3D.}
	\label{fig:poisson3D_meshes}
\end{figure}

The convergence study is performed for regular meshes of hexahedral, tetrahedral, prismatic and pyramidal elements. A cut through the meshes corresponding to the third level of refinement is represented in Figure~\ref{fig:poisson3D_meshes} for all element types.
Tetrahedral meshes are obtained from the corresponding hexahedral mesh by subdividing each hexahedron into 24 tetrahedrons. Similarly, prismatic meshes are obtained by subdividing each hexahedron into six prisms and pyramidal meshes are obtained by subdividing each hexahedron into four pyramids.

\begin{figure}[!tb]
	\centering
	\subfigure[$u$]{\includegraphics[width=0.4\textwidth]{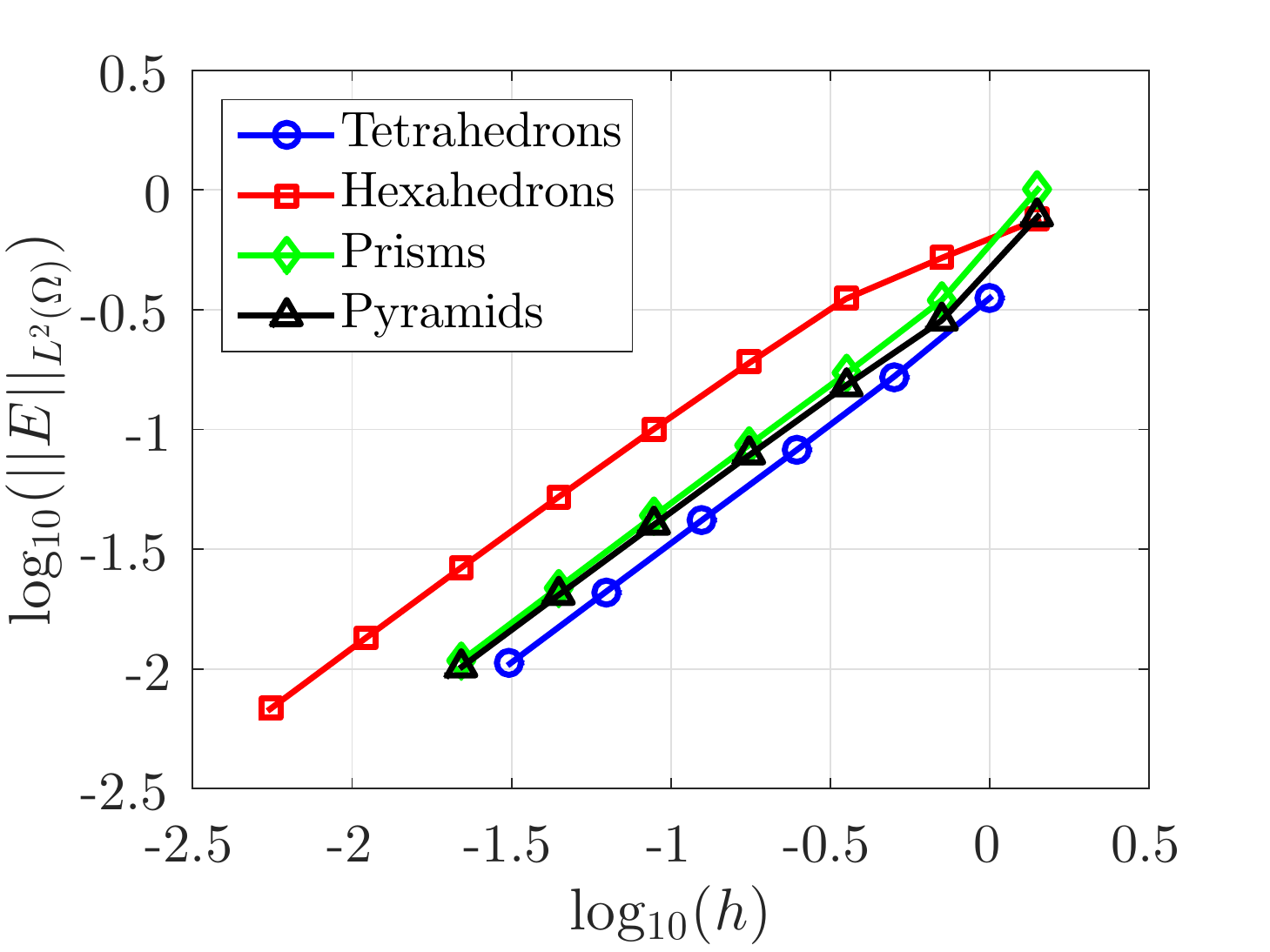}}
	\subfigure[$\bq$]{\includegraphics[width=0.4\textwidth]{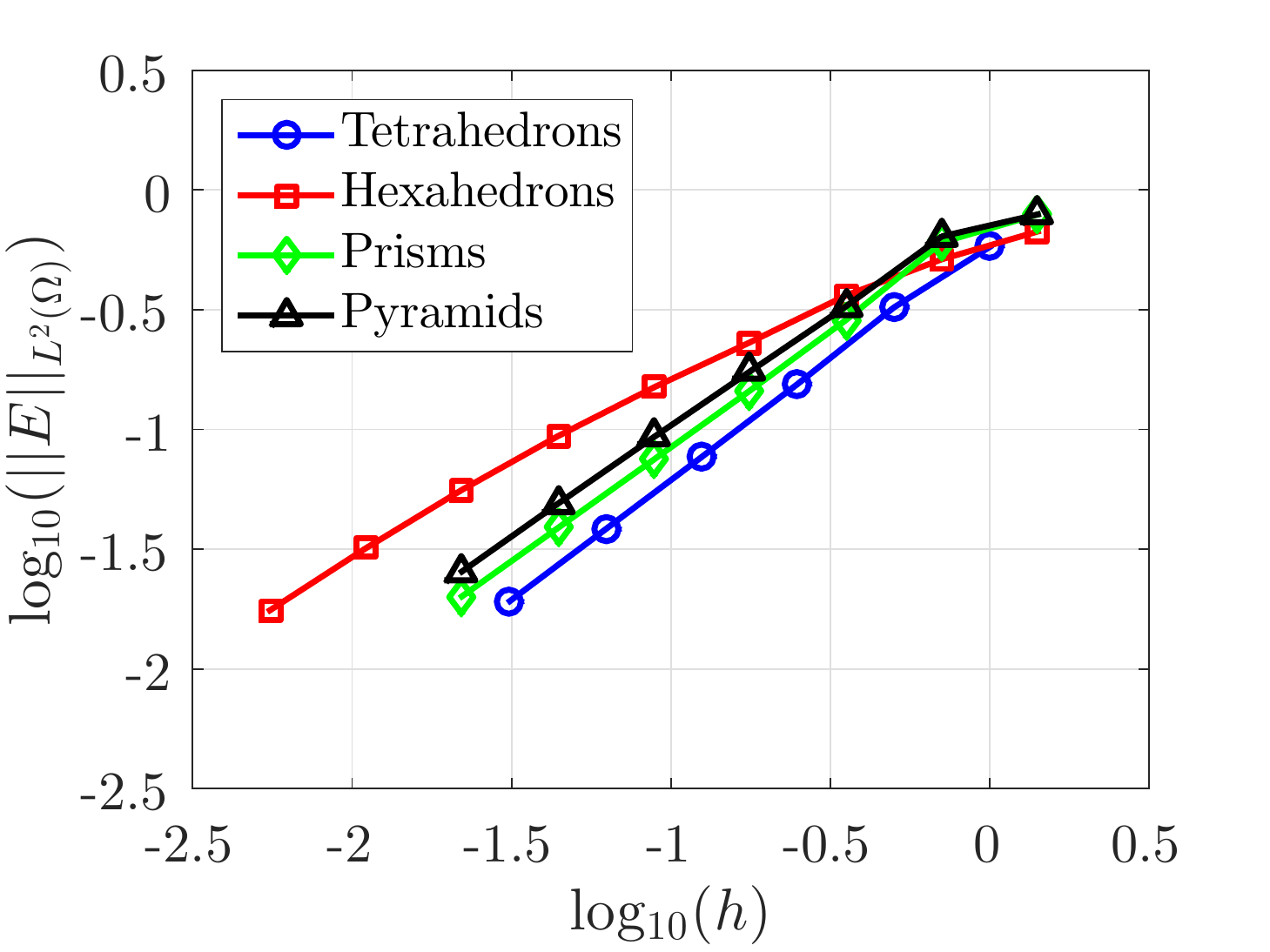}}
	\caption{Mesh convergence of the error of the solution and its gradient in the $\eltwo(\Omega)$ norm for the 3D Poisson problem.}
	\label{fig:poisson3D_hConv}
\end{figure}

The convergence of the error of the primal and dual variables, measured in the $\eltwo(\Omega)$ norm, as a function of the characteristic element size $h$ is depicted in Figure~\ref{fig:poisson3D_hConv}. 
The results show a linear rate of convergence for both the dual and primal variables and using all the different types of elements.

\subsection{Optimal convergence of the FCFV scheme for Stokes equation}
\label{sc:StokesVerification}

The Stokes problem given in Equation \eqref{eq:Stokes} is considered in two dimensions to verify the optimal convergence properties of the FCFV method for saddle-point problems. The convergence analysis for a classical benchmark case of 3D Stokes solvers is presented in Section~\ref{sc:stokesSphere}, not only analysing the convergence of the primal and dual variables but also the convergence of the drag force as the mesh is refined.

A two dimensional synthetic problem, taken from~~\cite{donea2003finite}, is considered in the domain $\Omega = [0,1]^2$. The boundary $\partial\Omega$ is split into two disjoint parts, namely $\Gamma_N = \{(x,y) \in \RR^2 \; | \; y=0\}$ where a pseudo-traction $\bm{t}$ is imposed and $\Gamma_D = \partial \Omega \setminus \Gamma_N$ where a velocity profile $\bu_D$ is set.
The viscosity parameter is set to $\nu=1$ and the source term $\bm{s}$ and the boundary data $\bm{t}$ and $\bu_D$ are chosen such that the analytical solution is
\begin{equation}
\left\{
\begin{aligned}
u_1(x,y) & = x^2 (1-x)^2 (2y-6y^2+4y^3) ,\\
u_2(x,y) & = - y^2 (1-y)^2 (2x-6x^2+4x^3) , \\
p(x,y) & = x (1-x) ,
\end{aligned}
\right.
\end{equation}
where $u_1$ and $u_2$ are the two components of the velocity field vector $\bu$.

\begin{figure}[!tb]
	\centering
	\subfigure[Velocity, Mesh 3]{\includegraphics[width=0.24\textwidth]{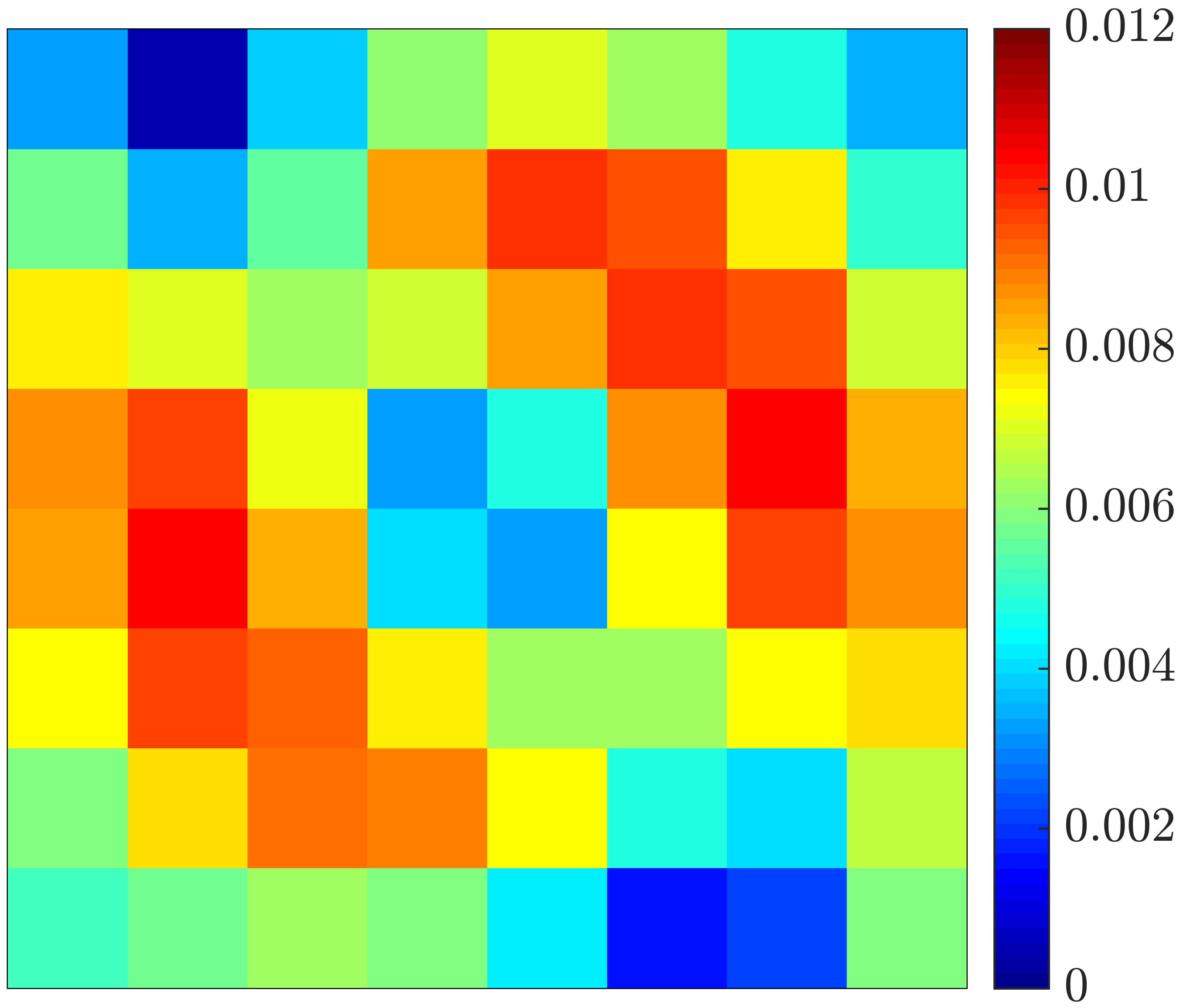}}
	\subfigure[Velocity, Mesh 5]{\includegraphics[width=0.24\textwidth]{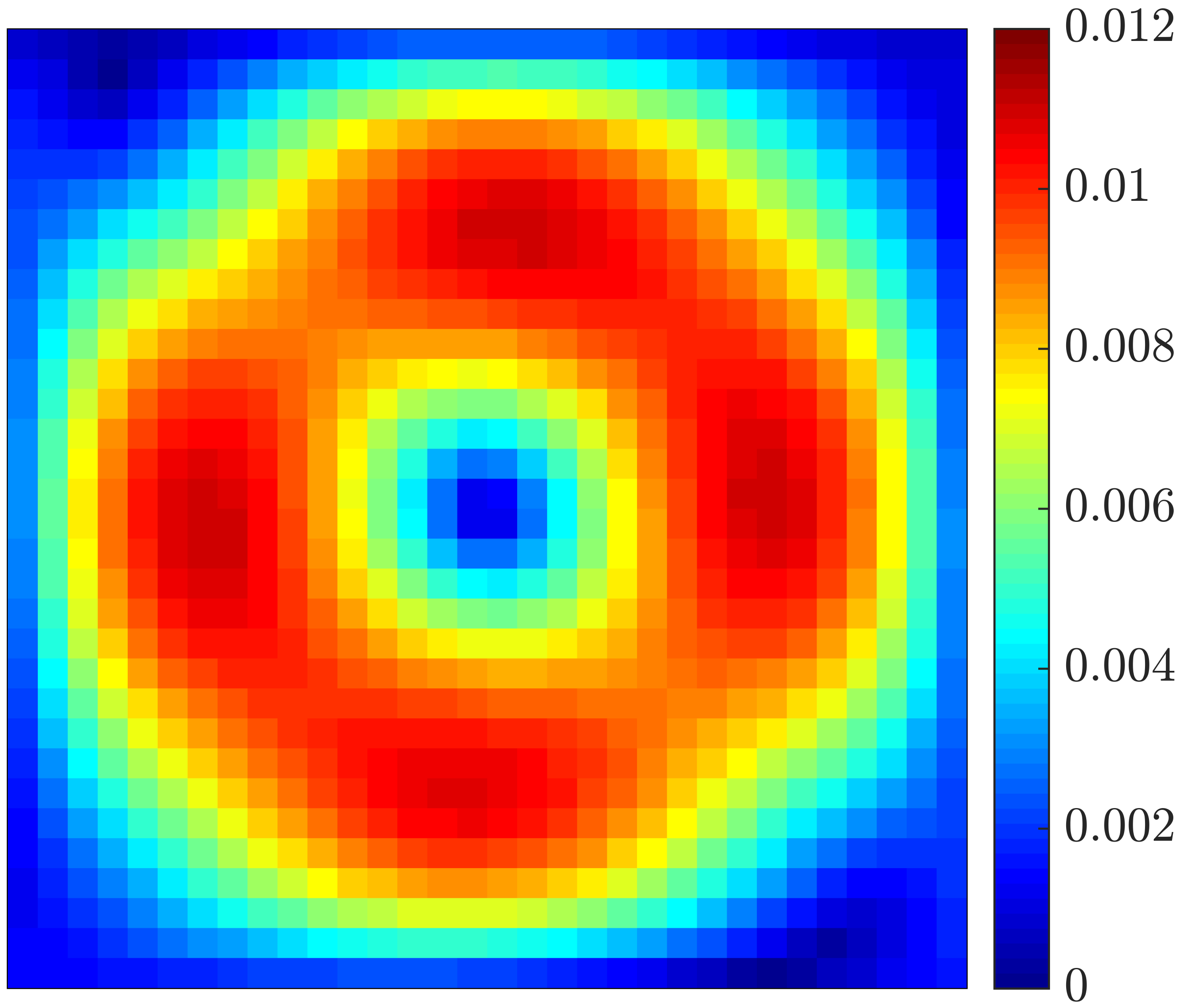}}
	\subfigure[Velocity, Mesh 7]{\includegraphics[width=0.24\textwidth]{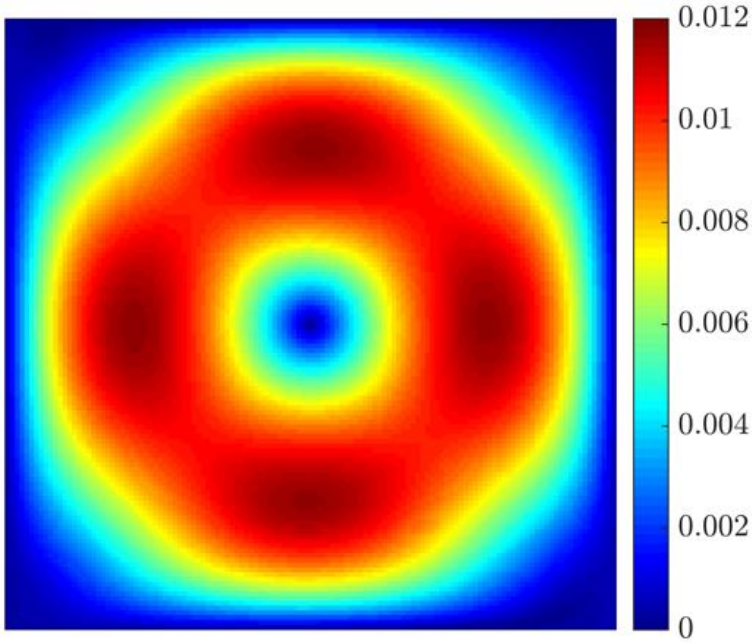}}
	\subfigure[Velocity, Mesh 9]{\includegraphics[width=0.24\textwidth]{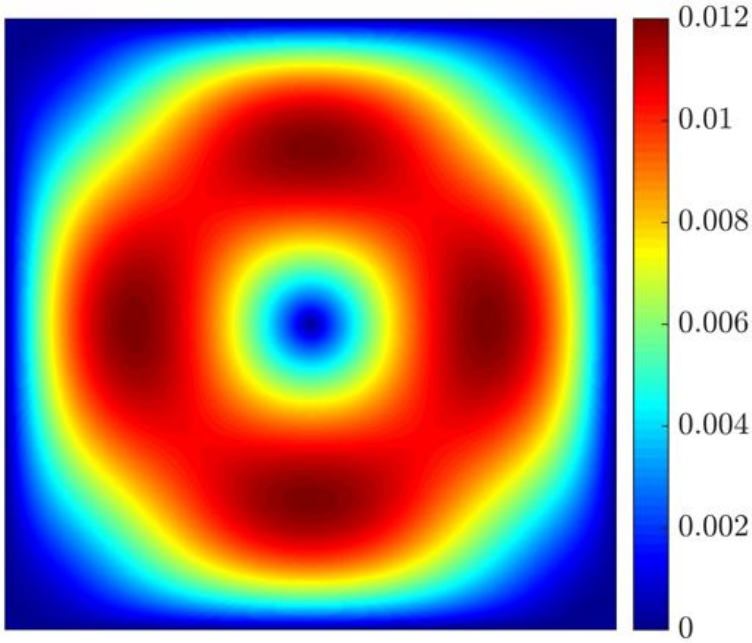}}	
	\subfigure[Pressure, Mesh 3]{\includegraphics[width=0.24\textwidth]{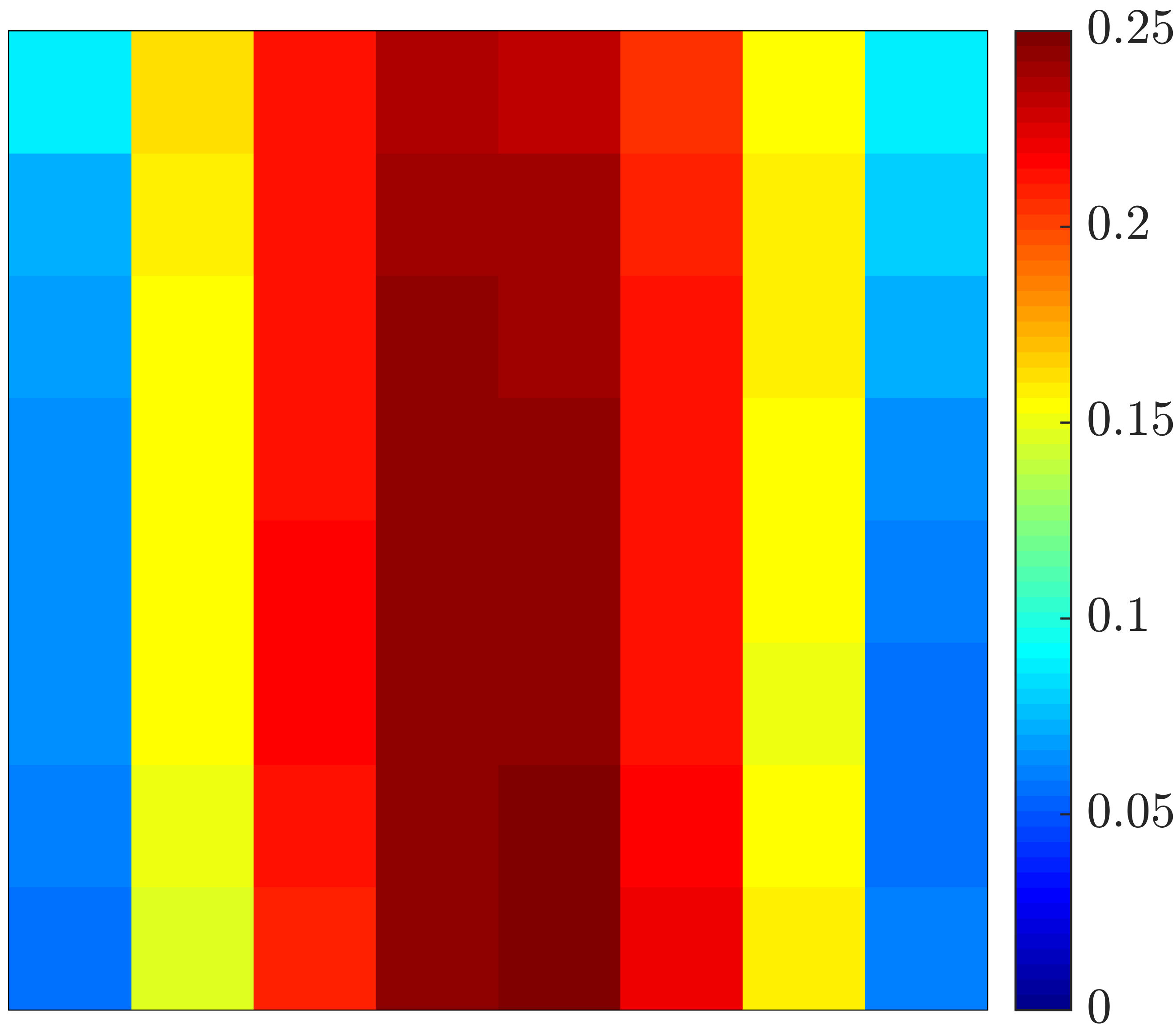}}
	\subfigure[Pressure, Mesh 5]{\includegraphics[width=0.24\textwidth]{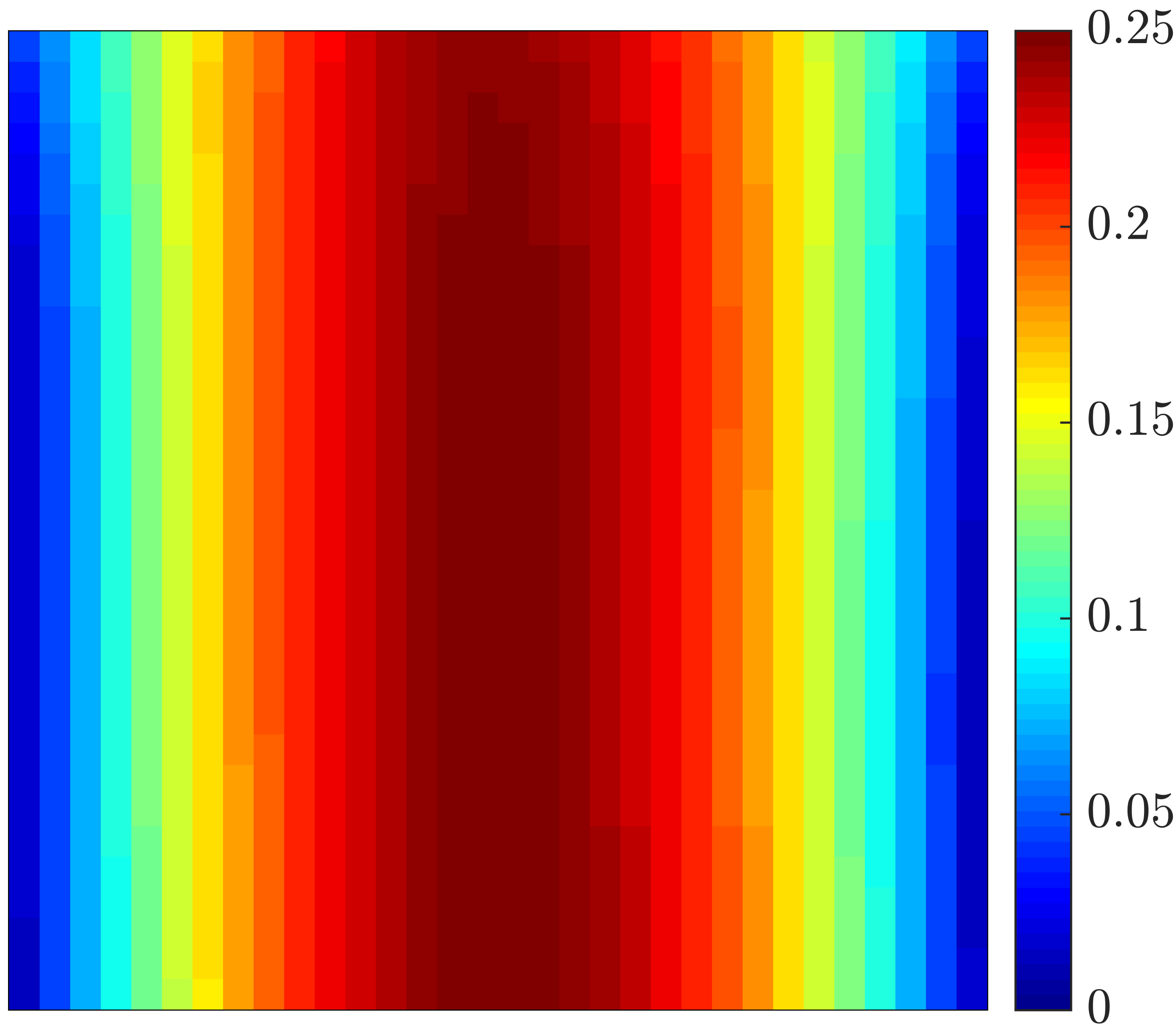}}
	\subfigure[Pressure, Mesh 7]{\includegraphics[width=0.24\textwidth]{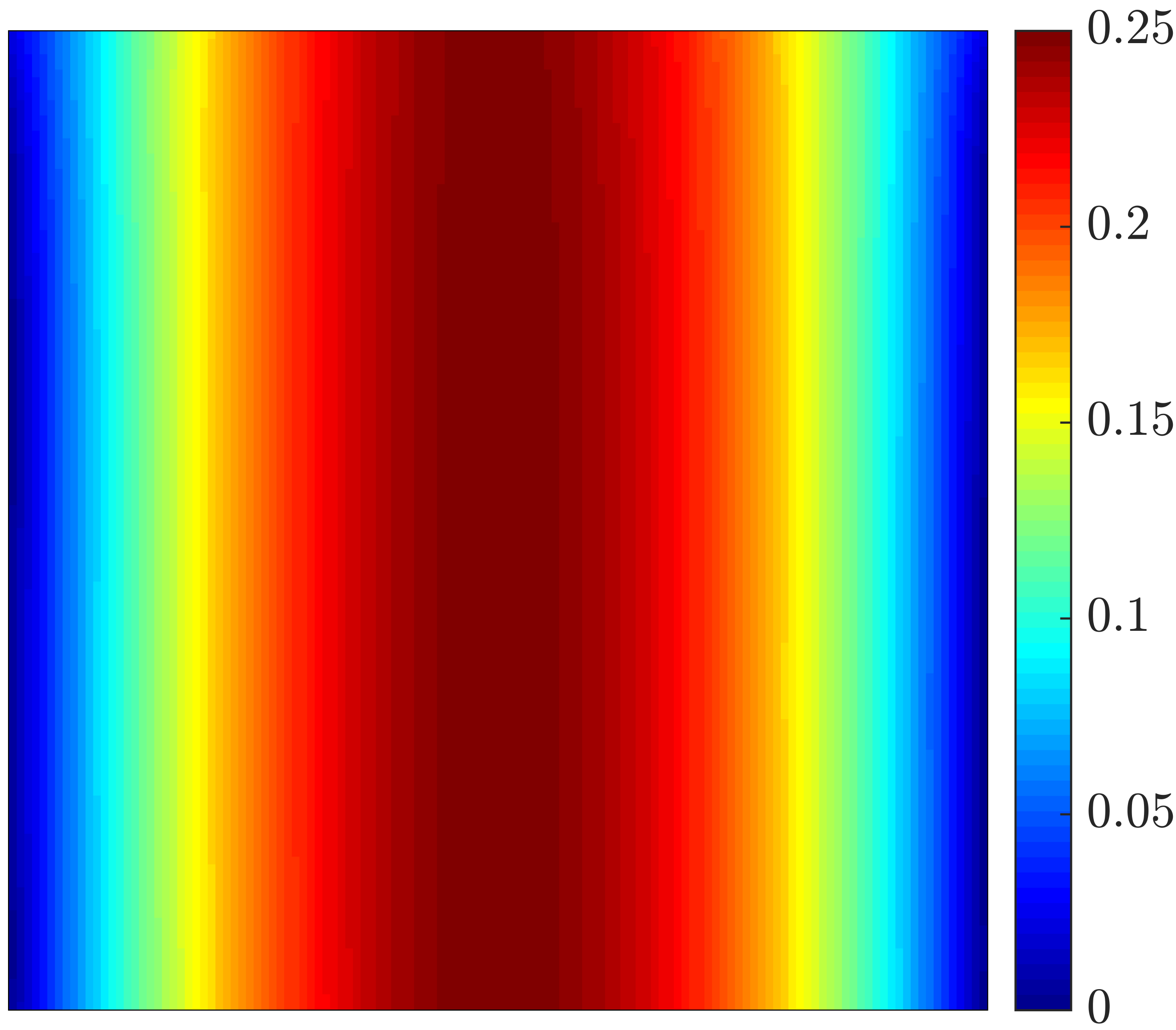}}
	\subfigure[Pressure, Mesh 9]{\includegraphics[width=0.24\textwidth]{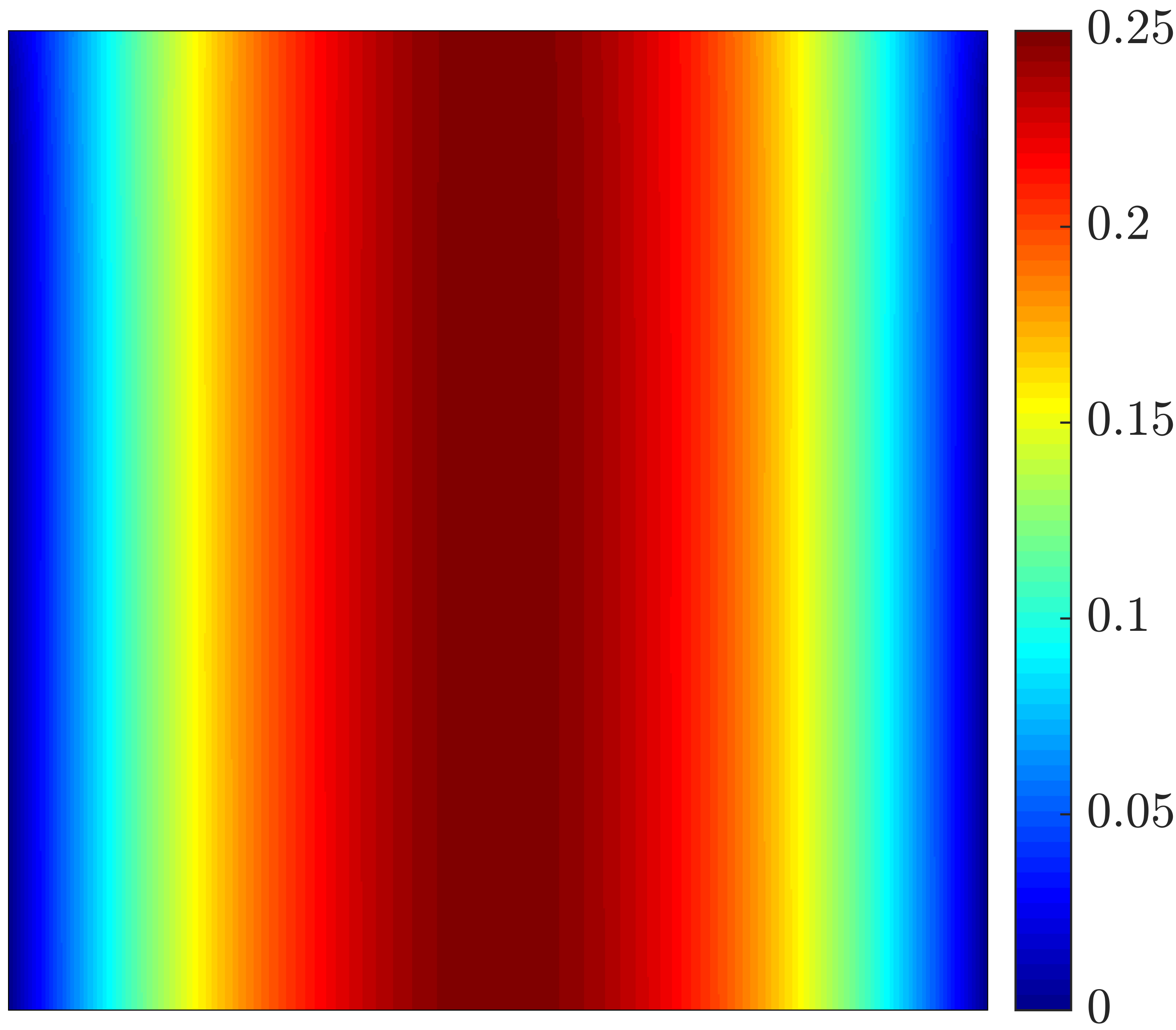}}	
	\caption{Solution of the 2D Stokes problem using quadrilateral meshes.}
	\label{fig:stokes2D_QUAsol}
\end{figure}

\begin{figure}[!tb]
	\centering
	\subfigure[Velocity, Mesh 3]{\includegraphics[width=0.24\textwidth]{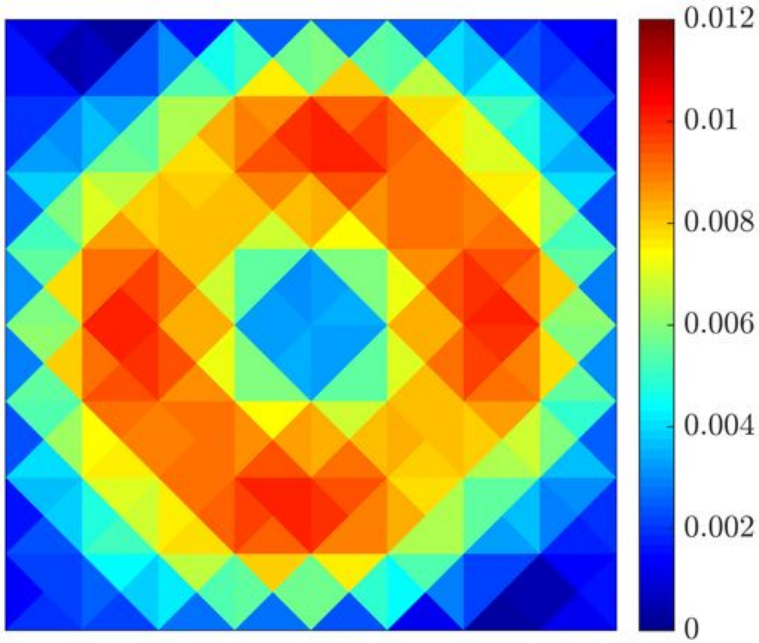}}
	\subfigure[Velocity, Mesh 5]{\includegraphics[width=0.24\textwidth]{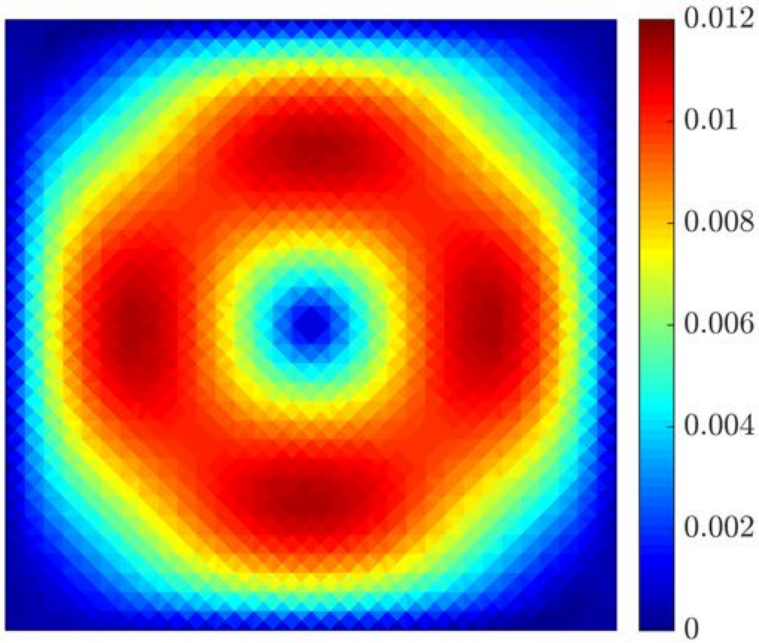}}
	\subfigure[Velocity, Mesh 7]{\includegraphics[width=0.24\textwidth]{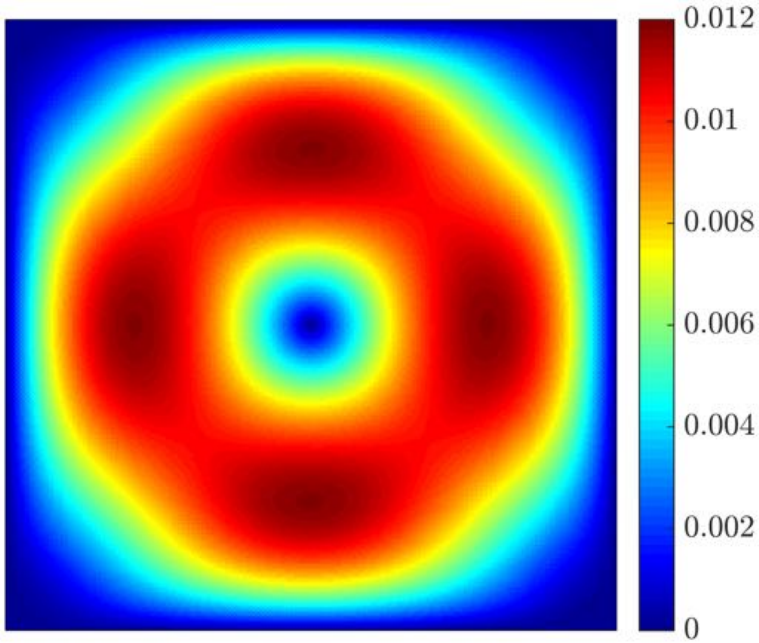}}
	\subfigure[Velocity, Mesh 9]{\includegraphics[width=0.24\textwidth]{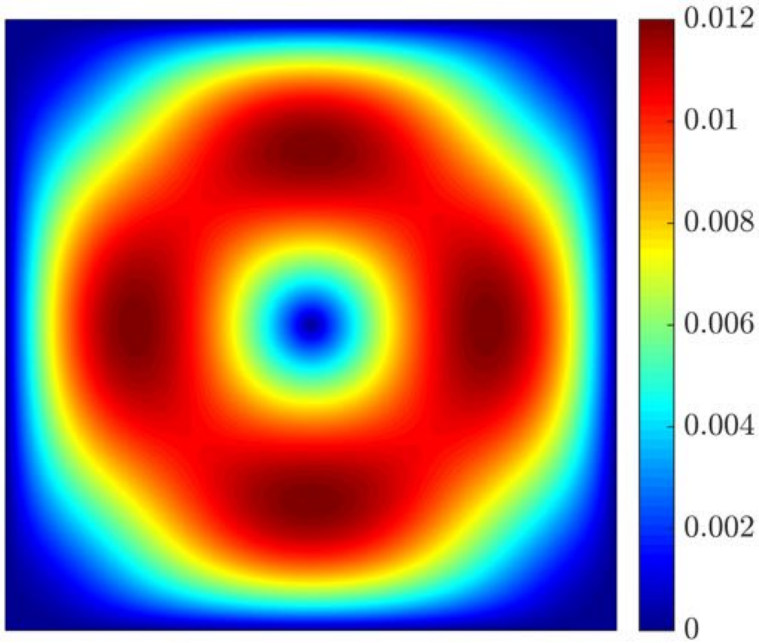}}
	\subfigure[Pressure, Mesh 3]{\includegraphics[width=0.24\textwidth]{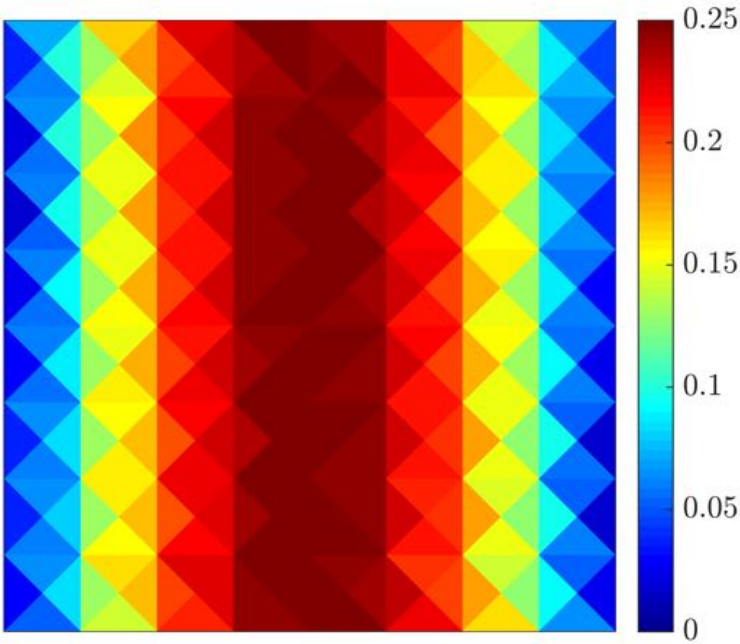}}
	\subfigure[Pressure, Mesh 5]{\includegraphics[width=0.24\textwidth]{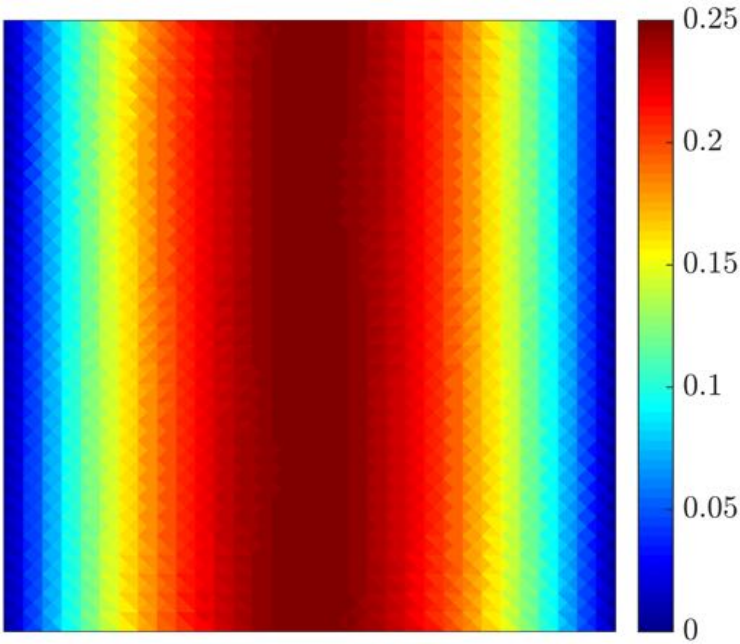}}
	\subfigure[Pressure, Mesh 7]{\includegraphics[width=0.24\textwidth]{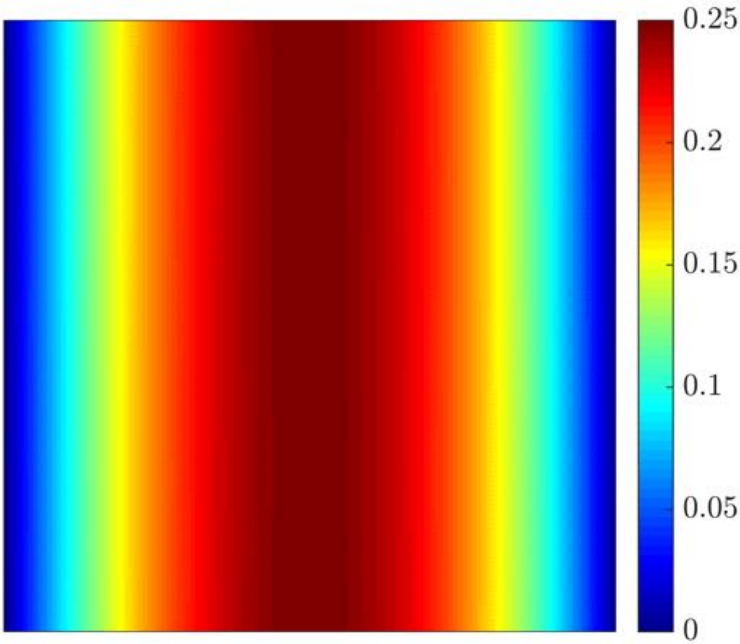}}
	\subfigure[Pressure, Mesh 9]{\includegraphics[width=0.24\textwidth]{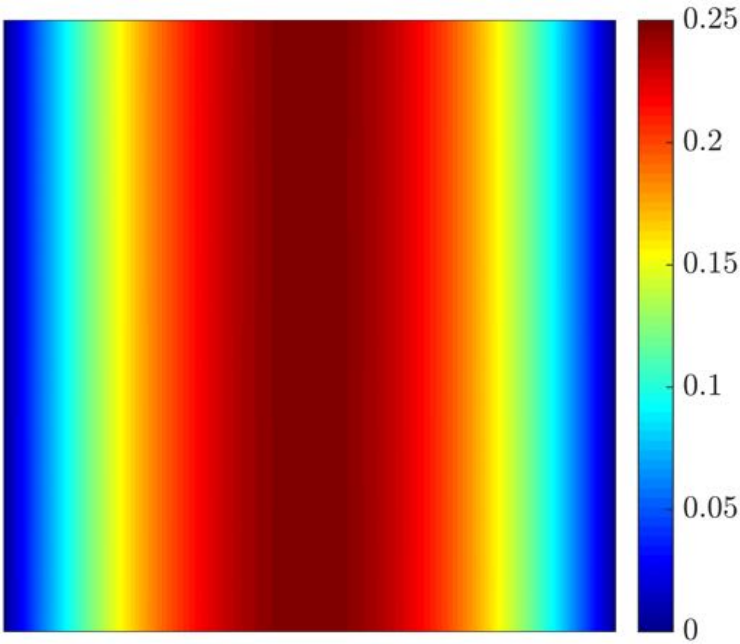}}
	\caption{Solution of the 2D Stokes problem using triangular meshes.}
	\label{fig:stokes2D_TRIsol}
\end{figure}

The same meshes used in the two dimensional verification example of Section~\ref{sc:PoissonVerification} are considered. 
Figures~\ref{fig:stokes2D_QUAsol} and \ref{fig:stokes2D_TRIsol} show the numerical solutions computed using the FCFV method on quadrilateral and triangular meshes respectively for different levels of mesh refinement. Once more, the piecewise constant nature of the FCFV approximation is clearly observed as well as the improved resolution as the mesh is refined.

\begin{figure}[!tb]
	\centering
	\subfigure[$p$]{\includegraphics[width=0.32\textwidth]{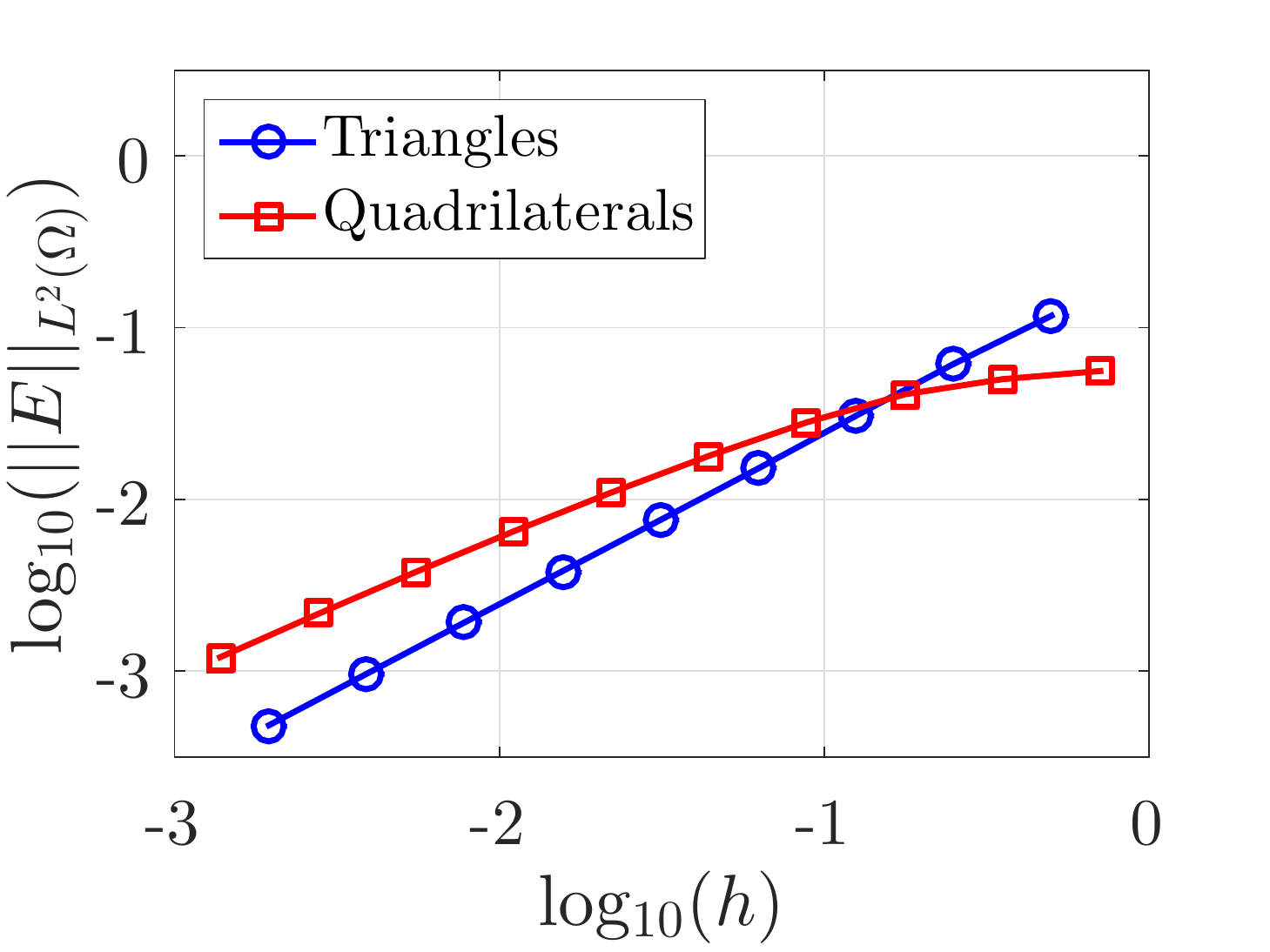}}
	\subfigure[$\bu$]{\includegraphics[width=0.32\textwidth]{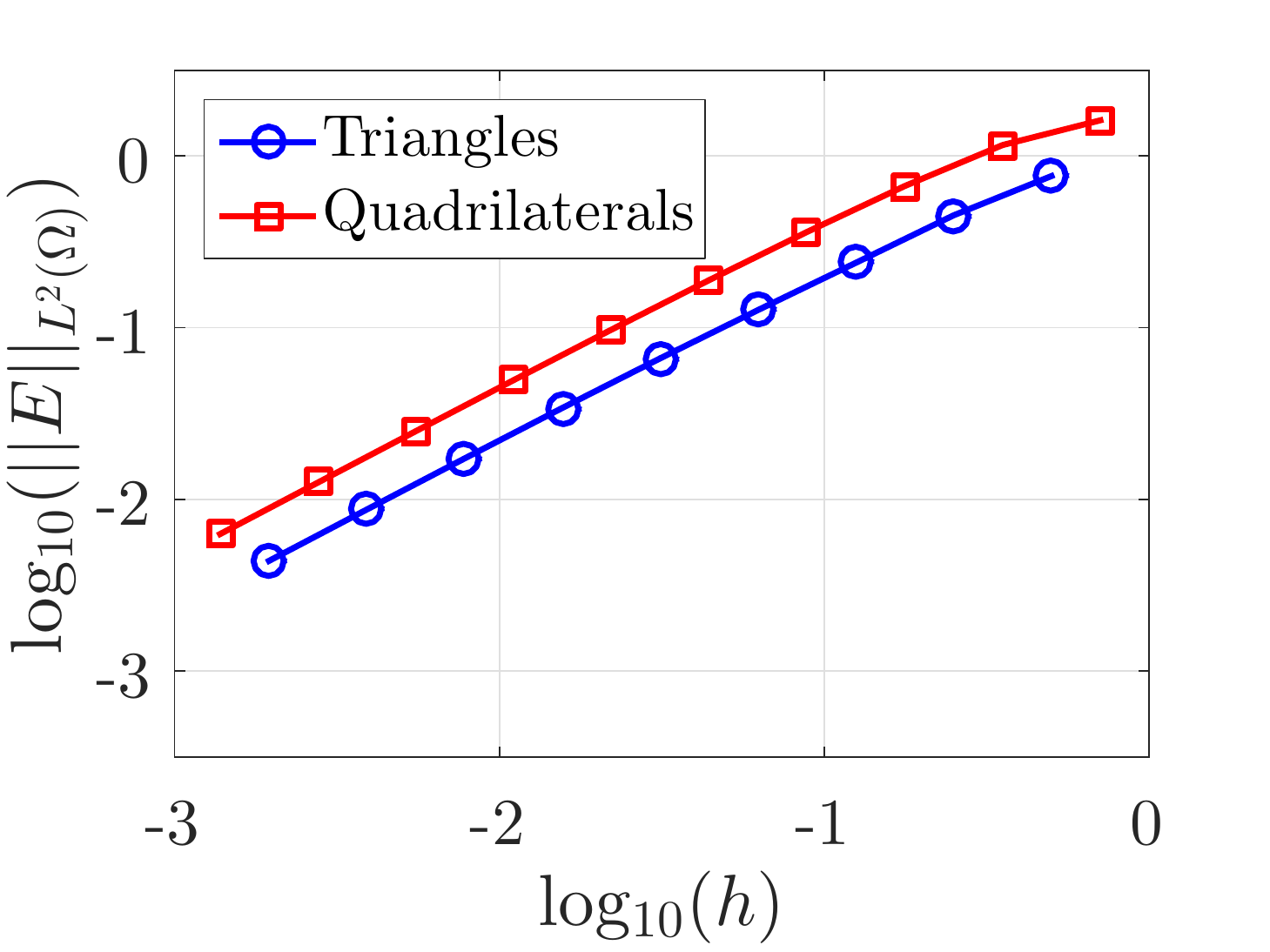}}
	\subfigure[$\bm{L}$]{\includegraphics[width=0.32\textwidth]{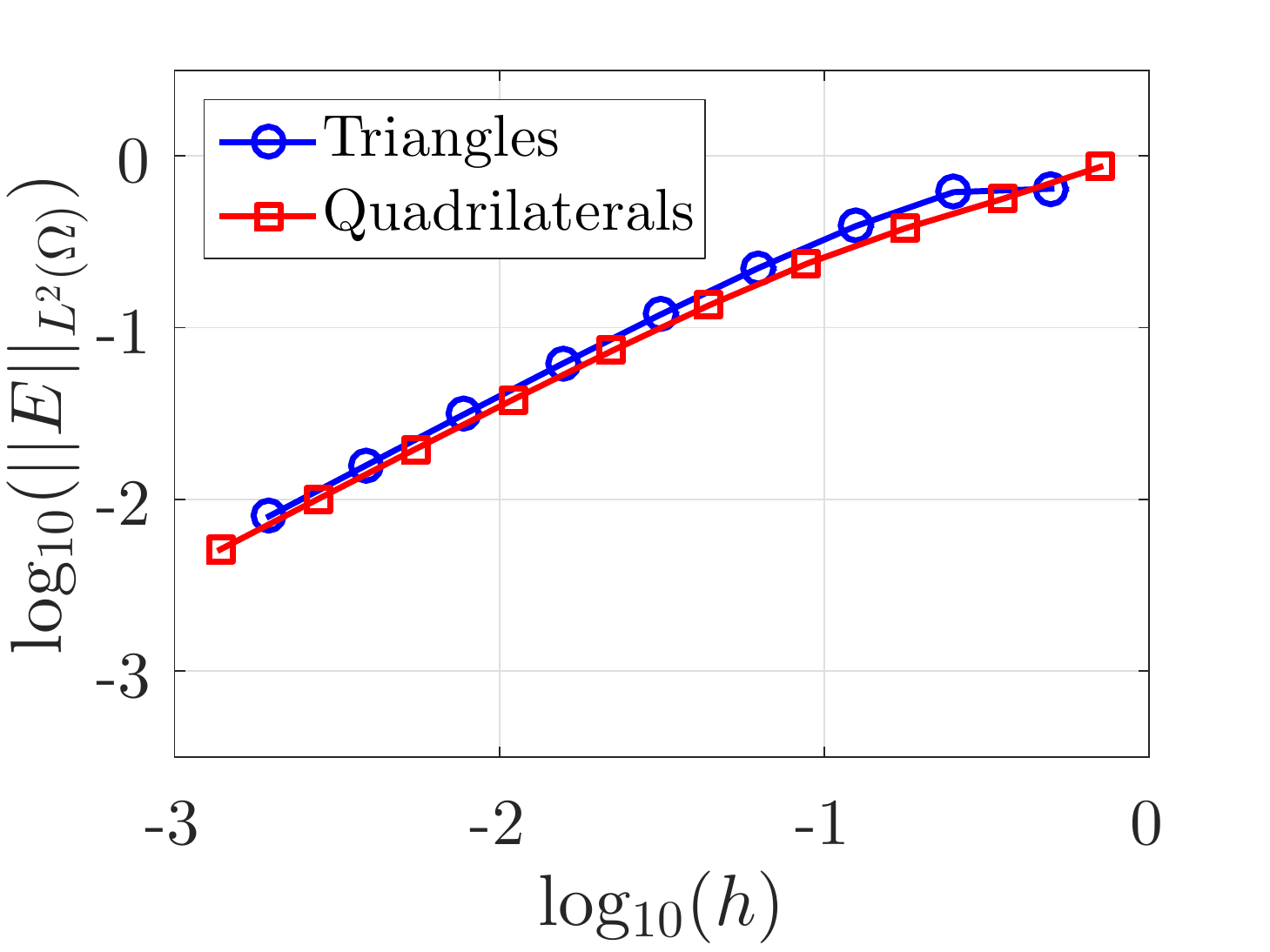}}
	\caption{Mesh convergence of the error of the pressure, the velocity and the velocity gradient in the $\eltwo(\Omega)$ norm for the 2D Stokes problem.}
	\label{fig:stokes2D_hConv}
\end{figure}

The convergence of the error of the pressure, velocity and velocity gradient, measured in the $\eltwo(\Omega)$ norm, as a function of the characteristic element size $h$ is depicted in Figure~\ref{fig:stokes2D_hConv} for both triangular and quadrilateral elements. 
The expected linear rate of convergence is observed for the primal, $\bu,p$ and for the dual, $\bm{L}$, variables using both triangular and quadrilateral meshes. It is worth noting that for quadrilateral meshes both $\bu$ and $\bm{L}$ converge with the optimal linear rate whereas the pressure converges with a slightly lower rate, 0.9 in this example.

\subsection{Computational cost for different element types}
\label{sc:ComputationalCost}

\begin{figure}[!tb]
	\centering
	\subfigure[$u$]{\includegraphics[width=0.4\textwidth]{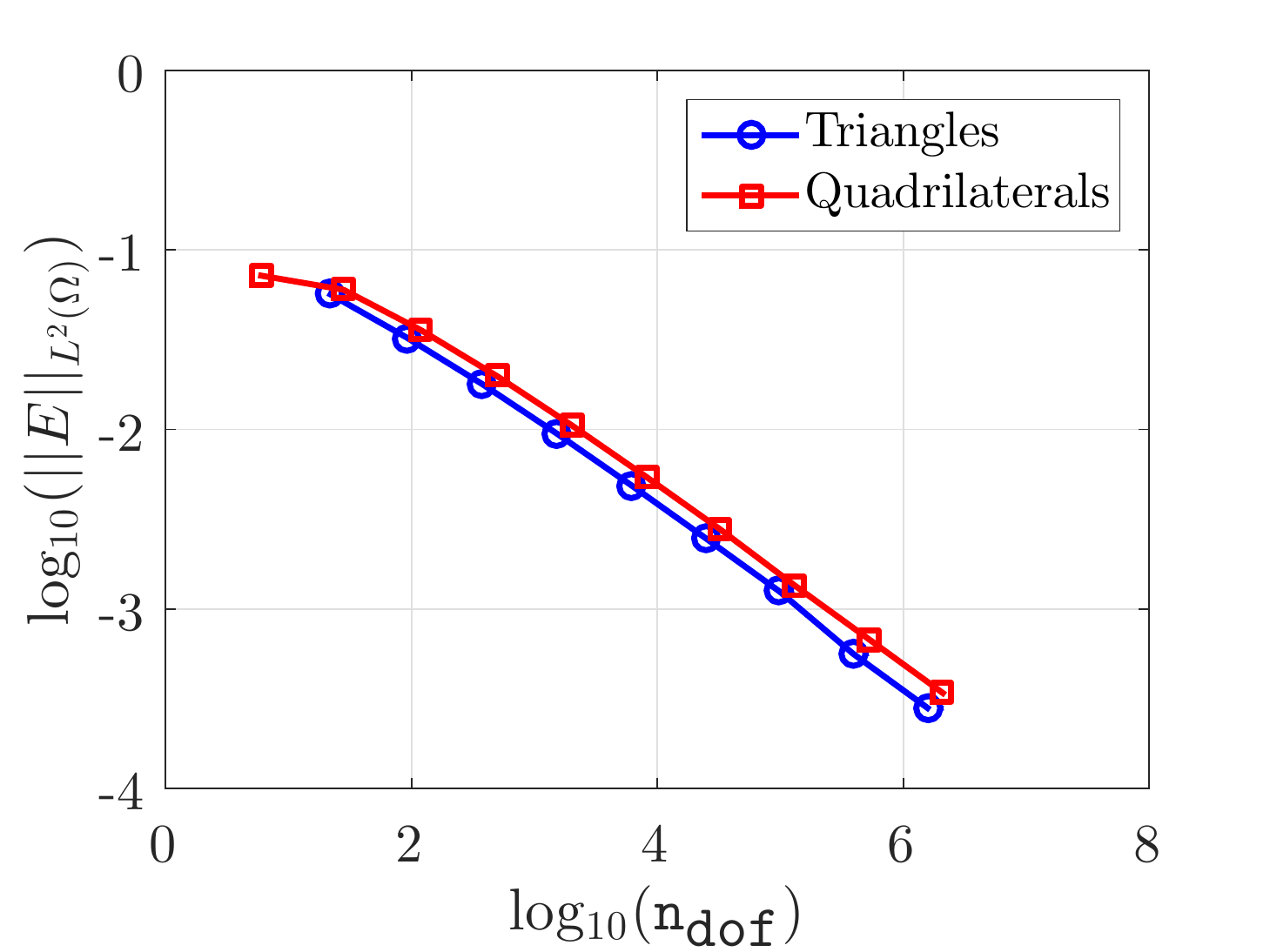}}
	\subfigure[$\bq$]{\includegraphics[width=0.4\textwidth]{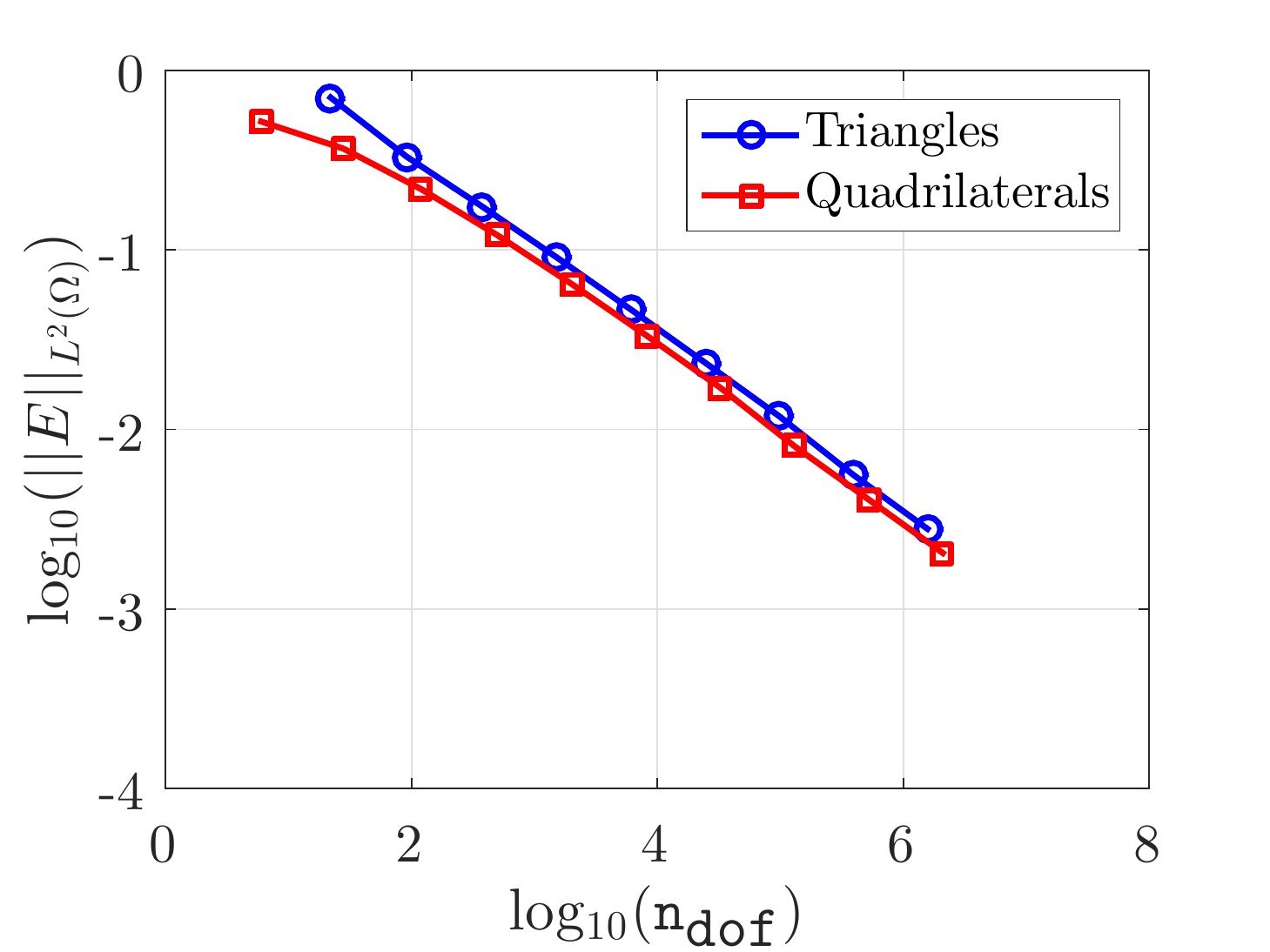}}
	\caption{Error of the solution and its gradient in the $\eltwo(\Omega)$ norm as a function of the number of degrees of freedom of the global system for the 2D Poisson problem.}
	\label{fig:poisson2D_ndof}
\end{figure}

The results in Figures~\ref{fig:poisson2D_hConv} and \ref{fig:stokes2D_hConv} show that, for the same level of mesh refinement, the FCFV method with triangular meshes provides more accurate results than using quadrilateral meshes, for the solution of both Poisson and Stokes problems. 
This is mainly because, for the same level of mesh refinement, triangular meshes have four times more internal faces than the corresponding quadrilateral meshes. This means that the second quadrilateral mesh has exactly the same number of internal faces as the first triangular mesh. Therefore, a fair comparison between triangular and quadrilateral meshes shows that for the same number of degrees of freedom (\ndof) of the global problem (i.e. number of internal faces plus number of faces on the Neumann boundary), triangular and quadrilateral meshes provide similar accuracy for both the primal and the dual variables. This can be observed in Figure~\ref{fig:poisson2D_ndof}, where the evolution of the error of the primal and dual variables is represented as a function of the number of degrees of freedom of the global system of equations.
From the results in Figure~\ref{fig:poisson2D_ndof}, it can be observed that triangular and quadrilateral elements provide almost the same accuracy for a given number of degrees of freedom, with triangular elements providing a marginal extra accuracy for the primal variable and with quadrilateral elements providing a marginal extra accuracy for the dual variable. 

\begin{figure}[!tb]
	\centering
	\subfigure[$u$]{\includegraphics[width=0.4\textwidth]{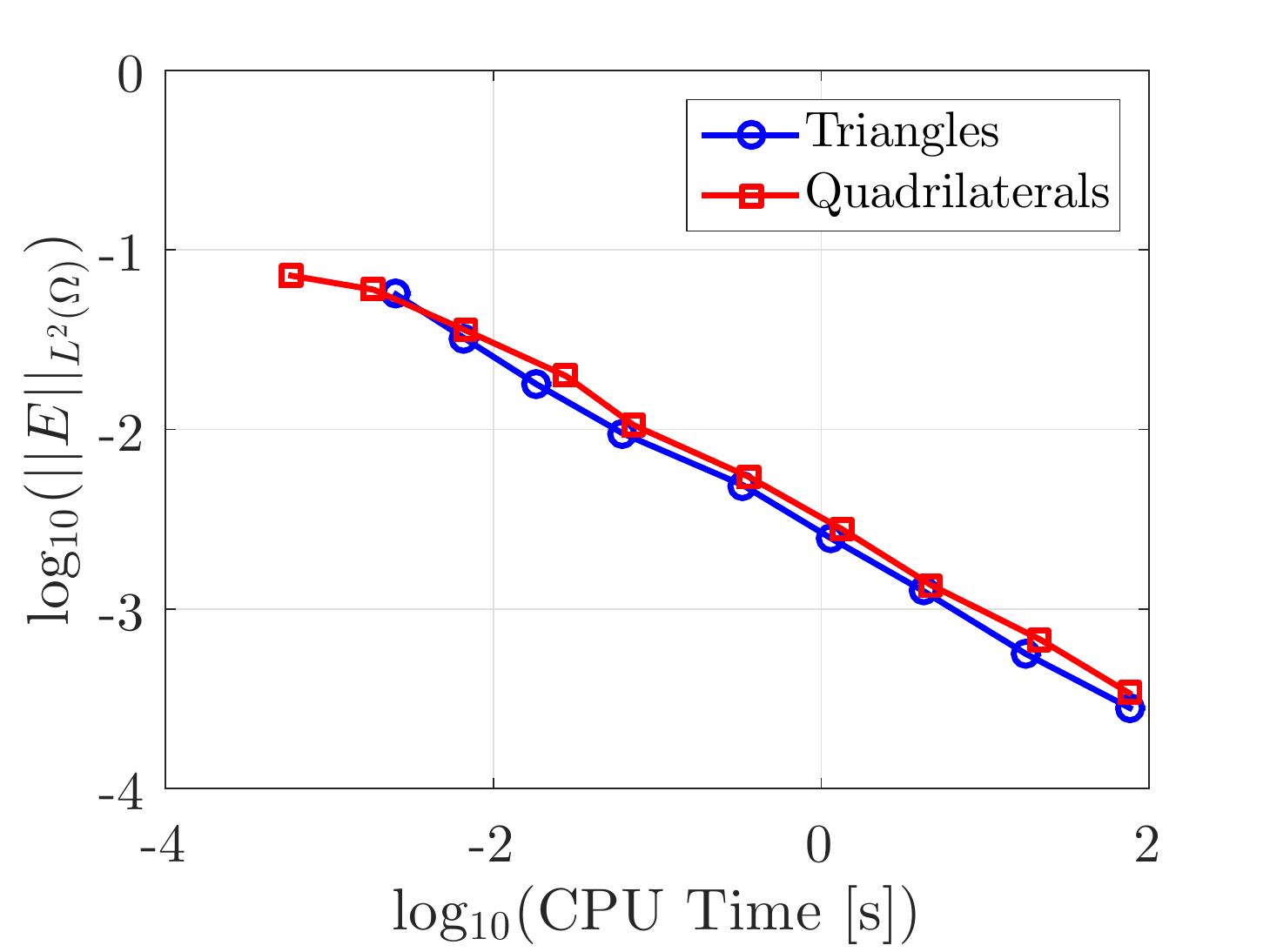}}
	\subfigure[$\bq$]{\includegraphics[width=0.4\textwidth]{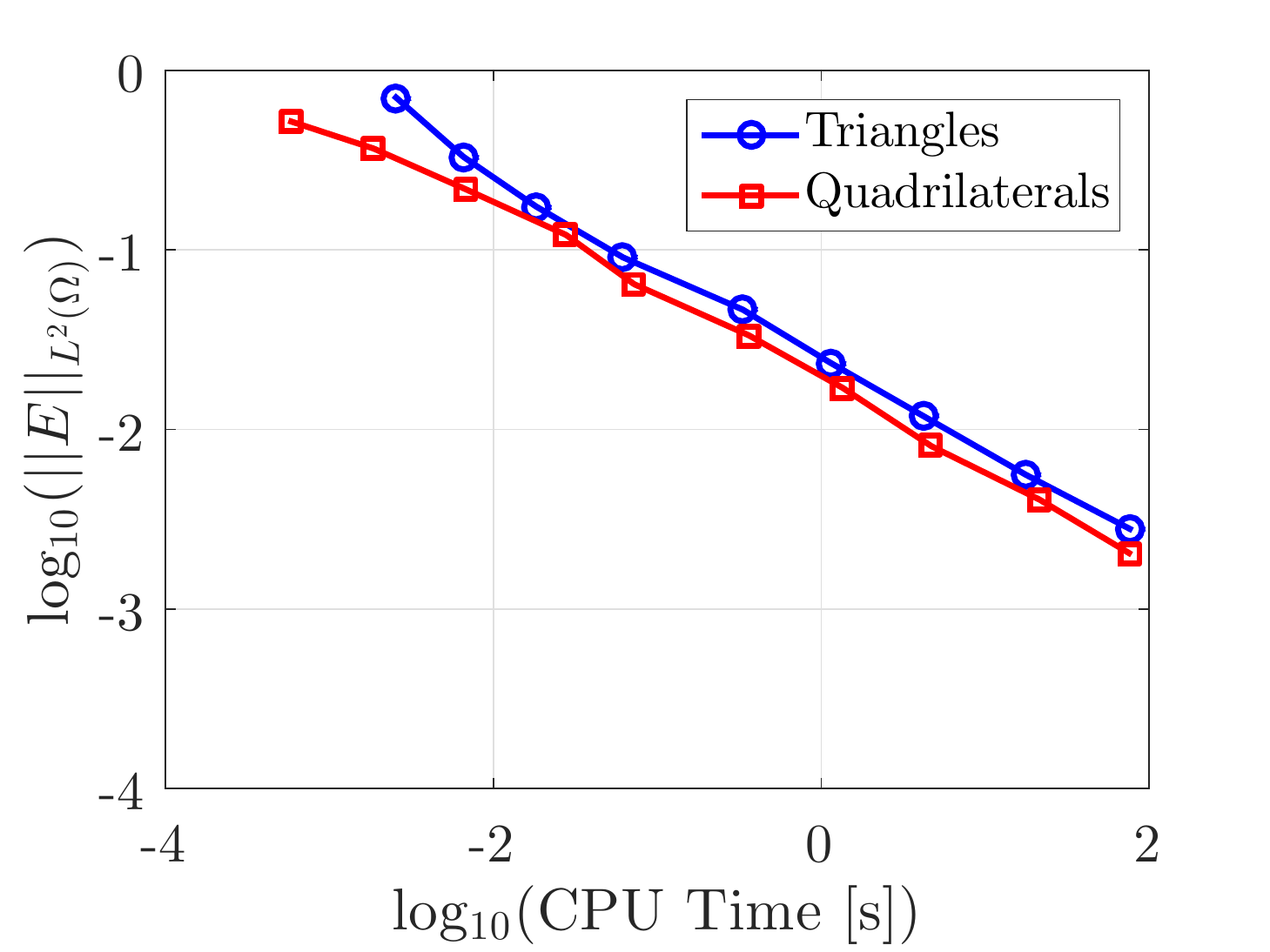}}
	\caption{Error of the solution and its gradient in the $\eltwo(\Omega)$ norm as a function of the CPU time for the 2D Poisson problem.}
	\label{fig:poisson2D_cpu}
\end{figure}

To further study the performance of the FCFV method with triangular and quadrilateral meshes, Figure~\ref{fig:poisson2D_cpu} shows a comparison of triangular and quadrilateral elements in terms of the CPU time. 
The CPU time is measured as the time required to assemble and solve the global system of equations. As described in Sections~\ref{sc:Poisson} and \ref{sc:Stokes}, the cost associated to the global problem is the dominant cost of the proposed FCFV scheme because the solution of the local problems involves the evaluation of an explicit expression and can be performed element-by-element. 
It is clear that the slight superiority in terms of number of degrees of freedom directly translates in a marginal better performance in terms of CPU time. The analysis also indicates the efficiency of the method as solving a problem with more than one million degrees of freedom takes under 100 seconds.

Similarly, for the three dimensional Poisson problem, Figure~\ref{fig:poisson3D_hConv} shows that for a given element size, tetrahedral elements provide the maximum accuracy but this is simply due to the higher number of internal faces (i.e. degrees of freedom of the global problem) compared to meshes of other element types.

\begin{figure}[!tb]
	\centering
	\subfigure[$u$]{\includegraphics[width=0.4\textwidth]{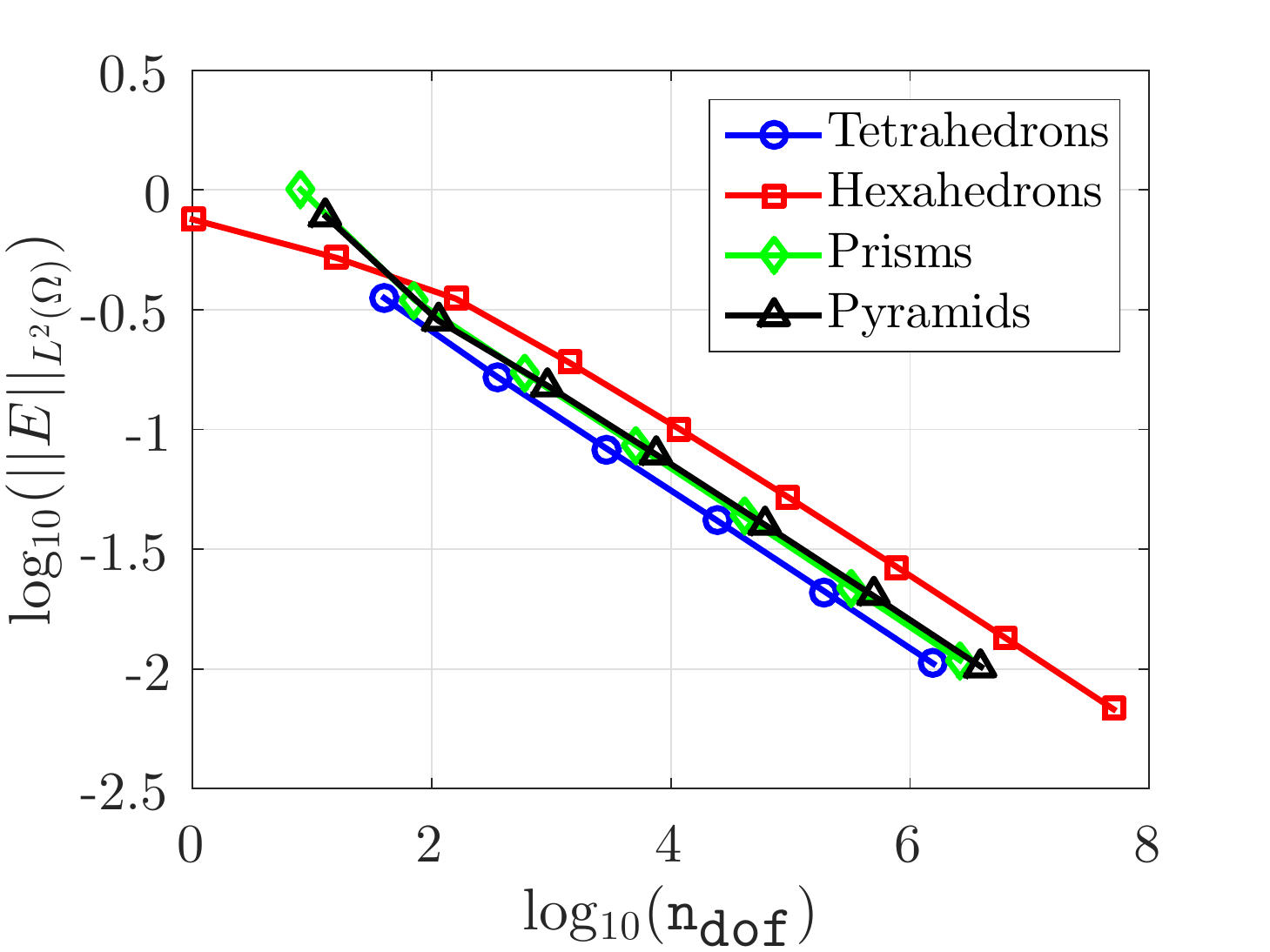}}
	\subfigure[$\bq$]{\includegraphics[width=0.4\textwidth]{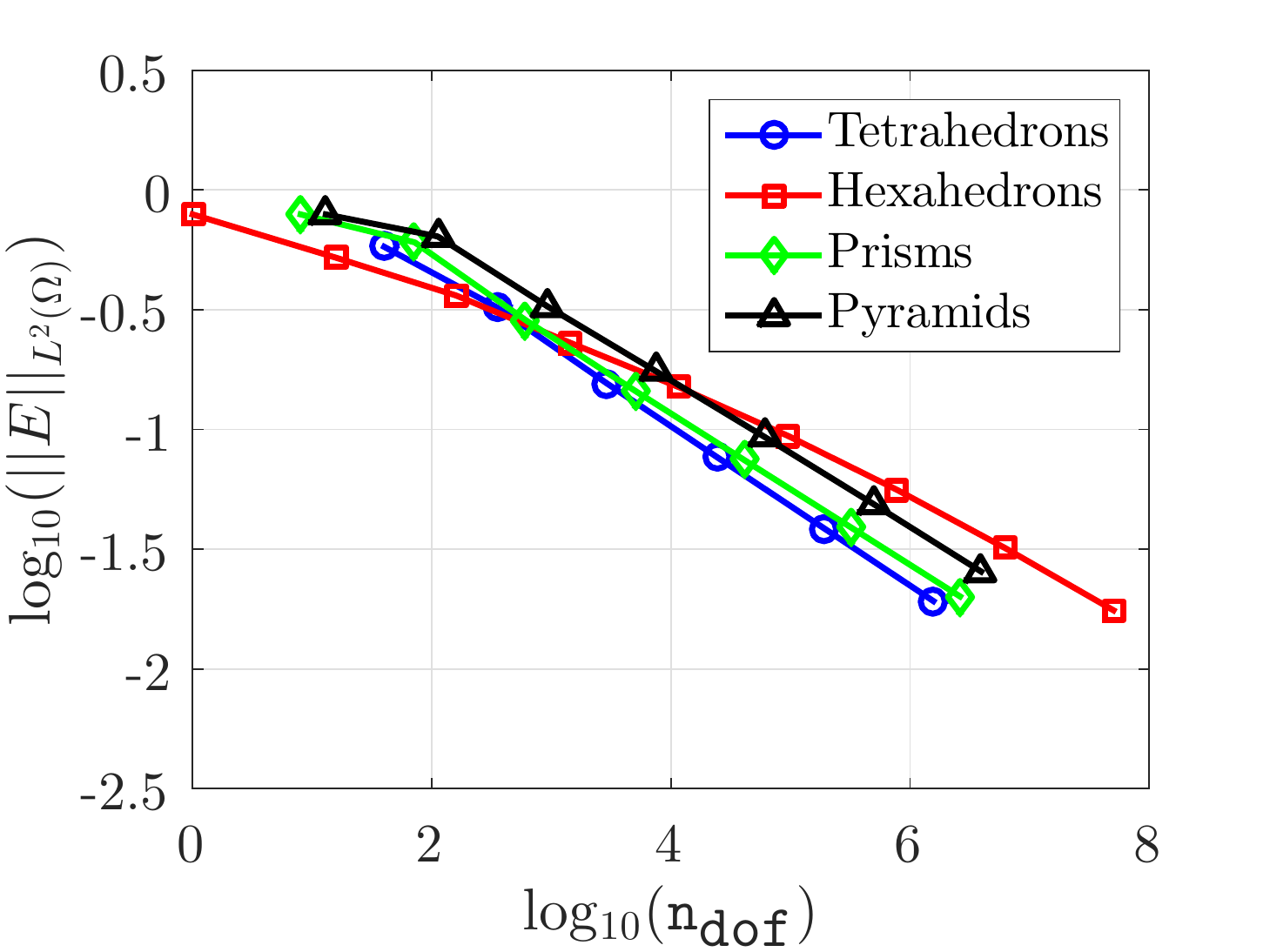}}
	\caption{Error of the solution and its gradient in the $\eltwo(\Omega)$ norm as a function of the number of degrees of freedom of the global system for the 3D Poisson problem.}
	\label{fig:poisson3D_ndof}
\end{figure}

Figure~\ref{fig:poisson3D_ndof} shows the evolution of the error of the primal and dual variables as a function of the number of degrees of freedom of the global problem. 
The results show that tetrahedral elements are able to attain a given error with slightly less degrees of freedom than other element types whereas hexahedral elements require the maximum number of degrees of freedom. It is worth noting that when the error on the dual variable is of interest, the advantages of tetrahedral elements can be better appreciated. For instance, tetrahedral elements provide an error of 0.0192 in the sixth mesh, with 1,564,672 degrees of freedom, whereas hexahedral elements require 50,200,576 in order to provide a similar error, 0.0176 in this example.
	
\begin{figure}[!tb]
	\centering
	\subfigure[$u$]{\includegraphics[width=0.4\textwidth]{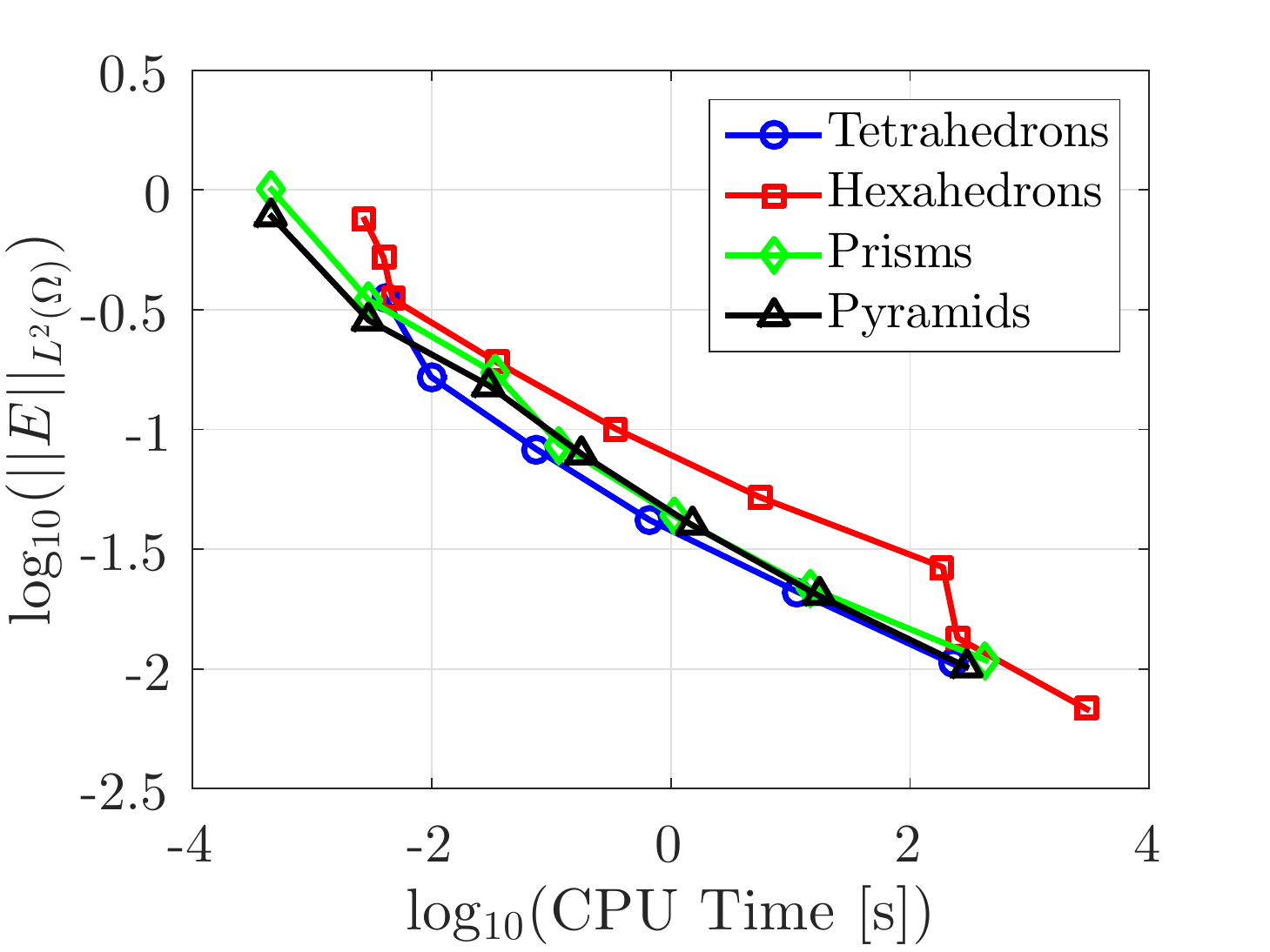}}
	\subfigure[$\bq$]{\includegraphics[width=0.4\textwidth]{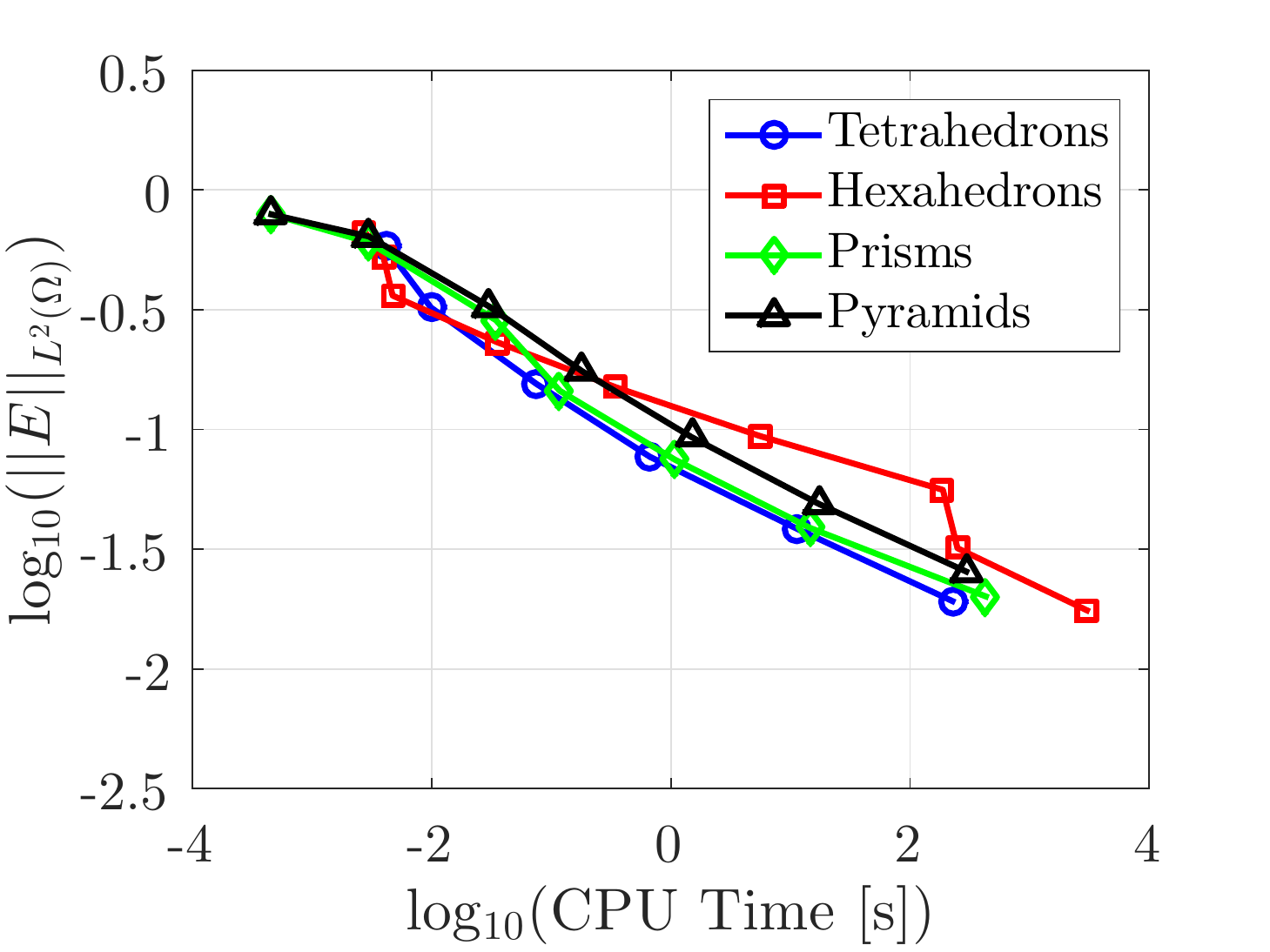}}
	\caption{Error of the solution and its gradient in the $\eltwo(\Omega)$ norm as a function of the CPU time for the 3D Poisson problem.}
	\label{fig:poisson3D_cpu}
\end{figure}

To study if the reduction of degrees of freedom also translates in more efficient computations, Figure~\ref{fig:poisson3D_cpu} shows a comparison of the different types of elements in terms of the CPU time. 
The results show that a similar performance is obtained by the FCFV scheme with tetrahedrons, prisms and pyramids. The worst performance is observed for hexahedral elements, especially if a low error in the dual variable is required.

\begin{figure}[!tb]
	\centering
	\subfigure[$p$]{\includegraphics[width=0.32\textwidth]{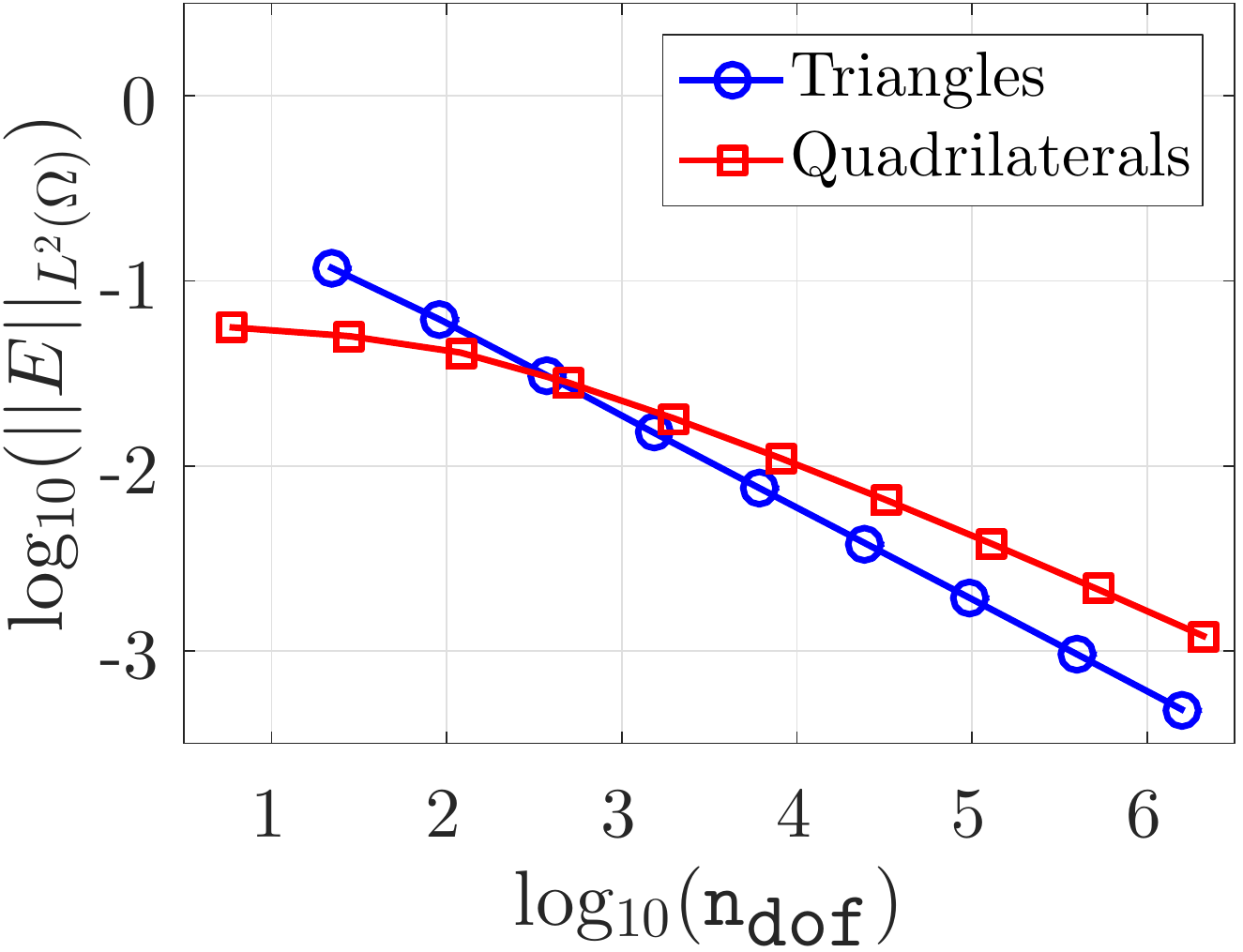}}
	\subfigure[$\bu$]{\includegraphics[width=0.32\textwidth]{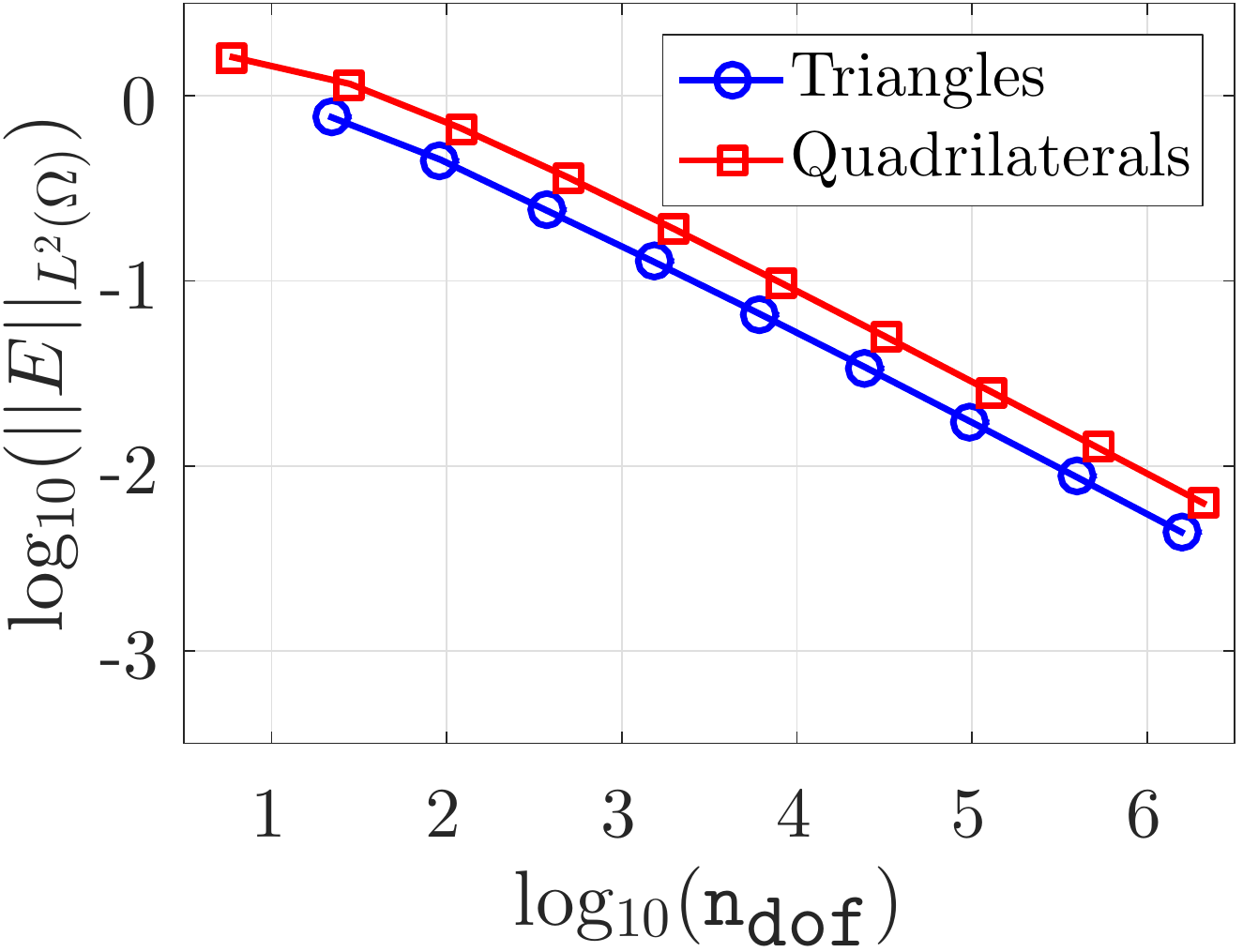}}
	\subfigure[$\bm{L}$]{\includegraphics[width=0.32\textwidth]{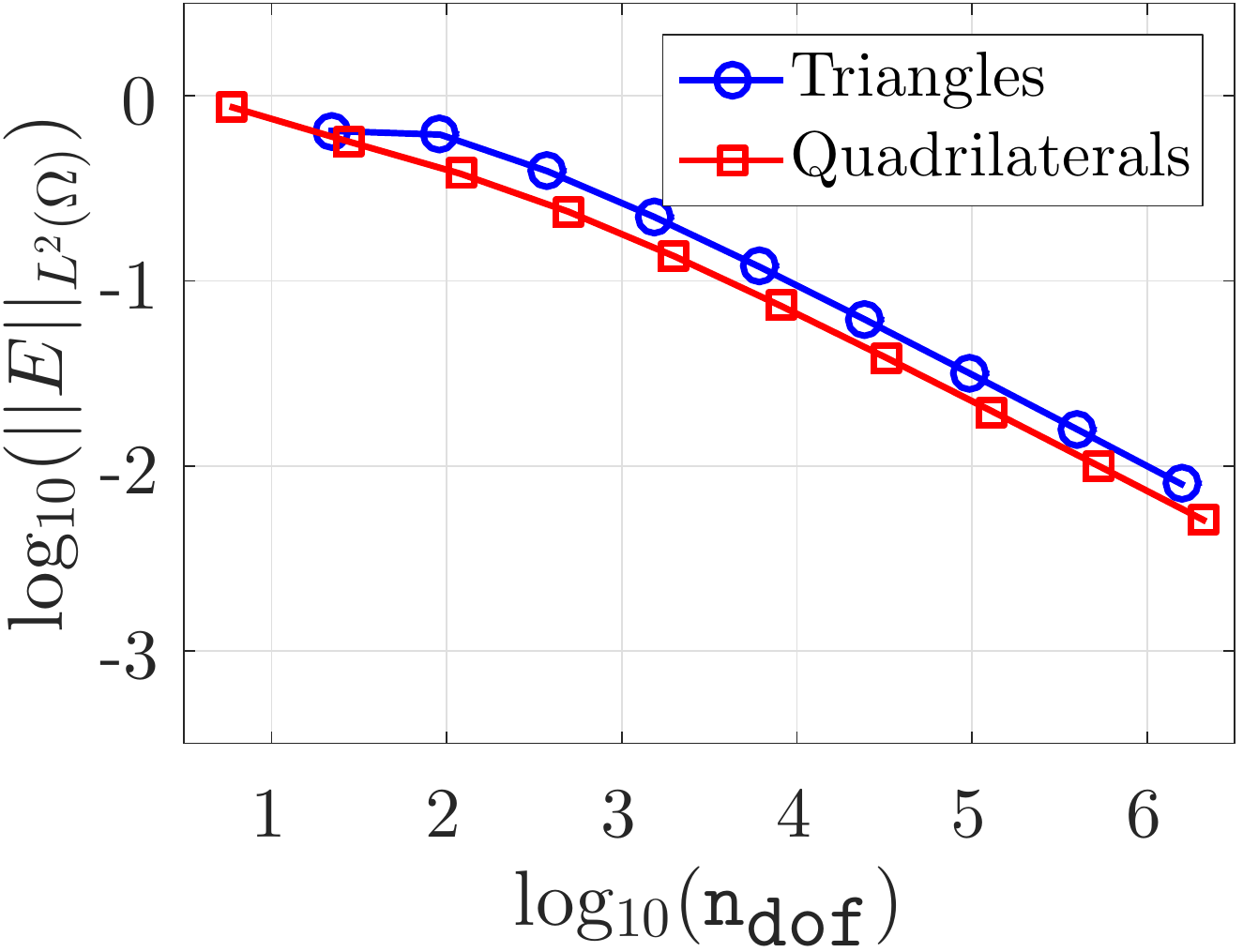}}
	\caption{Error of the pressure, the velocity and its gradient in the $\eltwo(\Omega)$ norm as a function of the number of degrees of freedom of the global problem for the 2D Stokes problem.}
	\label{fig:stokes2D_ndof}
\end{figure}

To study if the same conclusions are obtained for other problems, a similar analysis is performed for the solution of the Stokes equation in two dimensions. Figure~\ref{fig:stokes2D_ndof} shows the evolution of the error of the pressure, velocity and velocity gradient, measured in the $\eltwo(\Omega)$ norm, as a function of the number of degrees of freedom.
It is worth remarking that for the Stokes problem the size of the global system of equations corresponds to the number of internal and Neumann faces times the number of spatial dimensions plus the total number of elements, as described in Section~\ref{sc:StokesFV}.
The results show that triangular elements offer an extra accuracy compared to quadrilateral elements when the error of the pressure and the velocity is considered, whereas quadrilateral elements provide a slightly more accurate representation of the gradient of the velocity, corroborating the conclusions obtained for the Poisson problem. The more sizeable difference is observed in the pressure when fine meshes are considered.

\begin{figure}[!tb]
	\centering
	\subfigure[$p$]{\includegraphics[width=0.32\textwidth]{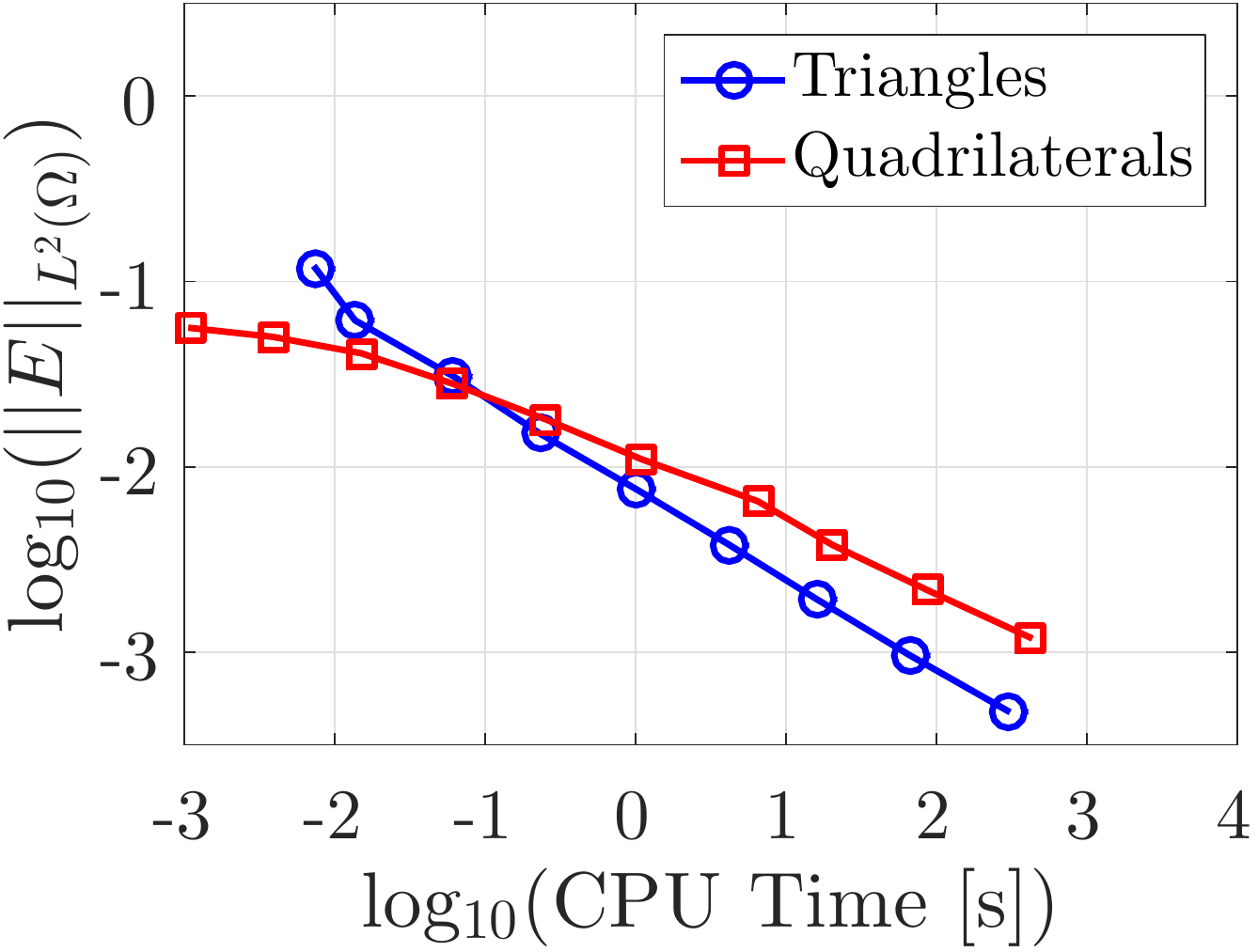}}
	\subfigure[$\bu$]{\includegraphics[width=0.32\textwidth]{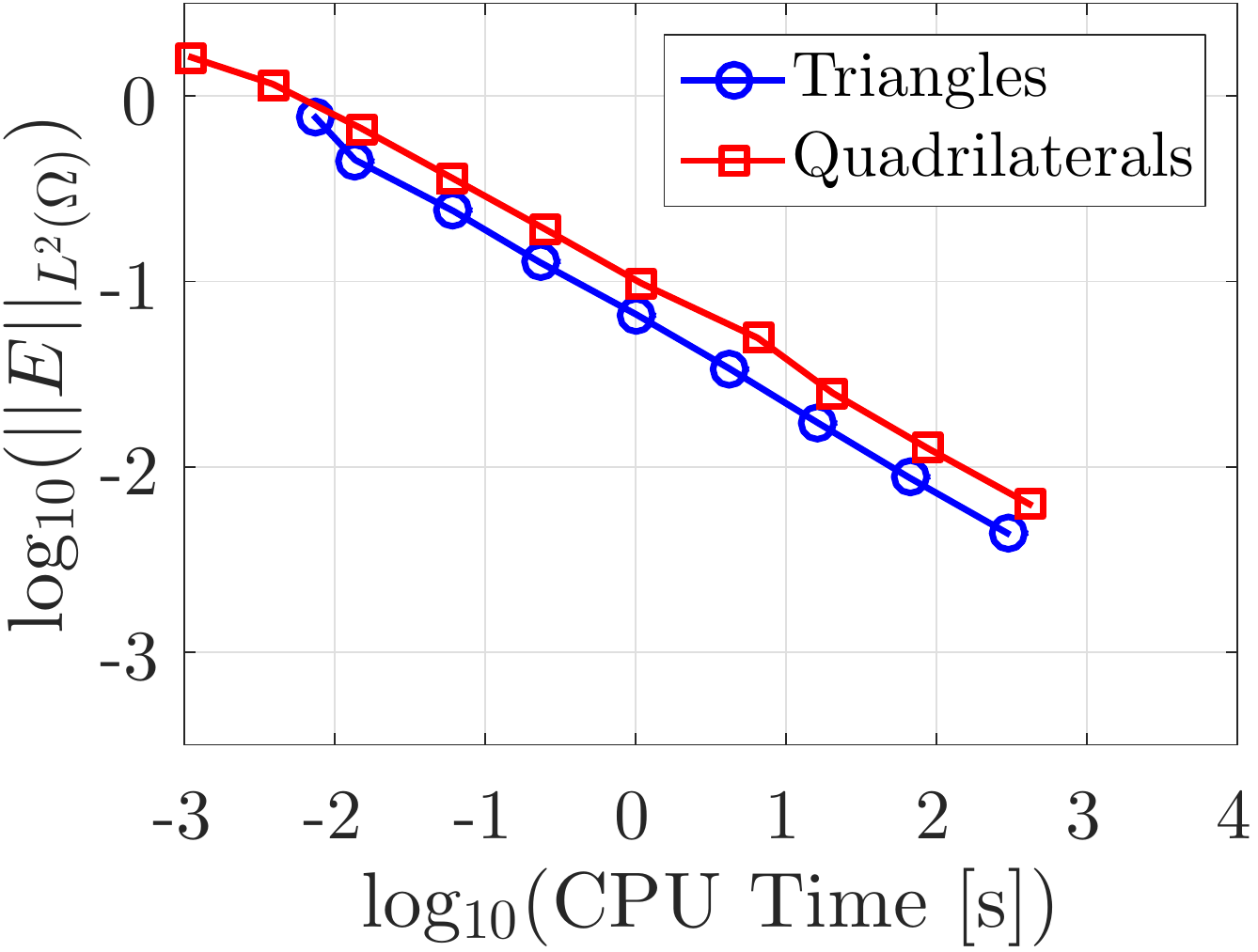}}
	\subfigure[$\bm{L}$]{\includegraphics[width=0.32\textwidth]{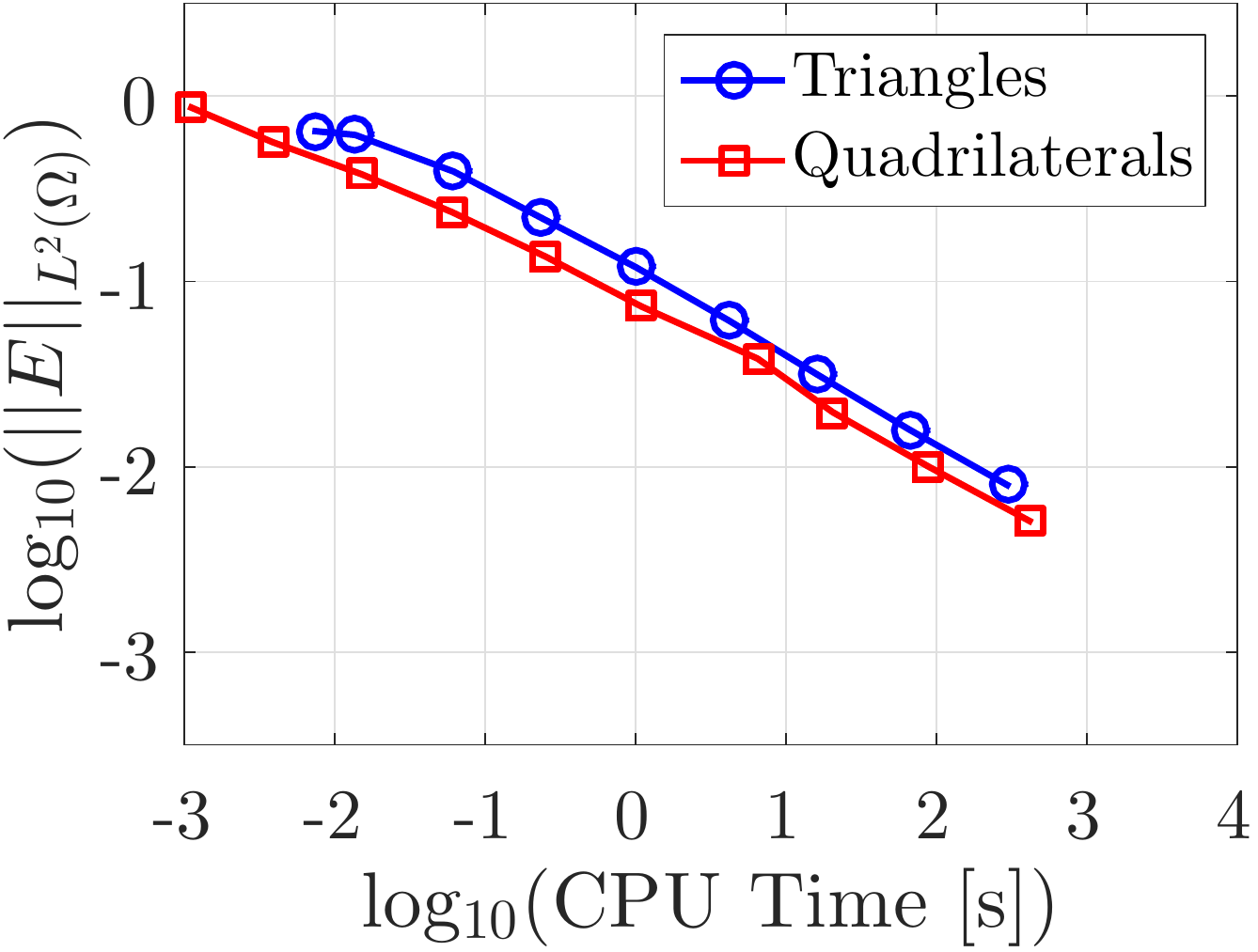}}
	\caption{Error of the pressure, the velocity and its gradient in the $\eltwo(\Omega)$ norm as a function of the CPU time for the 2D Stokes problem.}
	\label{fig:stokes2D_cpu}
\end{figure}

The comparison in terms of CPU time is shown in Figure~\ref{fig:stokes2D_cpu}.
The study shows, once more, that the slight superiority observed in terms of number of degrees of freedom also translates in a marginal better performance in terms of CPU time. 
It is important to note that the conclusions obtained for the Stokes problem are consistent with the observations made for the Poisson problem. 

\subsection{Influence of the stabilization parameter}
\label{sc:InfluenceTau}

In all the examples shown in previous Sections, the stabilization parameter $\tau$ has been fixed to a constant value to perform the mesh convergence studies. 	Next, the influence of the stabilization parameter $\tau$ is studied for the solution of Poisson and Stokes problems and for different element types.

\begin{figure}[!tb]
	\centering
	\subfigure[Tetrahedrons]{\includegraphics[width=0.4\textwidth]{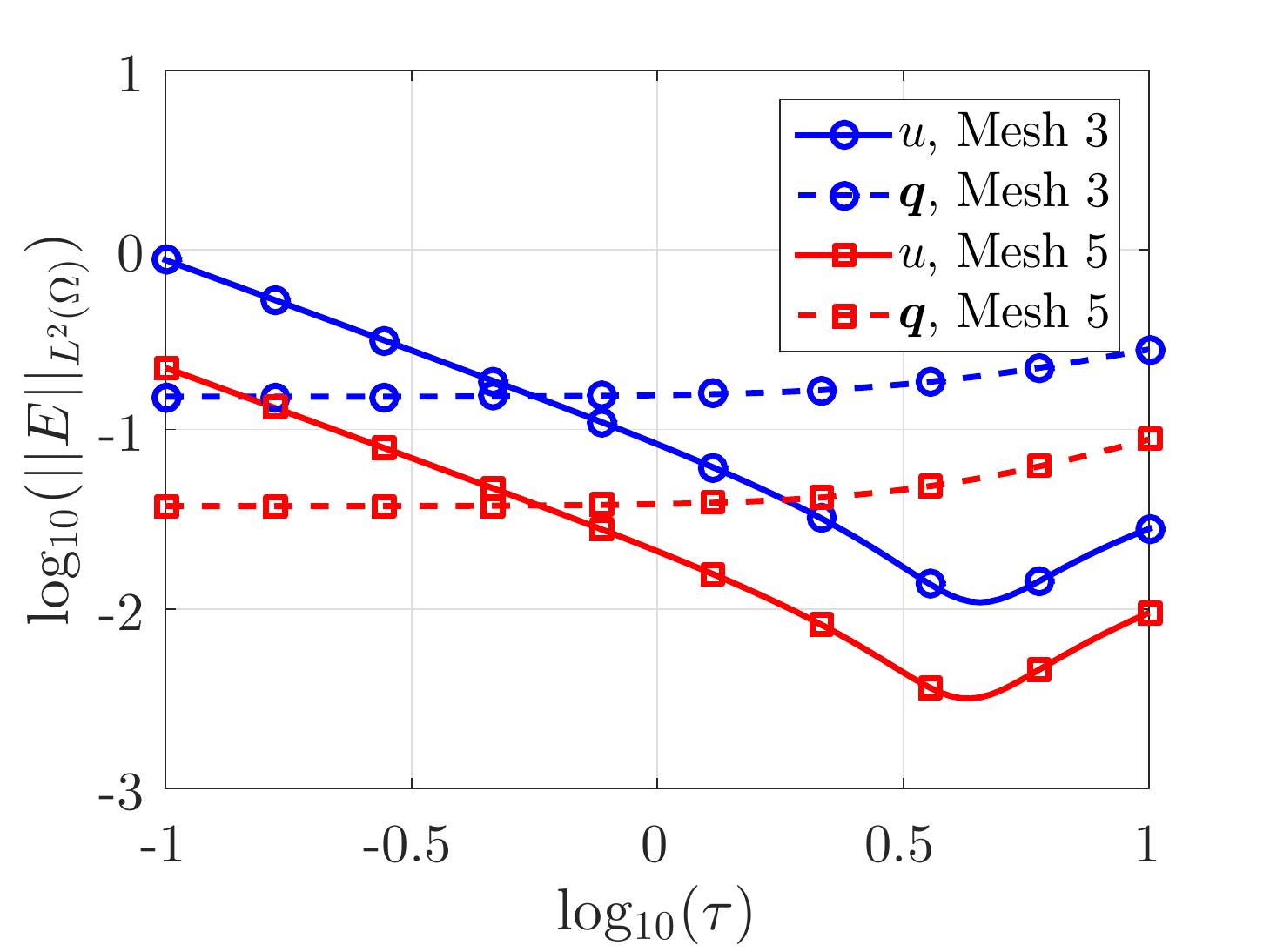}}
	\subfigure[Hexahedrons]{\includegraphics[width=0.4\textwidth]{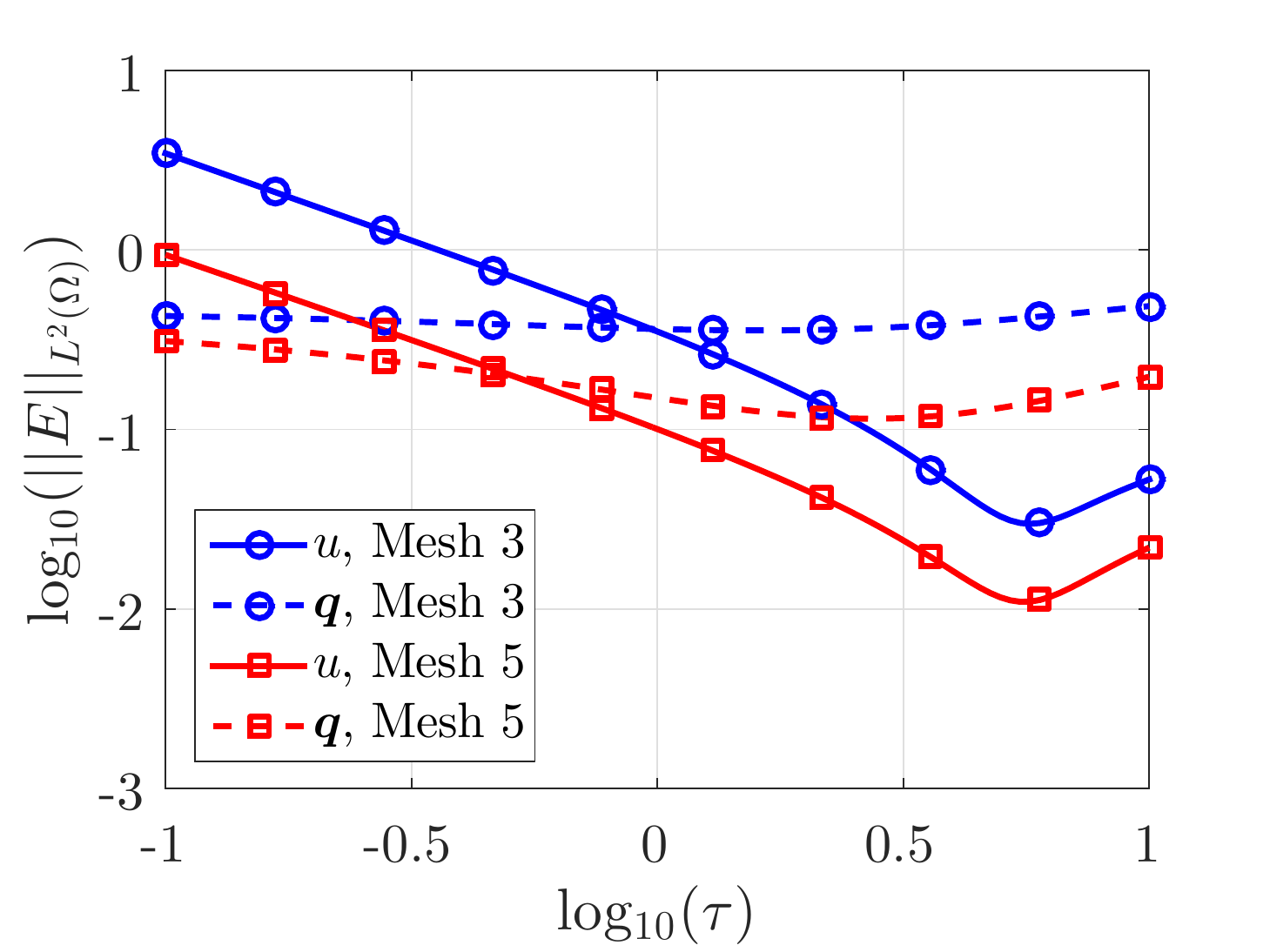}}
	\subfigure[Prisms]{\includegraphics[width=0.4\textwidth]{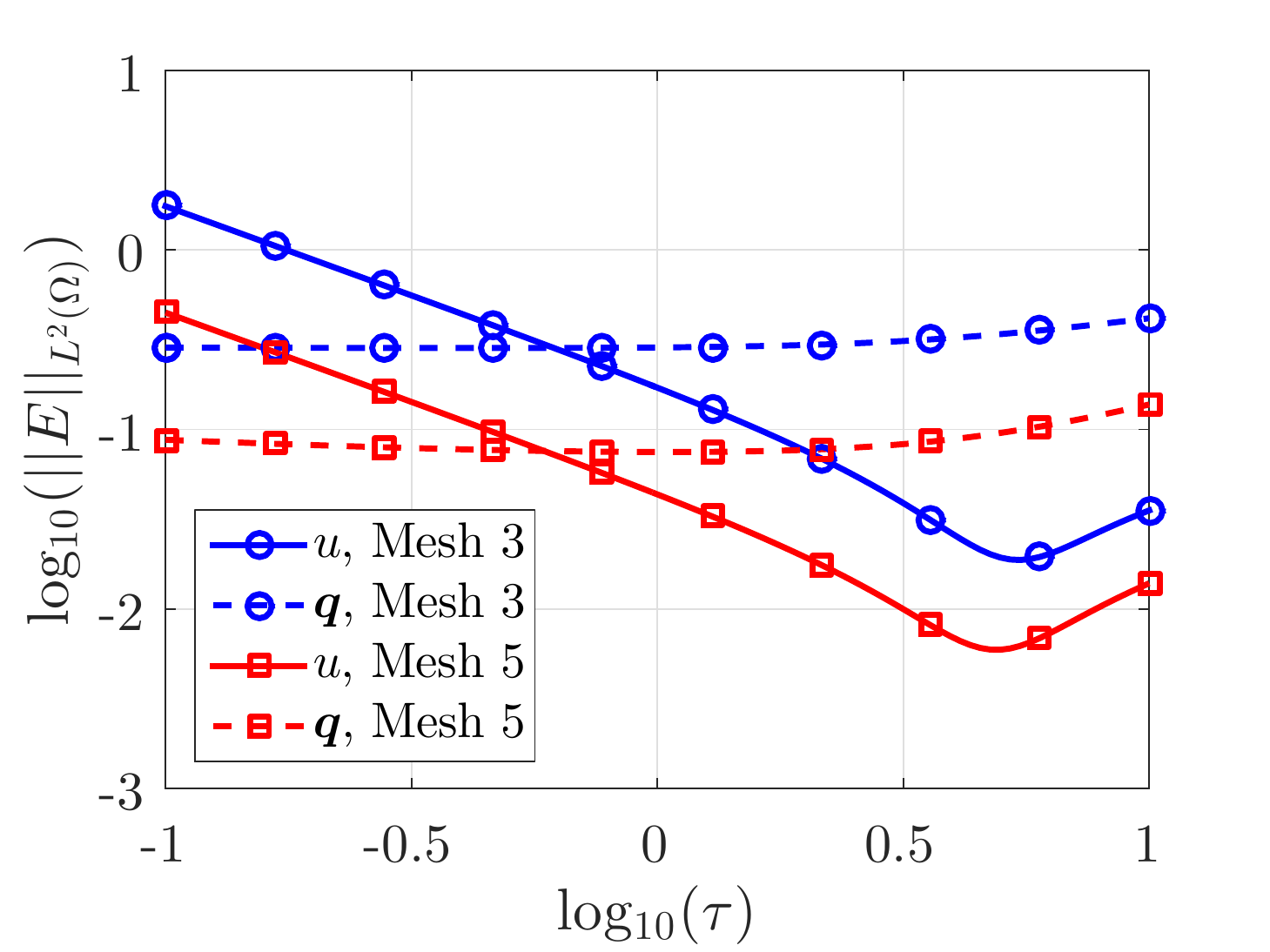}}
	\subfigure[Pyramids]{\includegraphics[width=0.4\textwidth]{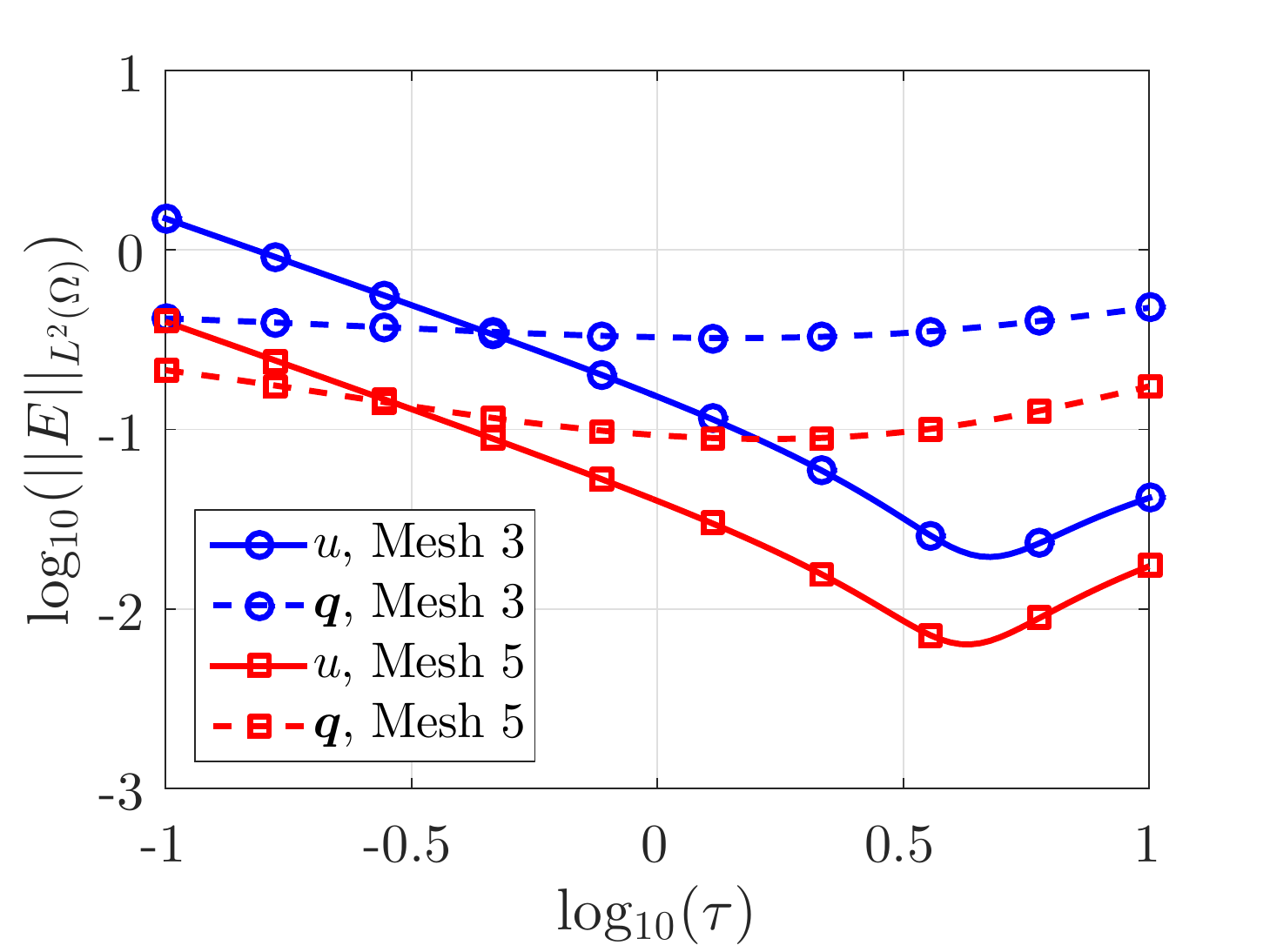}}
	\caption{Error of the solution and its gradient in the $\eltwo(\Omega)$ norm as a function of the stabilization parameter $\tau$ for a 3D Poisson problem.}
	\label{fig:poisson3DTau}
\end{figure}

Figure~\ref{fig:poisson3DTau} shows the evolution of the error of the solution and its gradient in the $\eltwo(\Omega)$ norm as a function of the stabilization parameter $\tau$ for a three dimensional Poisson problem using four different element types. The simulation is performed using two different meshes and the value of $\tau$ varies from 0.1 to 10. 
%
%
The results suggest that there is a value of $\tau$ that provides the minimum error on the primal solution. This value seems to be independent on the level of mesh refinement and only slightly dependent on the type of element. In this example, the FCFV scheme provides the minimum error in the primal solution for a value of $\tau \approx 3$. The same conclusions are also obtained for the Poisson problem in two dimensions (the results are not presented for brevity).

It can be observed that the error of the gradient of the solution $\bm{q}$ is less sensitive to the value of the stabilization parameter. For triangular, tetrahedral, prisms and pyramids, a value of $\tau$ between 0.1 and 3 does not induce a significant variation on the accuracy whereas a higher value induces an increase in the error of $\bm{q}$. For quadrilateral and hexahedral elements the minimum error is obtained for the value $\tau \approx 3$ whereas lower or higher values induce a loss of accuracy in $\bm{q}$.

\begin{figure}[!tb]
	\centering
	\subfigure[Triangles]{\includegraphics[width=0.4\textwidth]{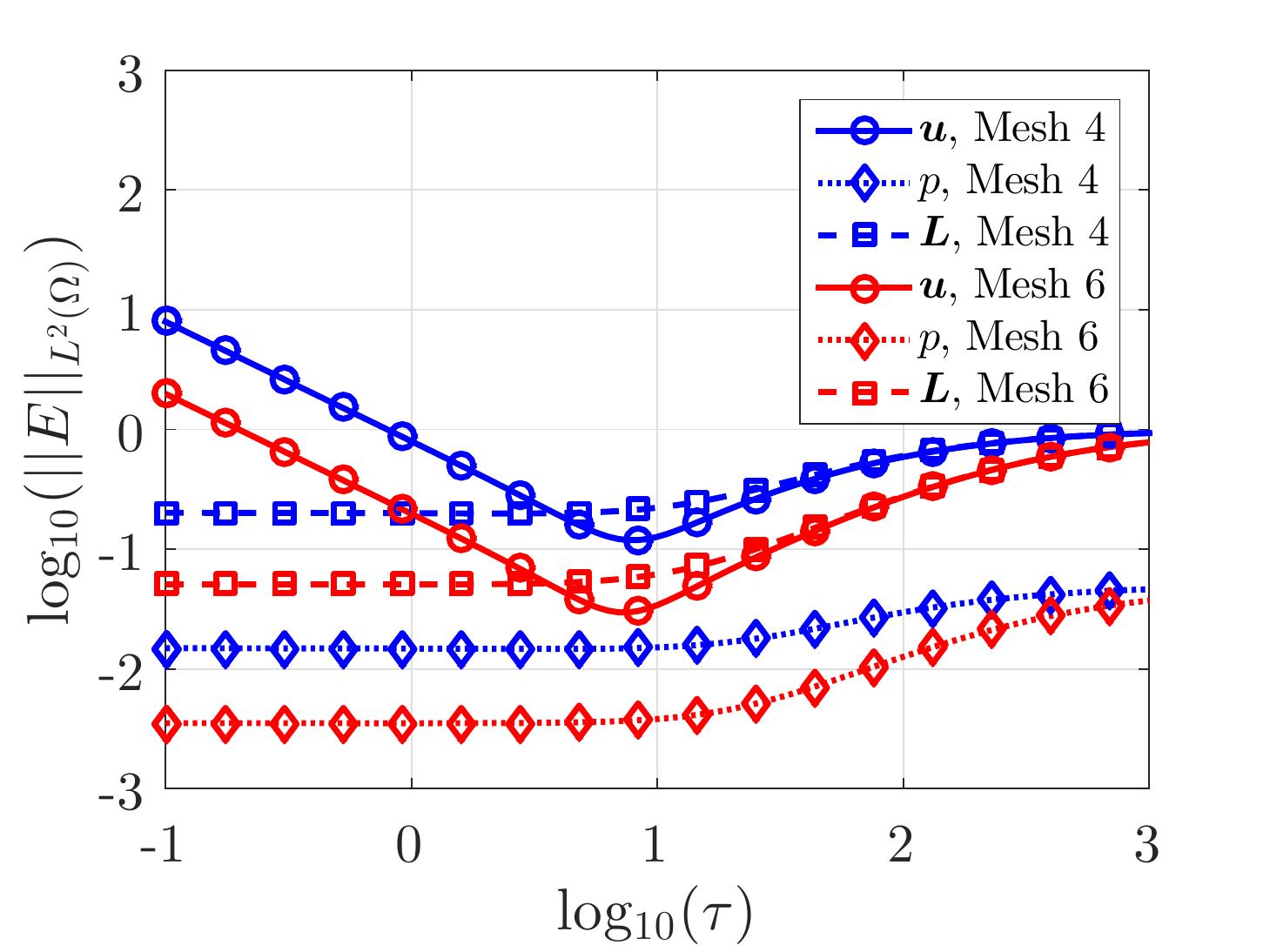}}
	\subfigure[Quadrilaterals]{\includegraphics[width=0.4\textwidth]{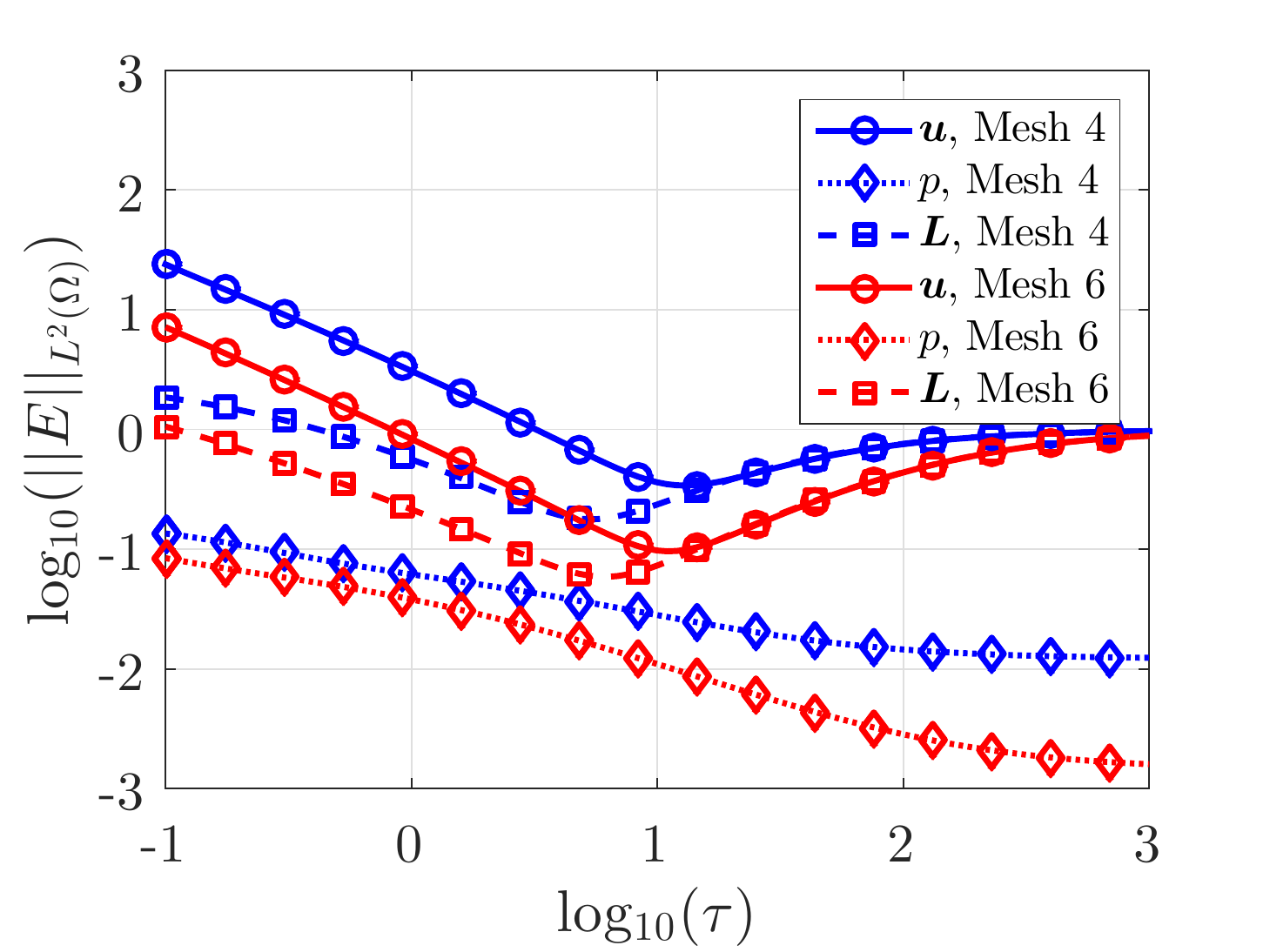}}
	\caption{Error of the solution and its gradient in the $\eltwo(\Omega)$ norm as a function of the stabilization parameter $\tau$ for a 2D Stokes problem.}
	\label{fig:stokes2DTau}
\end{figure}

A similar study has also been performed for the Stokes problem. The results represented in Figure~\ref{fig:stokes2DTau} correspond to a two dimensional Stokes problem and they show that the influence of $\tau$ is similar for both Poisson and Stokes problems. 
For both triangular and quadrilateral meshes a value of $\tau \approx 10$ provides the maximum accuracy for the velocity. When the error on the velocity gradient is of interest, the conclusions are identical to the ones discussed for the Poisson problem. In terms of the error in pressure, a different behaviour is observed for triangular and quadrilateral meshes. For triangular meshes, the accuracy on the pressure is not affected by the stabilization parameter if a value of $\tau \leq 10$ is selected. For higher values of $\tau$, the error in pressure increases with a significant impact in the accuracy. With quadrilateral elements the behaviour of the error on the pressure is different as higher values of the stabilisation parameter provide a lower error.
It is interesting to observe that the accuracy obtained for the velocity and its gradient in two different meshes is almost identical when a large value of the stabilization parameter is considered (e.g. $\tau \approx 1,000$). This means that all the error is controlled by the value of $\tau$ and not by the level of mesh refinement. This behaviour is attributed to the definition of the numerical fluxes in Equation~\eqref{eq:traceStokes}. A large value of $\tau$ implies that the numerical normal flux receives a negligible contribution from the physical normal flux.

\subsection{Influence of the element distortion}
\label{sc:InfluenceDistortion}

\begin{figure}[!tb]
	\centering
	\subfigure[Mesh 3]{\includegraphics[width=0.24\textwidth]{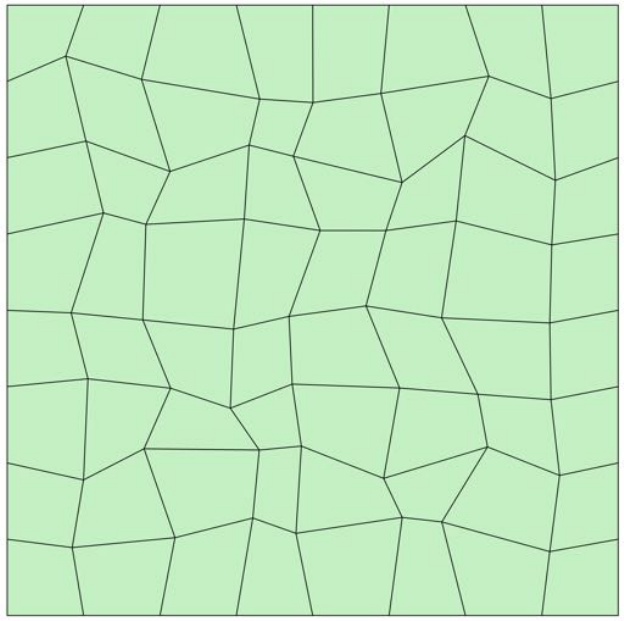}}
	\subfigure[Mesh 5]{\includegraphics[width=0.24\textwidth]{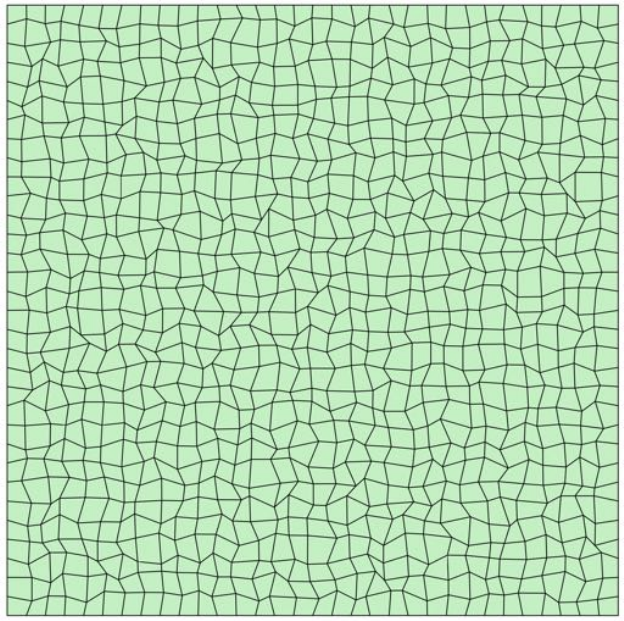}}
	\subfigure[Mesh 3]{\includegraphics[width=0.24\textwidth]{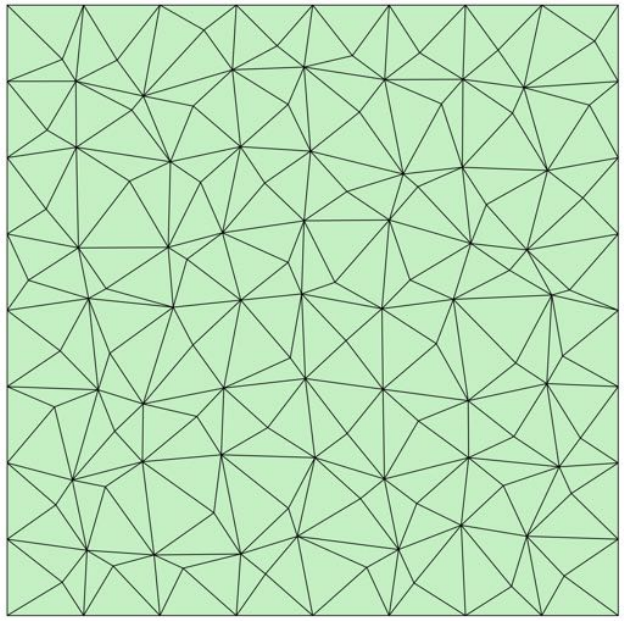}}
	\subfigure[Mesh 5]{\includegraphics[width=0.24\textwidth]{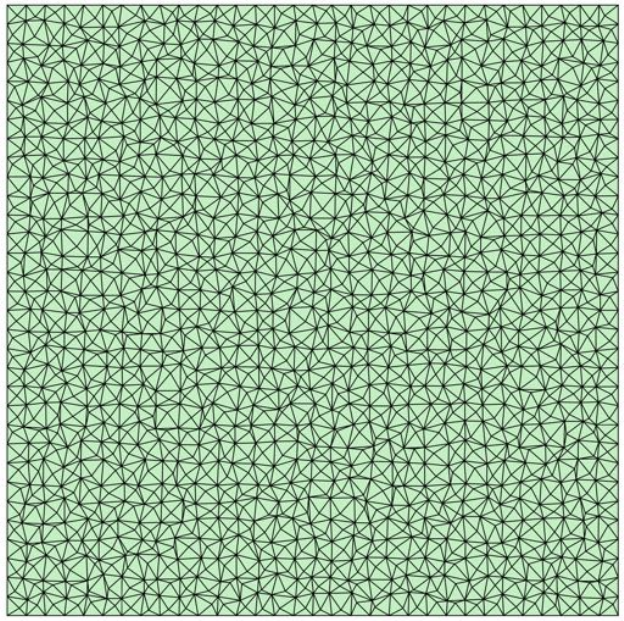}}
	\caption{Irregular quadrilateral and triangular meshes of $\Omega=[0,1]^2$.}
	\label{fig:poisson2D_meshesIrregular}
\end{figure}

The next numerical study involves exploring the influence of the element distortion on the accuracy of the proposed FCFV scheme. To this end, a new set of meshes is produced by perturbing the position of the interior nodes of the regular meshes employed in Sections~\ref{sc:PoissonVerification} and \ref{sc:StokesVerification}. In all cases the new position of the $i$-th node is simply defined as $\tilde{\bx}_i = \bx_i + \bm{r}_i$, where $\bm{r}_i$ is a vector of dimension $\nsd$ where each component is a randomly generated number within the interval $[-\ell_{\text{min}}/3,\ell_{\text{min}}/3]$ and $\ell_{\text{min}}$ denote the minimum edge length of the regular mesh.
Two of the resulting irregular quadrilateral and triangular meshes are represented in Figure~\ref{fig:poisson2D_meshesIrregular}.	

\begin{figure}[!tb]
	\centering
	\subfigure[$u$]{\includegraphics[width=0.4\textwidth]{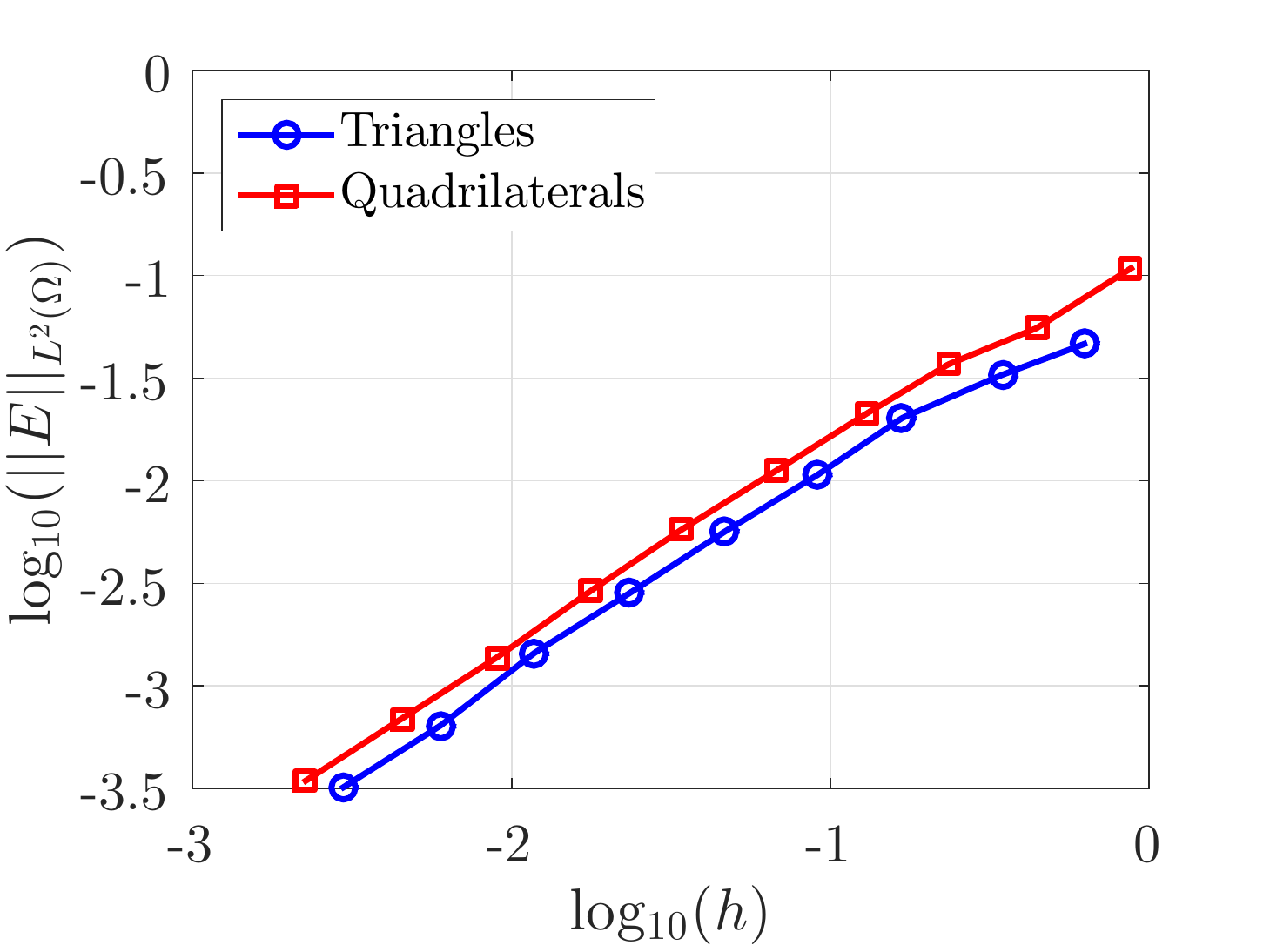}}
	\subfigure[$\bq$]{\includegraphics[width=0.4\textwidth]{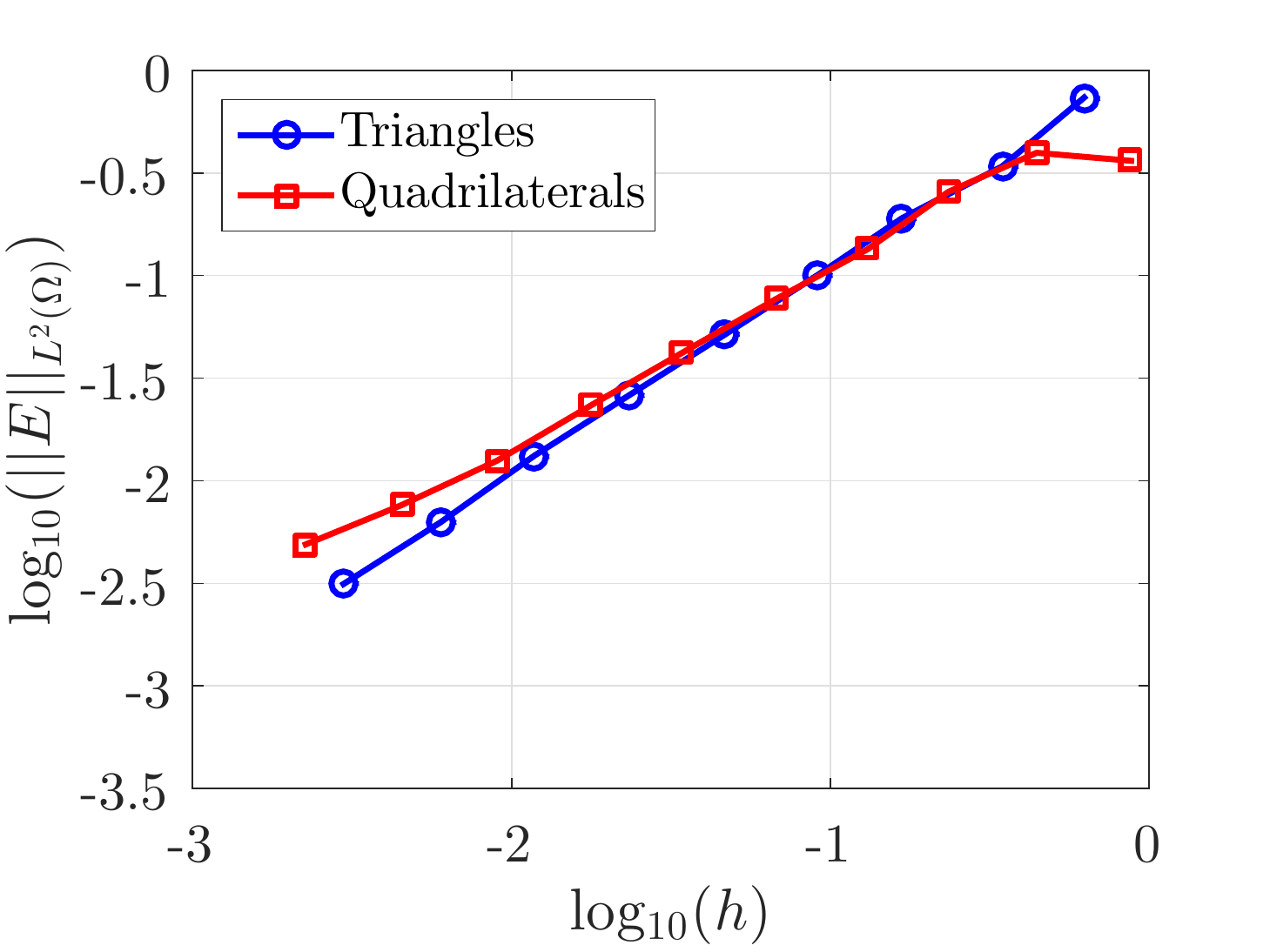}}
	\caption{Mesh convergence of the error of the solution and its gradient in the $\eltwo(\Omega)$ norm for the 2D Poisson problem with irregular meshes.}
	\label{fig:poisson2D_hConvIrregular}
\end{figure}

The Poisson solver is considered first. The convergence of the error of the primal variable $u$, measured in the $\eltwo(\Omega)$ norm, as a function of the characteristic element size $h$ is depicted in Figure~\ref{fig:poisson2D_hConvIrregular} (a) for both triangular and quadrilateral elements. Similarly, the convergence of the error of the dual variable $\bq$, measured in the $\eltwo(\Omega)$ norm, as a function of the characteristic element size $h$ is depicted in Figure~\ref{fig:poisson2D_hConvIrregular} (b) for both triangular and quadrilateral elements. 
The results show the optimal order of convergence in the primal variable for both elements. For the dual variable, triangular elements show an optimal rate of convergence whereas a slight loss of accuracy is obtained for quadrilateral elements in finer meshes. By comparing the results obtained in irregular meshes with the results for regular meshes presented in Section~\ref{sc:PoissonVerification}, it can be concluded that the accuracy of the FCFV method is not heavily dependent upon the distortion of the elements.

\begin{figure}[!tb]
	\centering
	\subfigure[$p$]{\includegraphics[width=0.32\textwidth]{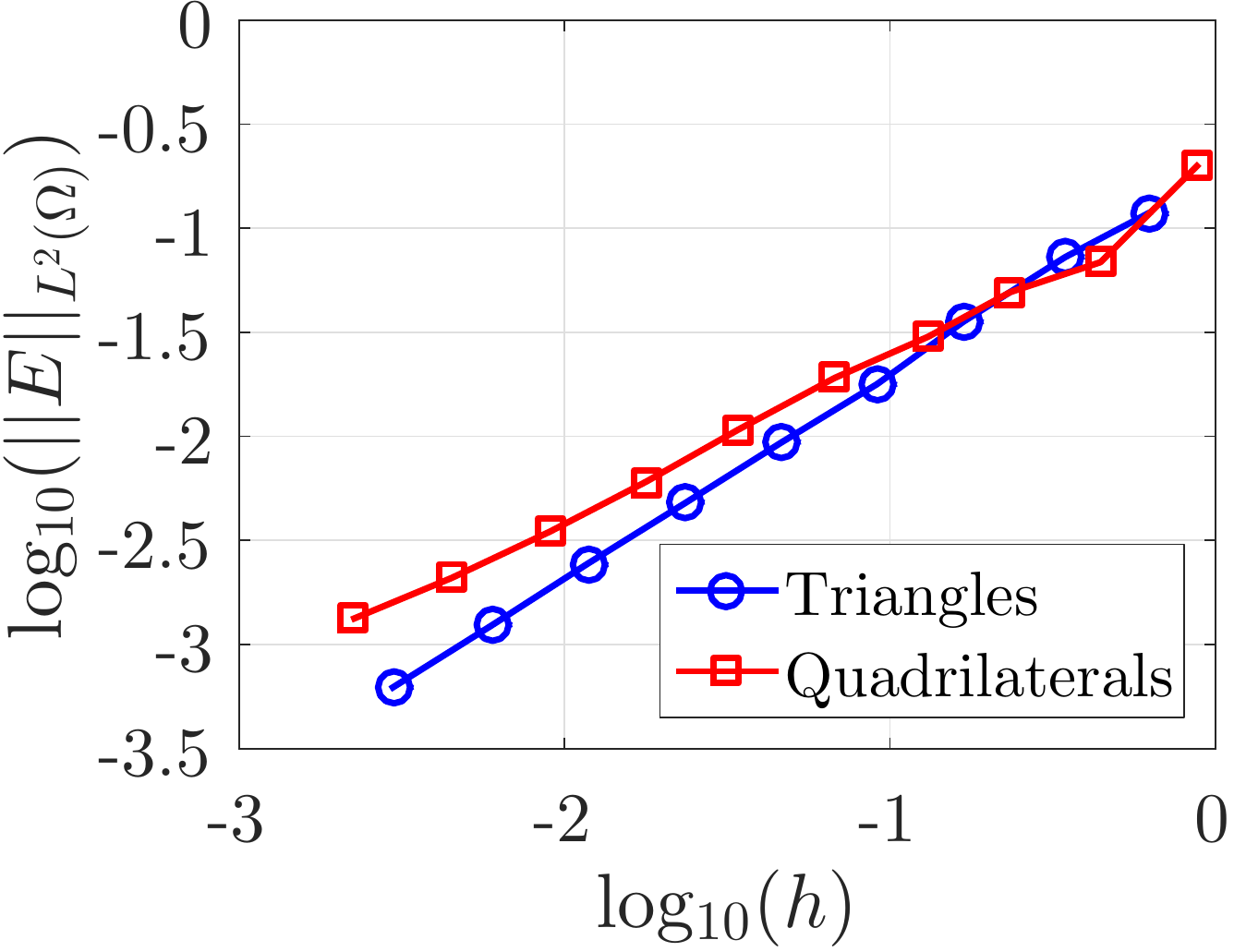}}
	\subfigure[$\bu$]{\includegraphics[width=0.32\textwidth]{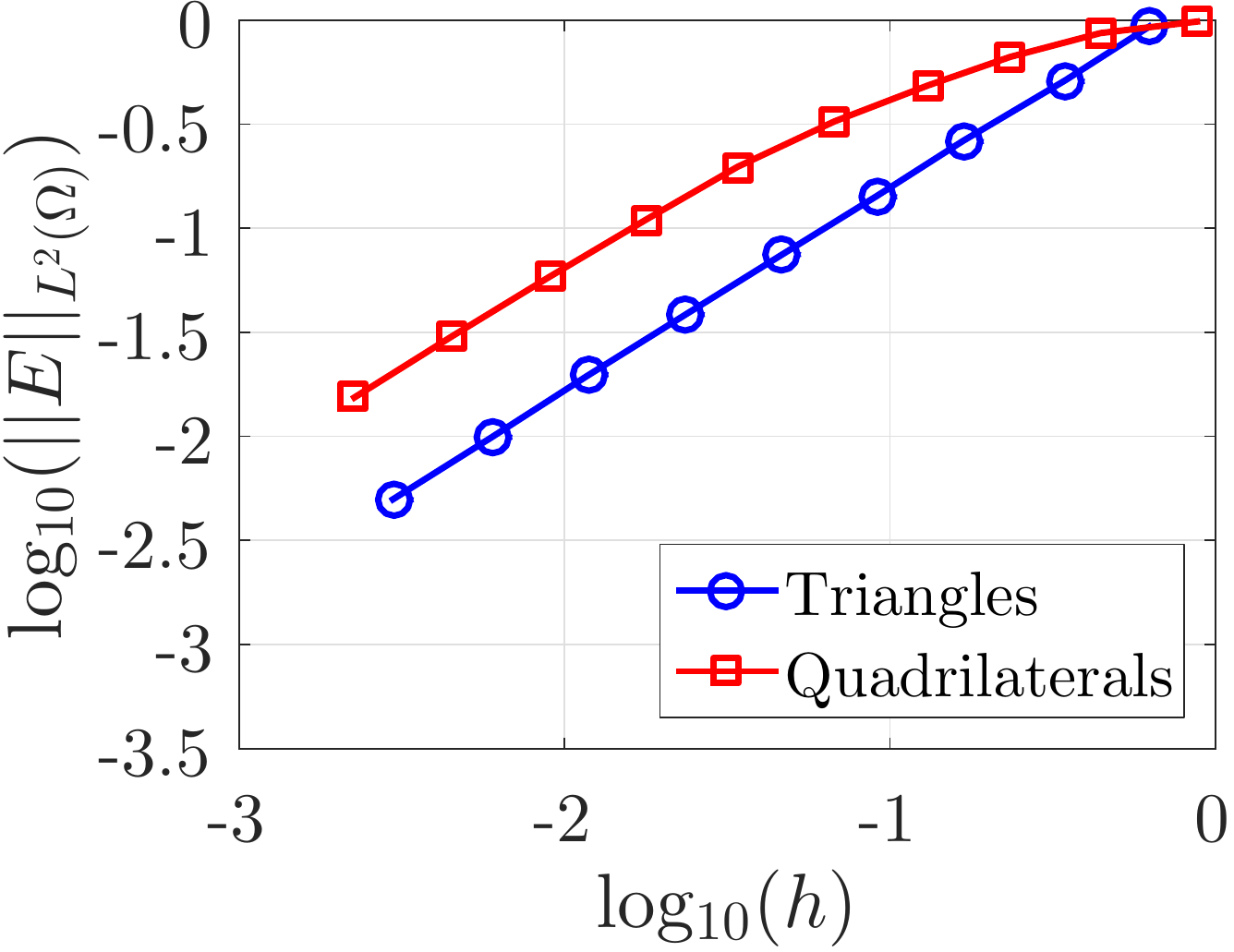}}
	\subfigure[$\bm{L}$]{\includegraphics[width=0.32\textwidth]{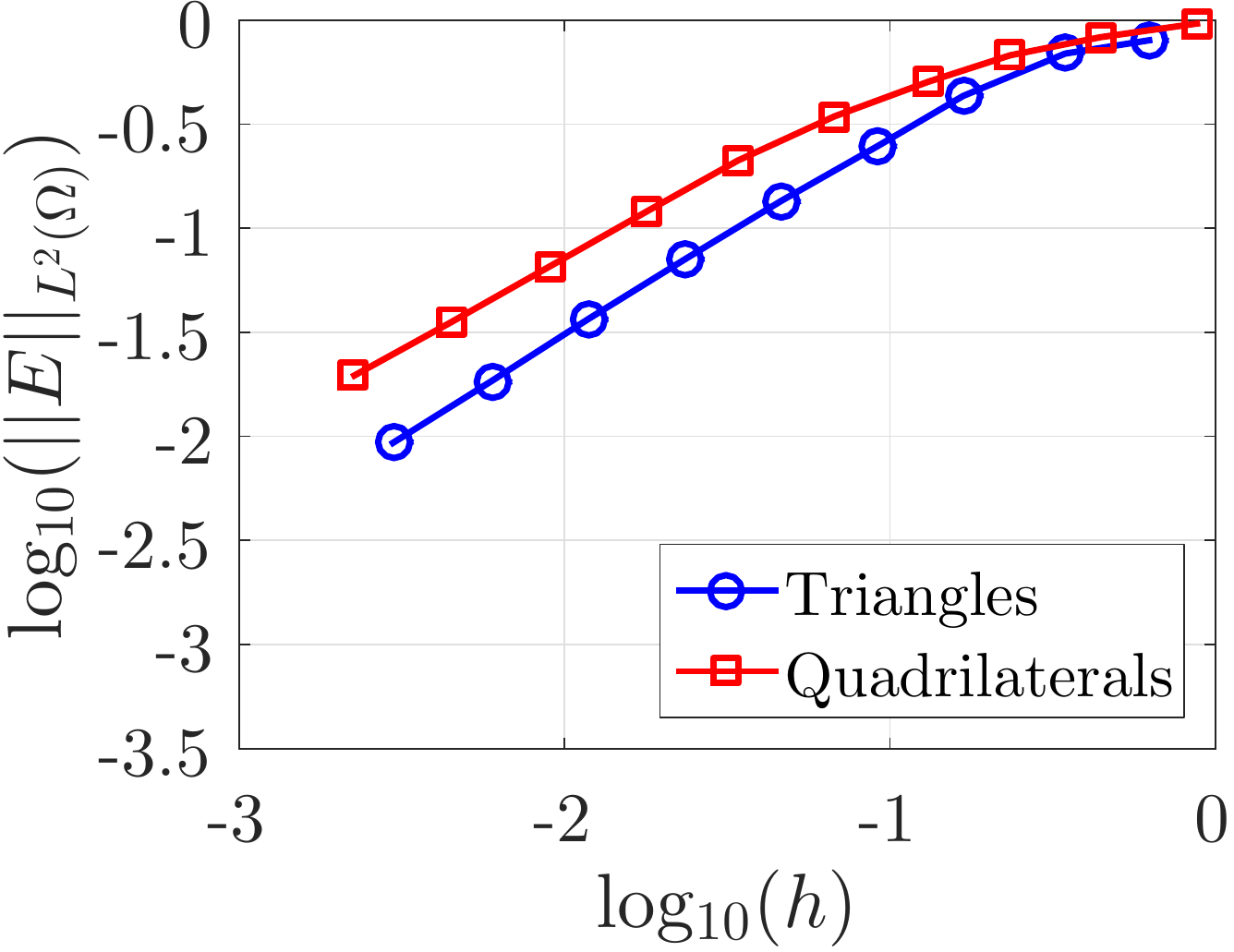}}
	\caption{Mesh convergence of the error of the pressure, the velocity and the velocity gradient in the $\eltwo(\Omega)$ norm for the 2D Stokes with irregular meshes.}
	\label{fig:stokes2D_hConvIrregular}
\end{figure}

The results of the mesh convergence study for the two dimensional Stokes solver are presented in Figure~\ref{fig:stokes2D_hConvIrregular}.
It can be observed that all variables converge with the optimal rate of convergence with triangular elements whereas the error in the pressure suffers a slight loss of accuracy for finer quadrilateral meshes.

\begin{figure}[!tb]
	\centering
	\subfigure[Hexahedrons]{\includegraphics[width=0.24\textwidth]{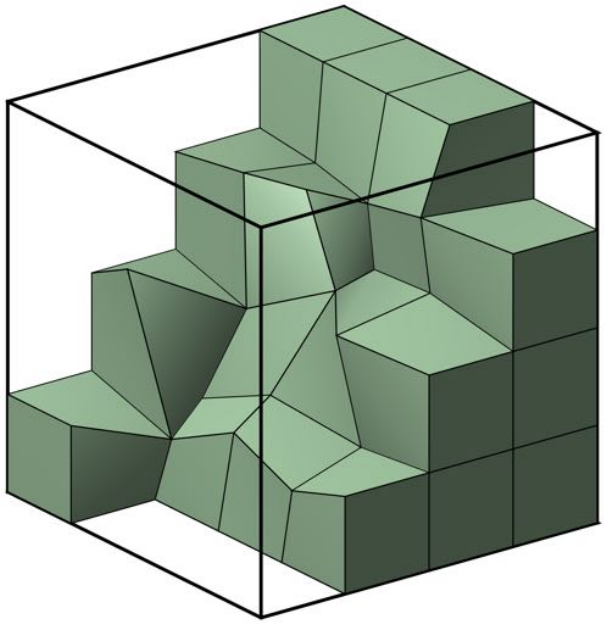}}
	\subfigure[Tetrahedrons]{\includegraphics[width=0.24\textwidth]{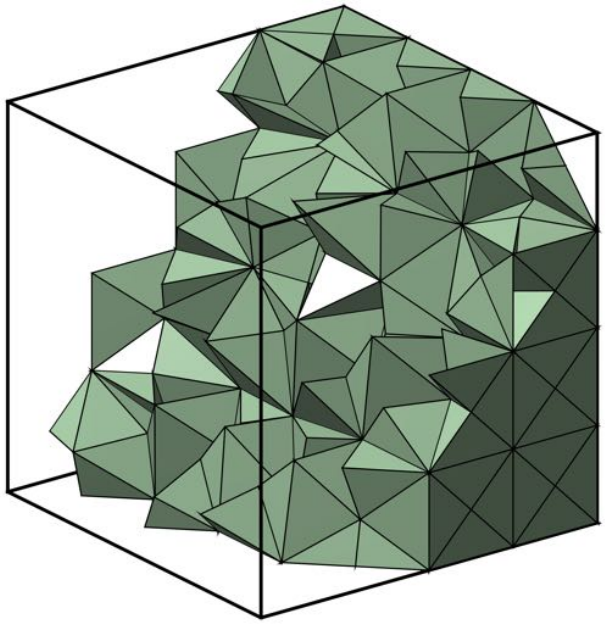}}		
	\subfigure[Prisms]{\includegraphics[width=0.24\textwidth]{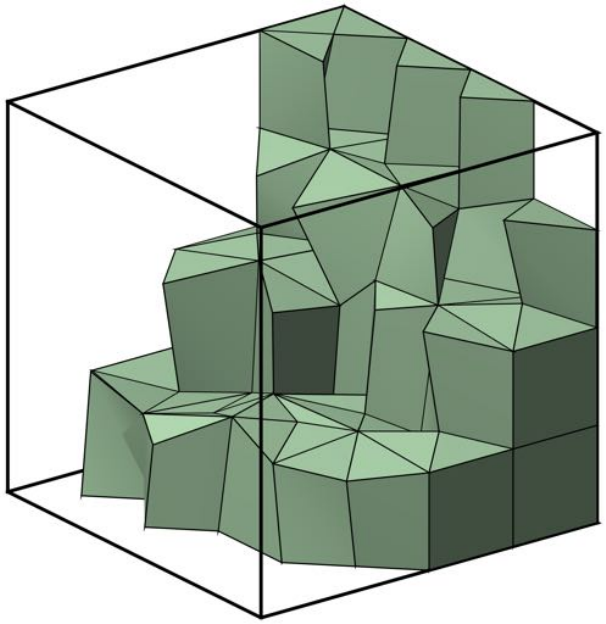}}
	\subfigure[Pyramids]{\includegraphics[width=0.24\textwidth]{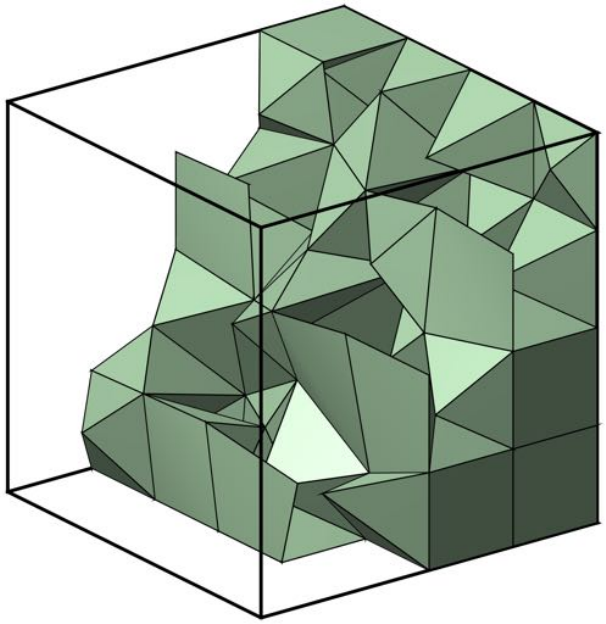}}
	\caption{Third level of mesh refinement for the irregular meshes of $\Omega=[0,1]^3$.}
	\label{fig:poisson3D_Irregularmeshes}
\end{figure}

To verify that the same behaviour is obtained for irregular three dimensional meshes, a new set of irregular three dimensional meshes is considered. A cut through the irregular meshes corresponding to the third level of refinement is represented in Figure~\ref{fig:poisson3D_Irregularmeshes} for all element types.

\begin{figure}[!tb]
	\centering
	\subfigure[$u$]{\includegraphics[width=0.4\textwidth]{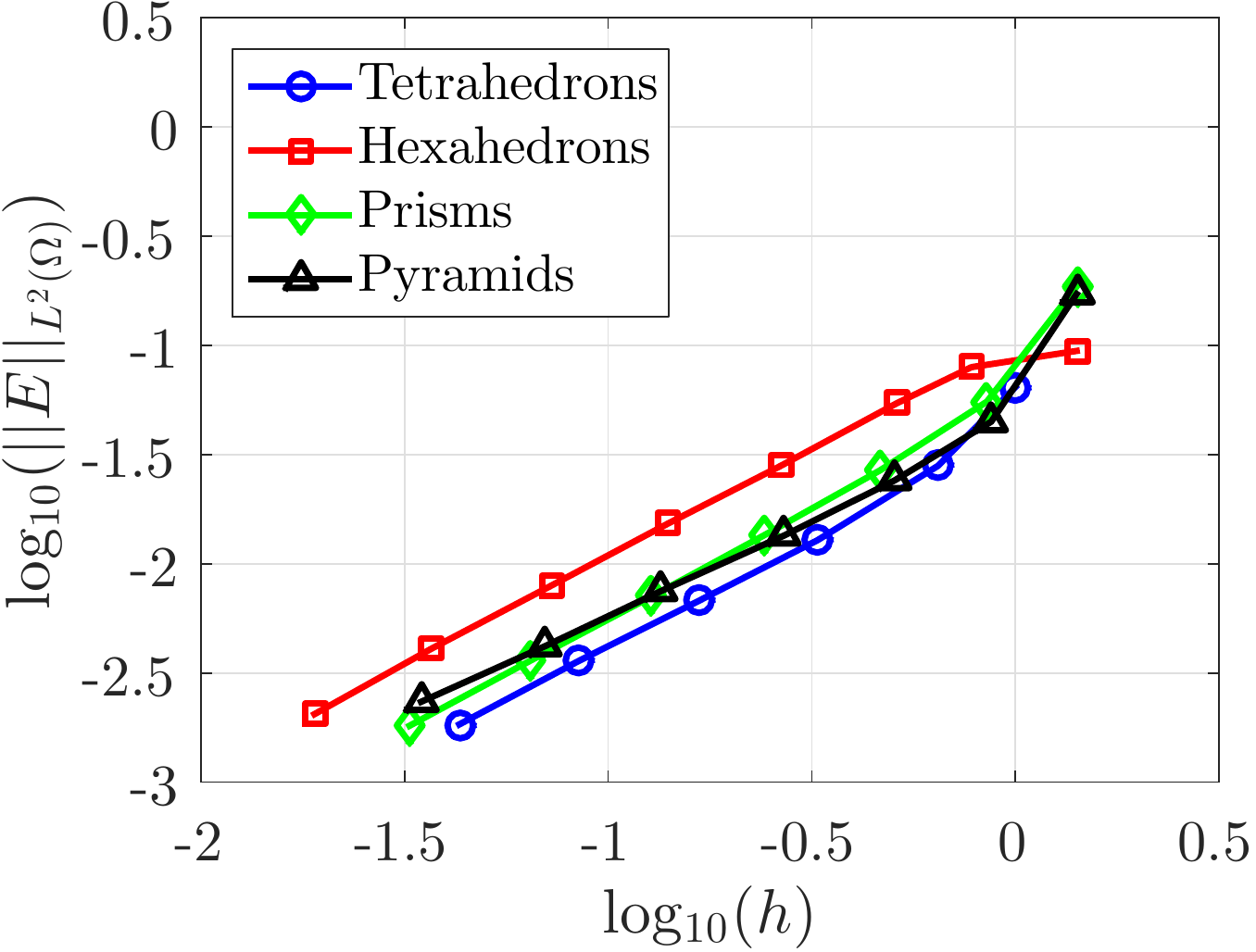}}
	\subfigure[$\bq$]{\includegraphics[width=0.4\textwidth]{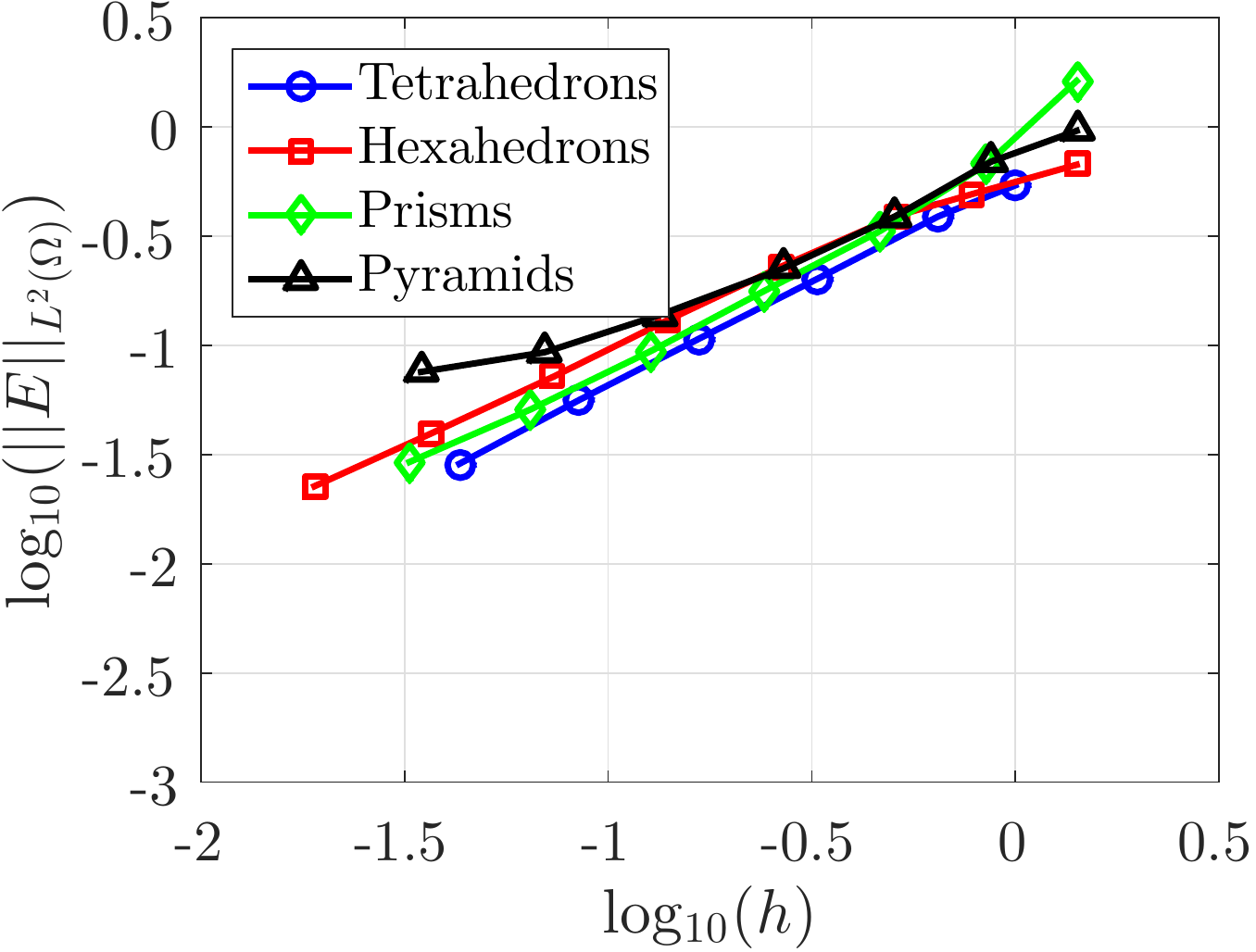}}
	\caption{Mesh convergence of the error of the solution and its gradient in the $\eltwo(\Omega)$ norm for the 3D Poisson problem with irregular meshes.}
	\label{fig:poisson3D_hConvIrregular}
\end{figure}

The results of the mesh convergence study for the three dimensional Poisson solver are presented in Figure~\ref{fig:poisson3D_hConvIrregular}.
All element types are able to provide the optimal rate of convergence for both the primal and dual variables, but a loss of accuracy is observed for the last pyramidal mesh. In all cases the loss of accuracy is observed for finer meshes and it is always observed for quadrilateral meshes in 2D or for meshes in 3D with elements containing quadrilateral faces. This is attributed to the large deformation introduced for the inner nodes, without any guarantee that quadrilateral elements (or elements with quadrilateral faces) remain convex.

\subsection{Influence of the mesh stretching}
\label{sc:InfluenceStretching}

\begin{figure}[!tb]
	\centering
	\subfigure[Mesh 3]{\includegraphics[width=0.24\textwidth]{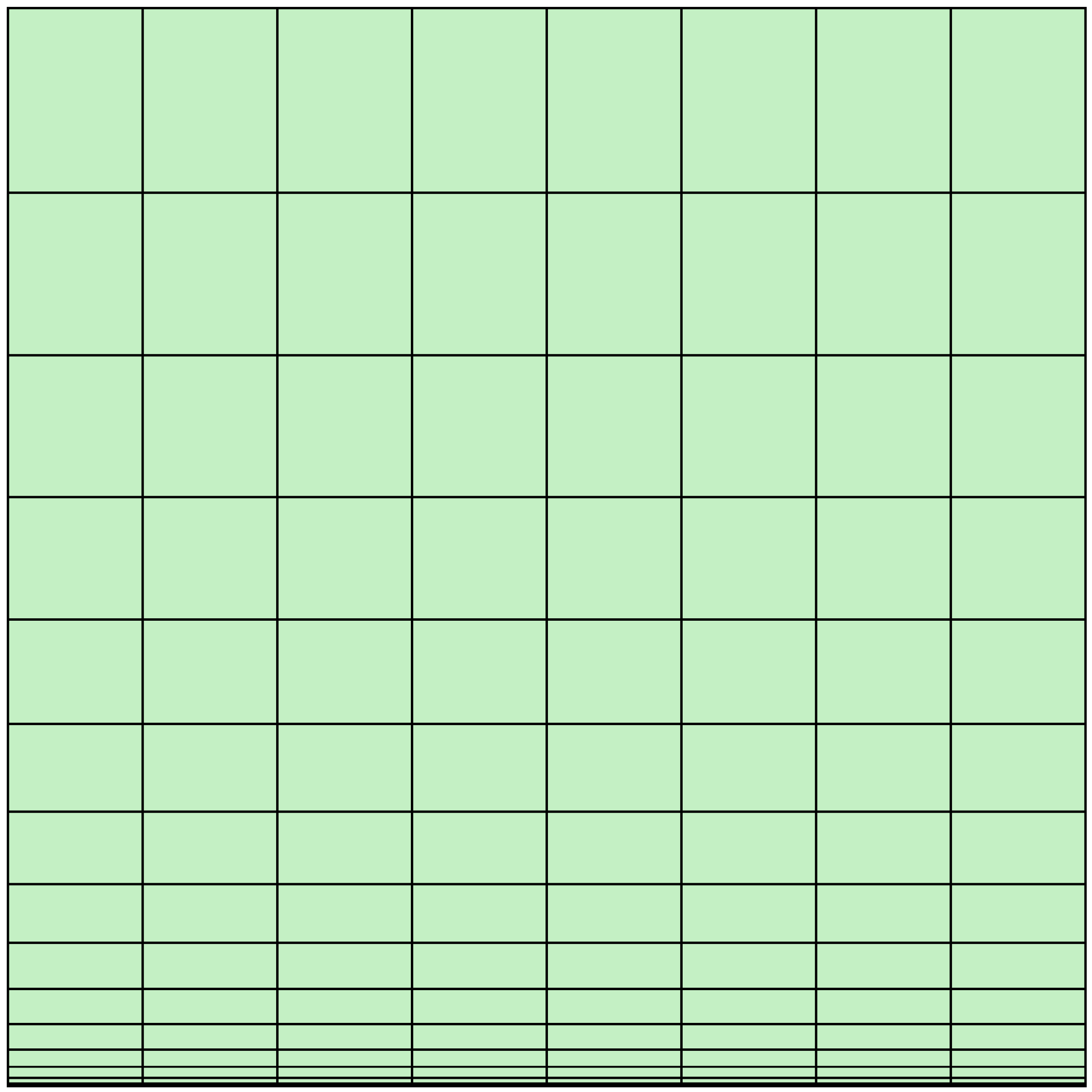}}
	\subfigure[Mesh 5]{\includegraphics[width=0.24\textwidth]{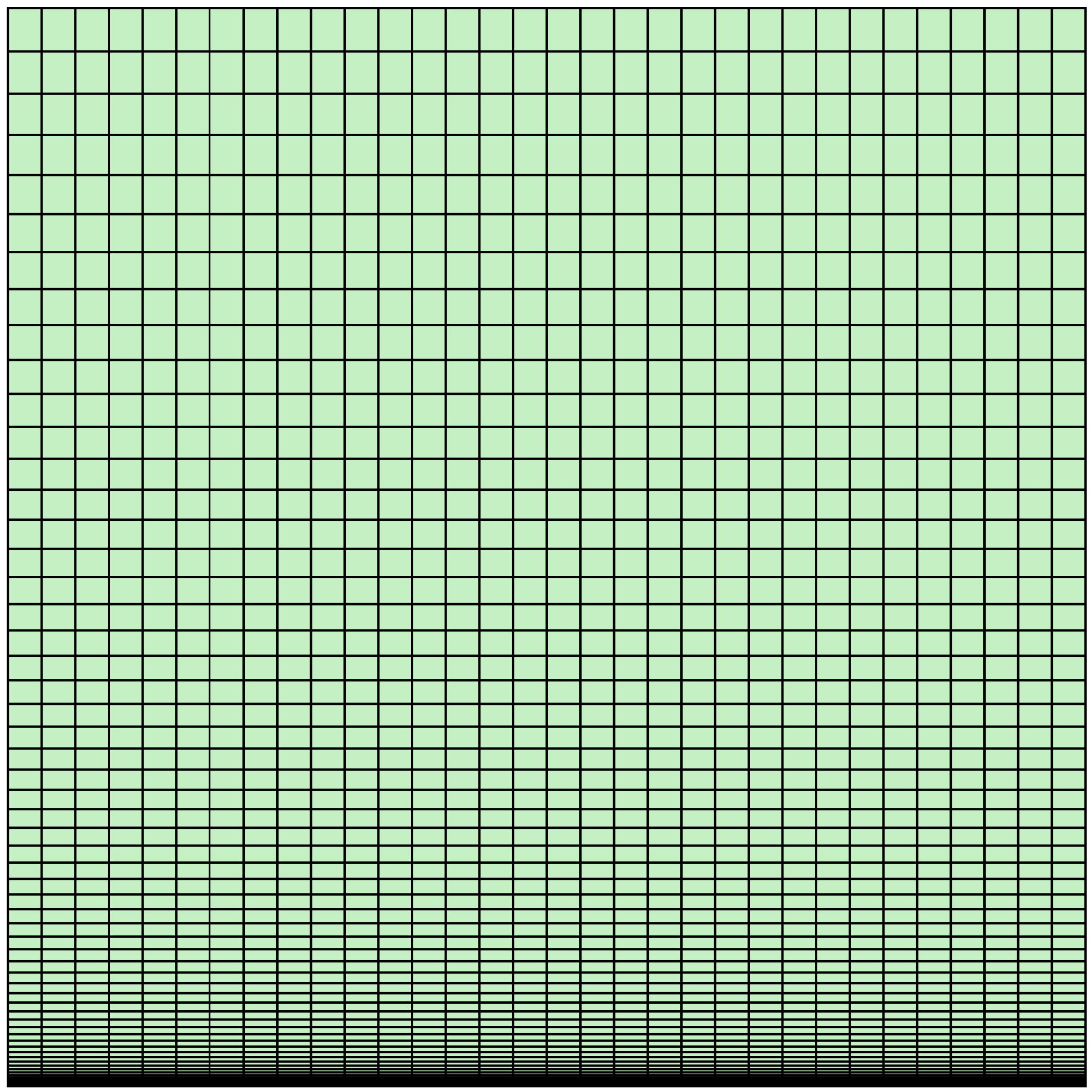}}
	\subfigure[Mesh 3]{\includegraphics[width=0.24\textwidth]{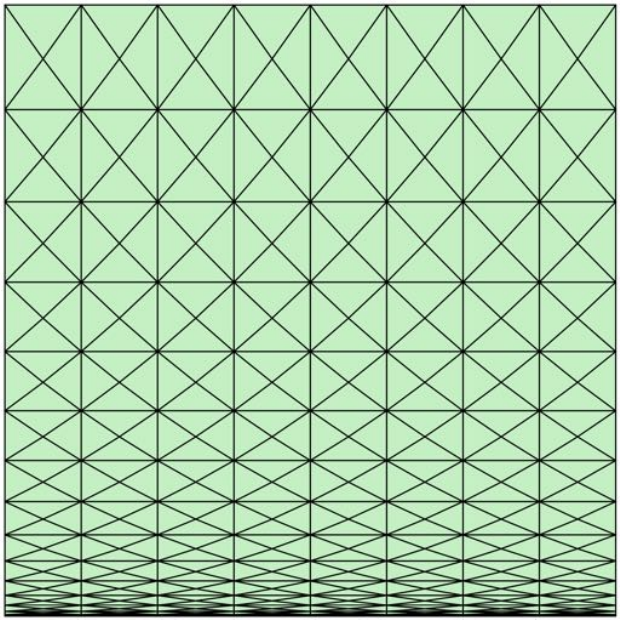}}
	\subfigure[Mesh 5]{\includegraphics[width=0.24\textwidth]{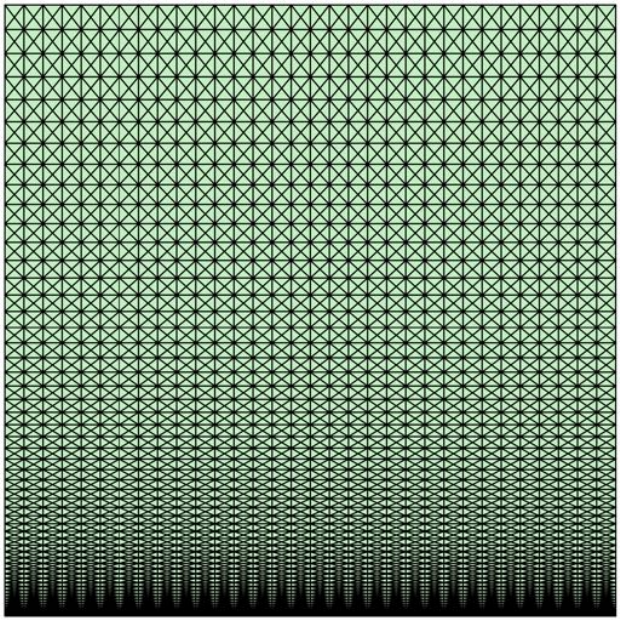}}
	\caption{Stretched quadrilateral and triangular meshes of $\Omega=[0,1]^2$.}
	\label{fig:poisson2D_meshesStretch}
\end{figure}

The last numerical study considers the influence of the mesh stretching on the accuracy and rate of convergence of the FCFV method. The regular meshes used in Sections~\ref{sc:PoissonVerification} and \ref{sc:StokesVerification} are modified to achieve a maximum given stretching $s$ near the bottom boundary. To construct the meshes in two dimensions, the vertical coordinate of the first layer is fixed to guarantee the desired stretching. The vertical coordinate of the subsequent layers is defined as
\begin{equation}
y_k = y_{k-1} + (h/s) \beta^{k-2}, \qquad \text{for $k=2,\ldots,N_y+1$}
\end{equation}
where $h$ is the maximum edge length of the regular mesh, $N_y$ is the number of elements in the vertical direction and the stretching factor $\beta$ is computed by imposing that the vertical coordinate of the last layer is one, that is finding the roots of 
\begin{equation}
(h/s) \beta^{N_y} - \beta + 1- (h/s) = 0.
\end{equation}
Two of the resulting stretched quadrilateral and triangular meshes are represented in Figure~\ref{fig:poisson2D_meshesStretch}, corresponding to a stretching  $s=100$.	

\begin{figure}[!tb]
	\centering
	\subfigure[$p$]{\includegraphics[width=0.32\textwidth]{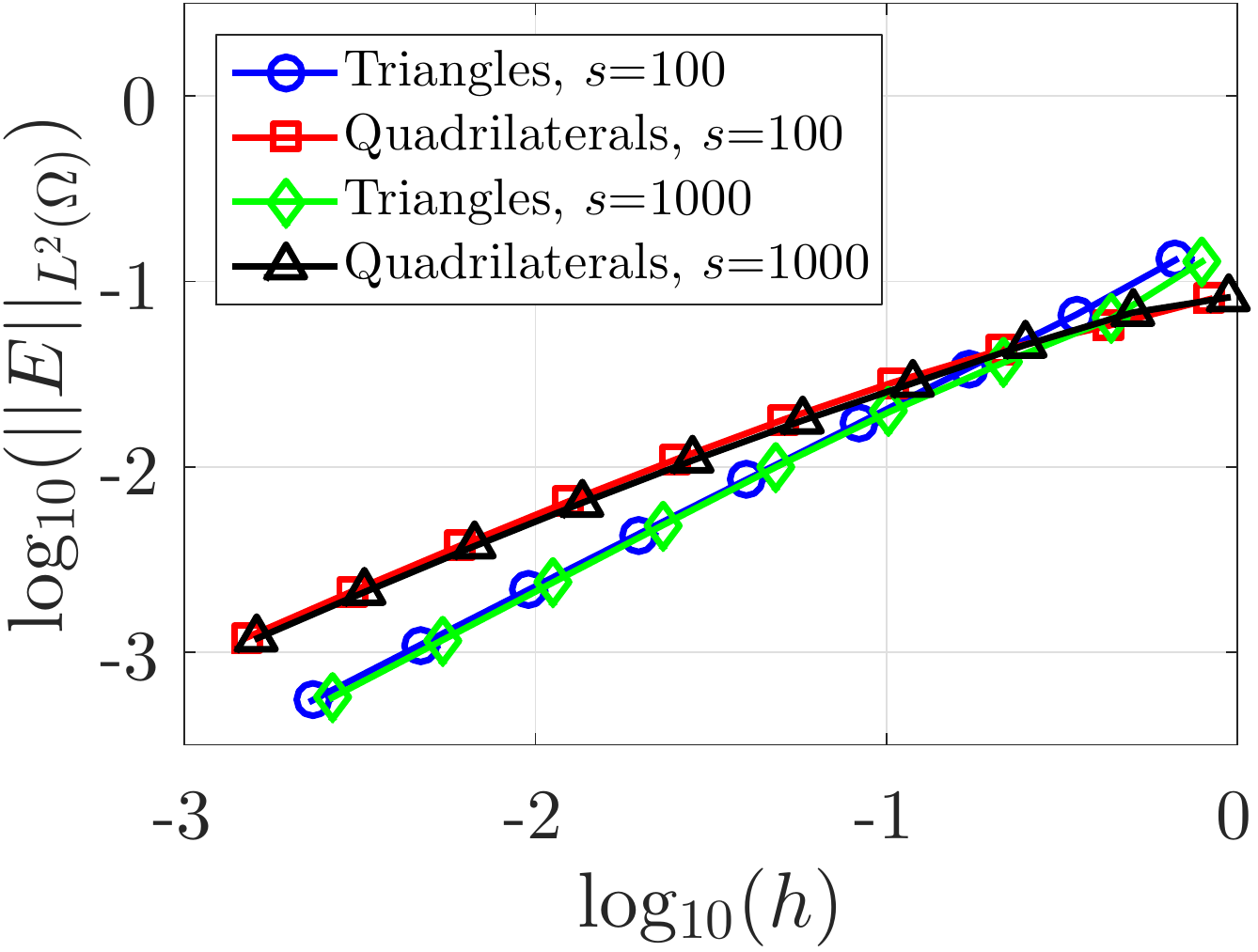}}
	\subfigure[$\bu$]{\includegraphics[width=0.32\textwidth]{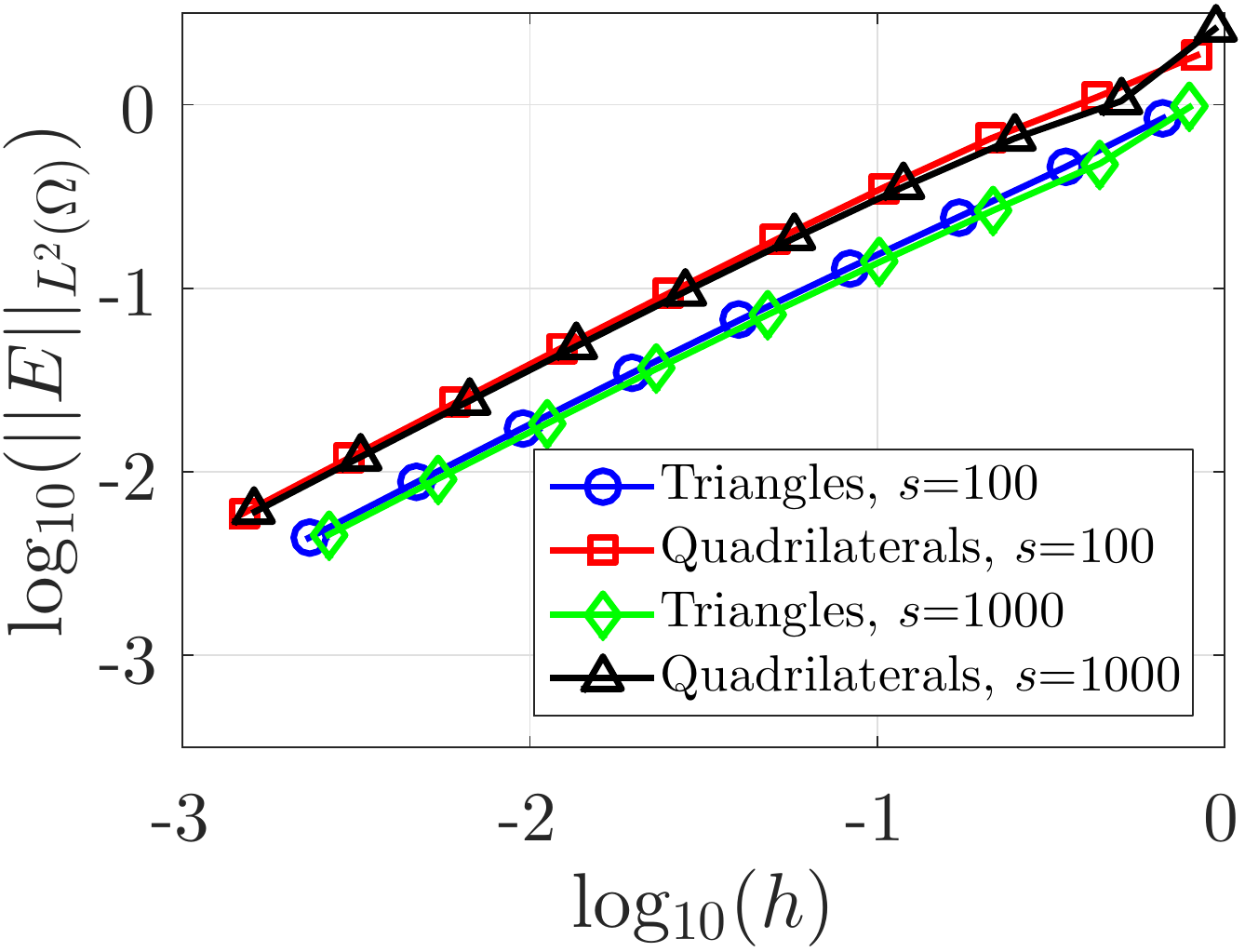}}
	\subfigure[$\bL$]{\includegraphics[width=0.32\textwidth]{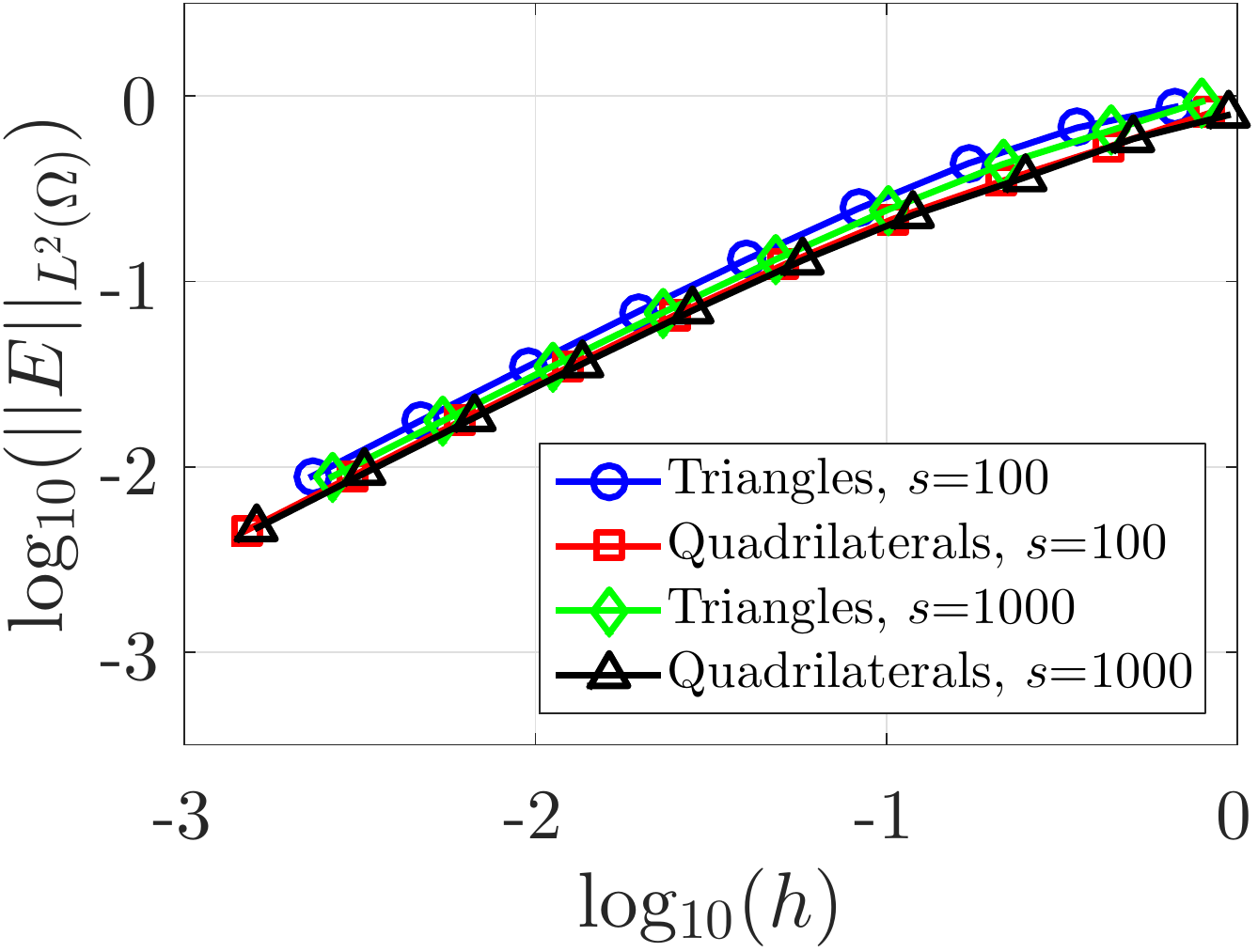}}
	\caption{Mesh convergence of the error of the solution and its gradient in the $\eltwo(\Omega)$ norm for the 2D Stokes problem with stretched meshes with stretching factor $s=100$ and $s=1,000$.}
	\label{fig:stokes2D_hConvStretch} 
\end{figure}

The convergence study on two dimensional stretched meshes is performed for the Stokes problem. Figure~\ref{fig:stokes2D_hConvStretch} depicts the convergence of the error of the pressure, velocity and gradient of the velocity for two different levels of stretching corresponding to $s=100$ and $s=1,000$ and for both triangular and quadrilateral elements. 
The results demonstrate the optimal convergence on highly stretched meshes. In addition, the results show that the accuracy of the FCFV method is not sensitive to the the level of stretching. In this example, the accuracy obtained with $s=100$ and $s=1,000$ is almost identical and comparable to the accuracy obtained in regular meshes in Section~\ref{sc:StokesVerification}.

\section{Numerical examples}
\label{sc:examples}

This section presents a series of large scale three dimensional examples to show the potential of the proposed FCFV methodology. The examples include the solution of the Poisson and the Stokes equations in 3D using tetrahedral meshes.

\subsection{Temperature distribution in a heat sink}
\label{sc:heatSink} 

\begin{figure}[!tb]
	\centering
	\subfigure[Geometry]{\includegraphics[width=0.32\textwidth]{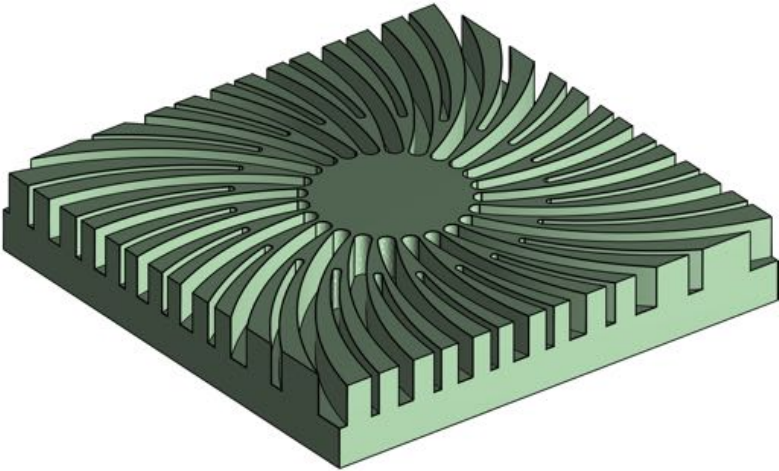}}
	\subfigure[Temperature]{\includegraphics[width=0.32\textwidth]{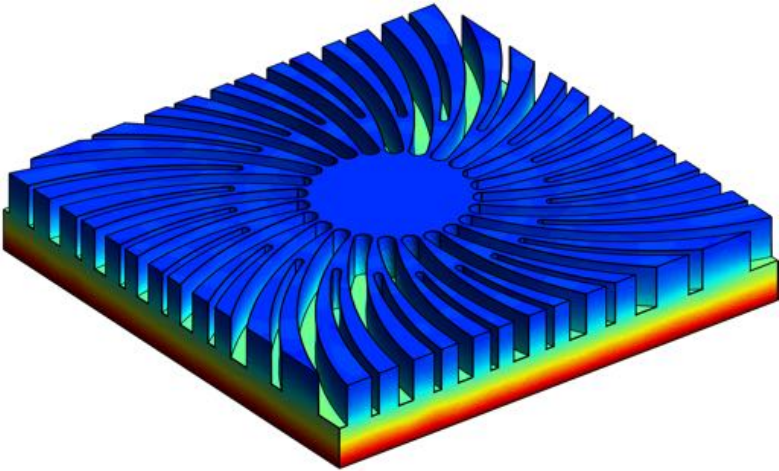}}
	\subfigure[Heat flux]{\includegraphics[width=0.32\textwidth]{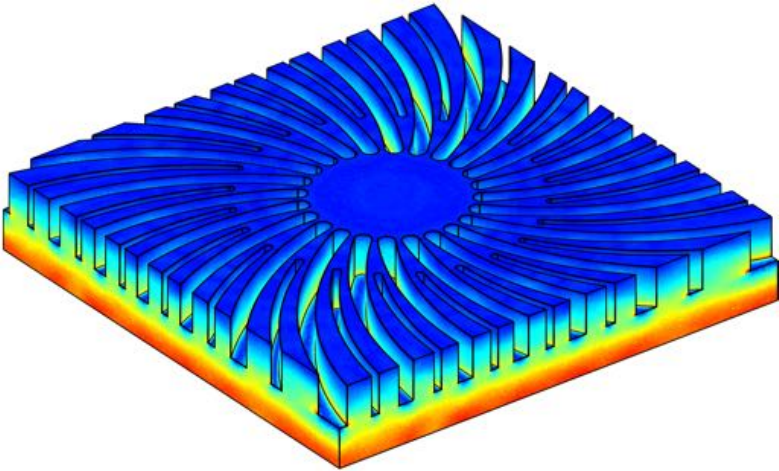}}
	\caption{a) Geometry of a heat sink, (b) temperature distribution and (c) magnitude of the heat flux vector on the surface of the heat sink.}
	\label{fig:heatSink}
\end{figure}

The first example considers the solution of the Poisson problem in a complex heat sink geometry. The analysis of such devices is of interest when designing and optimising heat sink modules that are fitted in many powerful electronic devices~\cite{chiang2005optimization,wang2012design}. The geometry of the heat sink is shown in Figure~\ref{fig:heatSink} (a). 
The geometry is defined by 198 NURBS surfaces and the generated mesh has 5,354,353 tetrahedral elements, 21,417,412 nodes, 10,490,943 internal faces and 435,526 external faces. A fixed temperature is imposed on the bottom part of the domain and a Neumann boundary condition corresponding to a negative heat flux in the rest of the boundary. The temperature distribution and the magnitude of the heat flux vector over the external faces of the mesh is represented in Figures~\ref{fig:heatSink} (b) and (c) respectively.

\subsection{Irrotational flow past complex aerodynamic configurations}
\label{sc:PoissonBoeing} 

Next, the irrotational flow around two complex three dimensional configurations is considered. The Poisson problem is solved with Neumann boundary conditions, corresponding to an imposed normal velocity. The potential function is imposed at one arbitrary point on the boundary to remove the indeterminacy of the potential. 

\begin{figure}[!tb]
	\centering
	\subfigure[Geometry]{\includegraphics[width=0.32\textwidth]{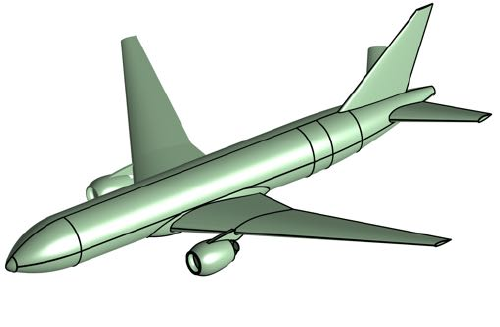}}
	\subfigure[Velocity]{\includegraphics[width=0.32\textwidth]{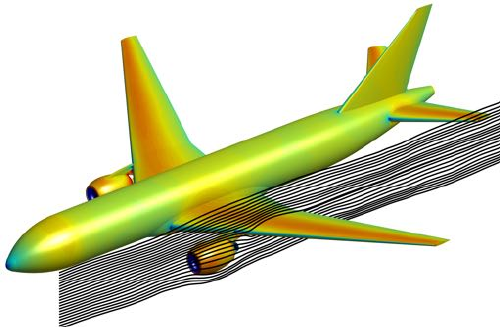}}
	\subfigure[Pressure]{\includegraphics[width=0.32\textwidth]{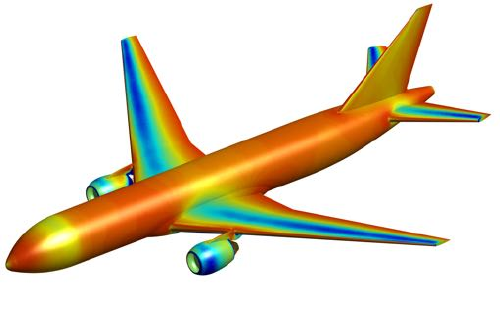}}
	\caption{Magnitude of the velocity and pressure distribution for the irrotational flow around a full aircraft configuration.}
	\label{fig:poissonBoeing}
\end{figure}

First, the flow around a full aircraft configuration is considered. The computational domain is meshed with 5,125,998 tetrahedral elements. The mesh contains 9,220,701 internal faces and 2,062,412 faces on the Neumann boundary, so the total number of unknowns of the global problem is 11,283,113. The magnitude of the velocity computed using the dual variable (i.e. $v = \|\bq\|_2$) and the pressure distribution computed from Bernoulli equation are represented in Figure~\ref{fig:poissonBoeing} on the surface of the aircraft. 
The computation took 3.7 minutes for the assembly of the global system of equations (computation of the elemental matrices plus assembly of the global matrix) and 5.7 minutes to solve the global problem using a direct solver. The developed code is written in Matlab and the computation was performed in an Intel$^{\tiny{\textregistered}}$
Xeon$^{\tiny{\textregistered}}$ CPU $@$ 3.70GHz and 32GB main memory available.

\begin{figure}[!tb]
	\centering
	\subfigure[Geometry]{\includegraphics[width=0.32\textwidth]{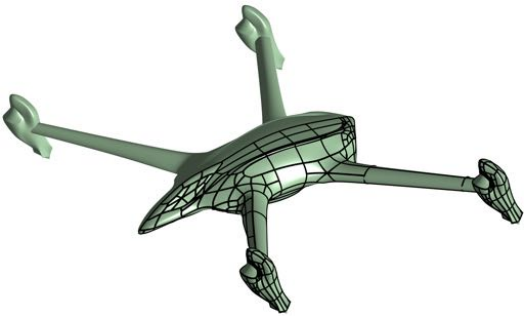}}
	\subfigure[Velocity]{\includegraphics[width=0.32\textwidth]{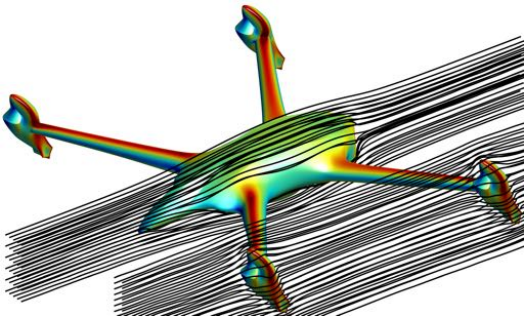}}
	\subfigure[Pressure]{\includegraphics[width=0.32\textwidth]{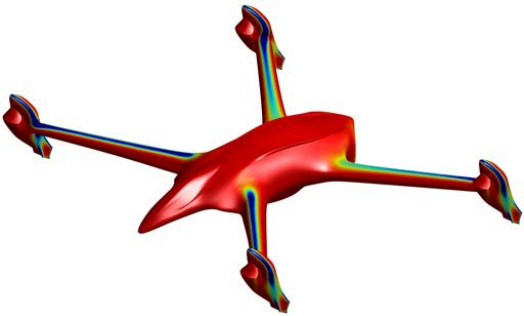}}	
	\caption{(a) Geometry of a generic drone, (b) magnitude of the velocity and streamlines and (c) pressure field over the surface of the drone.}
	\label{fig:poissonDrone}
\end{figure}

The second configuration considers a more challenging geometry corresponding to a generic drone represented in Figure~\ref{fig:poissonDrone} (a), where the 376 NURBS surfaces that define half of the geometry are highlighted.
The computational domain is meshed with 4,093,200 tetrahedral elements. The mesh contains 7,312,154 internal faces and 1,743,273 faces on the Neumann boundary, so the total number of unknowns of the global problem is 9,055,427. The magnitude of the velocity, computed using the dual variable and the pressure distribution are also represented in Figure~\ref{fig:poissonDrone} on the surface of the drone. The computation took 2.8 minutes for the assembly of the global system of equations (computation of the elemental matrices plus assembly of the global matrix) and 5.4 minutes to solve the global problem using a direct solver.

\subsection{Stokes flow past a sphere}
\label{sc:stokesSphere} 

The next example considers a classical test case for three dimensional Stokes solvers, the flow around a sphere. Using the known analytical solution~\cite{batchelor2000introduction}, this problem is used to show the optimal convergence of the FCFV for a three dimensional Stokes flow and to illustrate the accuracy of the proposed technique when evaluating quantities of interest such as the drag force.

The domain is defined as $\Omega = \left( [-7,15] \times [-5,5] \times [-5,5] \right)\setminus \mathcal{B}_{1,\mathbf{0}}$, where $\mathcal{B}_{1,\mathbf{0}}$ denotes a ball of unit radius centred at the origin. Seven tetrahedral meshes of the domain are considered with the number of elements, nodes, faces, characteristic element size ($h$) and induced number of degrees of freedom (\ndof) detailed in Table~\ref{tab:meshesSphere}.
\begin{table}[hbt]
\centering
\begin{tabular}[hbt]{| c | c | c | c | c |}
	\hline
	Elements & Nodes & Faces & $h$ & $\ndof$ \\
	\hline & & & &  
	\\ [-1em] 
	\hline
	3,107 & 12,428 & 6,560 & 2.2552 & 20,711 \\
	\hline
	10,680 & 42,720 & 22,197 & 1.5424 & 72,249 \\
	\hline
	43,682 & 174,728 & 89,530 & 1.0095 & 299,276 \\
	\hline
	204,099 & 816,396 & 414,457 & 0.6528 & 1,409,916 \\
	\hline
	686,853 & 2,747,412 & 1,387,771 & 0.4523 & 4,765,776 \\
	\hline
	2,516,099 & 10,064,396 & 5,065,404 & 0.3097 & 17,513,075 \\
	\hline
	7,604,928 & 30,419,712 & 15,279,422 & 0.2172 & 53,025,798 \\
	\hline
\end{tabular}
\caption{Details of the seven tetrahedral meshes to study the Stokes flow past a sphere.}
\label{tab:meshesSphere}
\end{table}

\begin{figure}[!tb]
	\centering
	\subfigure[Velocity]{\includegraphics[width=0.4\textwidth]{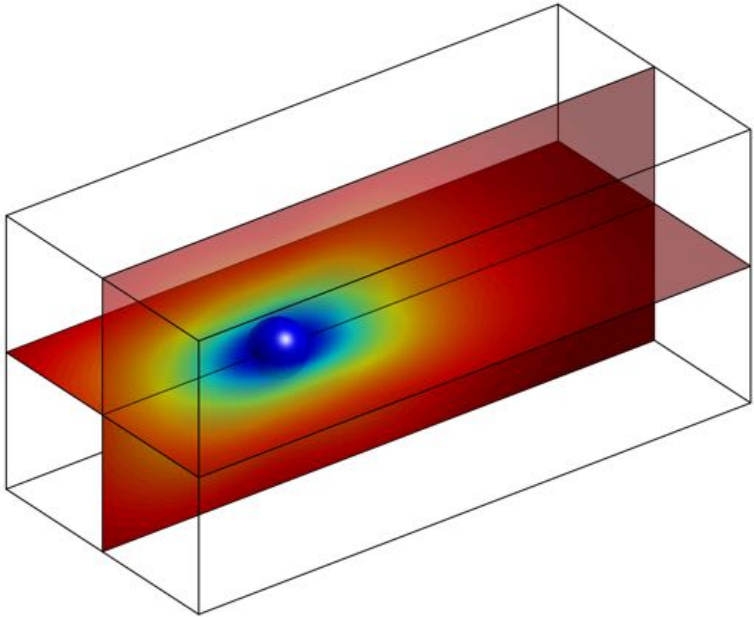}}
	\hspace{2cm}
	\subfigure[Pressure]{\includegraphics[width=0.4\textwidth]{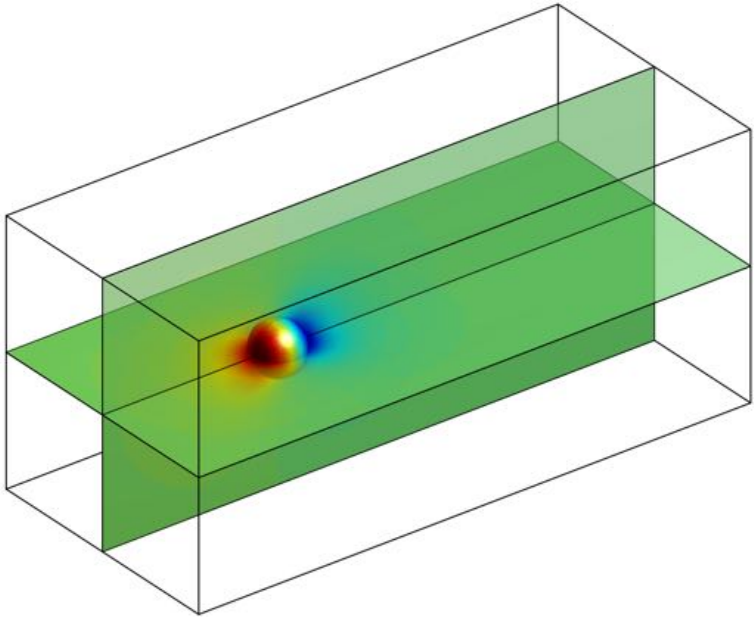}}
	\caption{Magnitude of the velocity and pressure distribution for the Stokes flow past a sphere.}
	\label{fig:stokesSphere}
\end{figure}

\begin{figure}[!tb]
	\centering
	\subfigure[]{\includegraphics[width=0.4\textwidth]{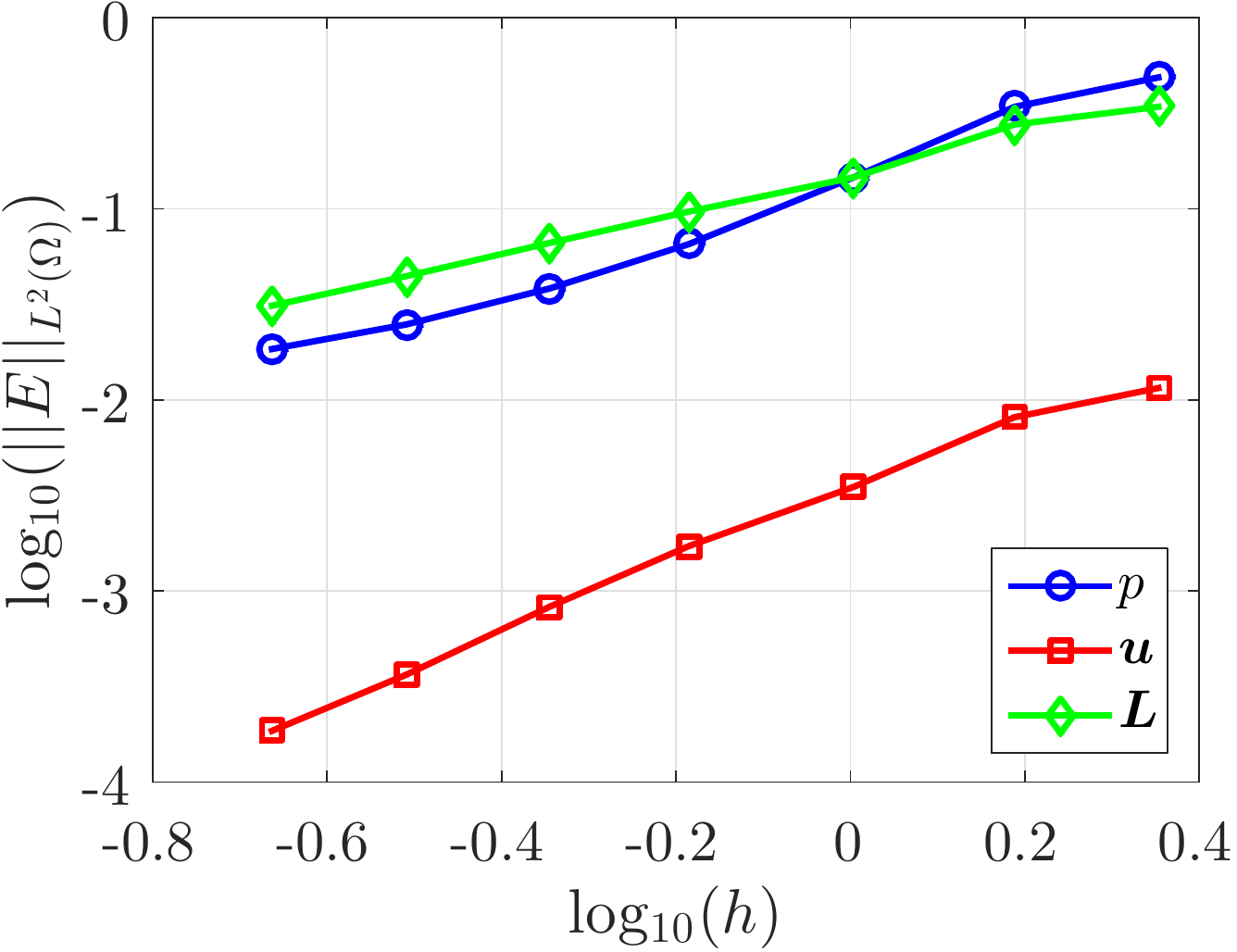}}
	\subfigure[]{\includegraphics[width=0.4\textwidth]{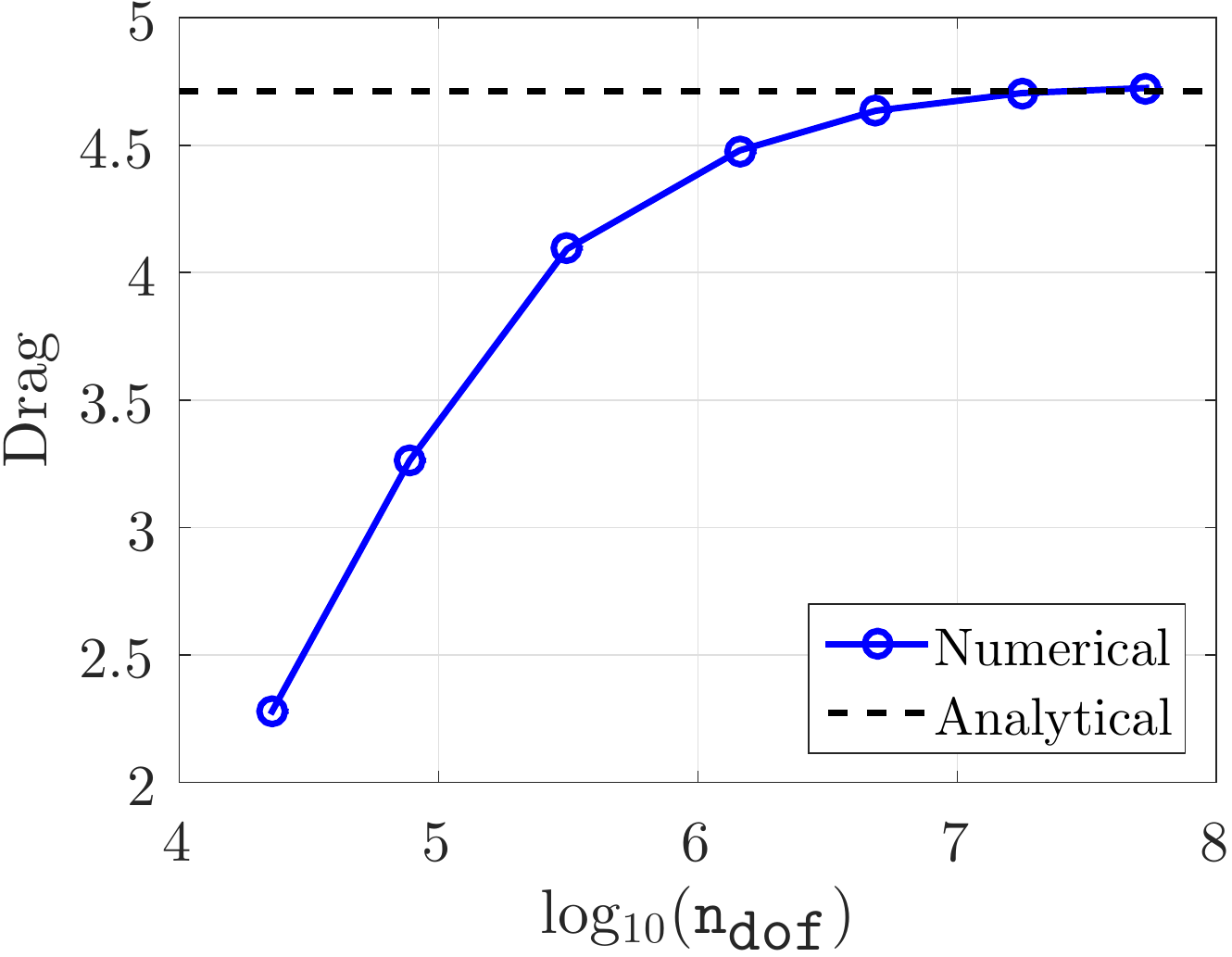}}
	\caption{(a) Mesh convergence of the error of the pressure, the velocity and the velocity gradient in the $\eltwo(\Omega)$ norm for the Stokes flow around a sphere and (b) convergence of the drag as a function of the number of degrees of freedom.}
	\label{fig:stokesSphereHConvDrag}
\end{figure}

The magnitude of the velocity and the pressure field are represented in Figure~\ref{fig:stokesSphere} over the surface of the sphere and the symmetry planes.
The evolution of the error in pressure, velocity and the gradient of the velocity, measured in the $\eltwo(\Omega)$ norm is represented in Figure~\ref{fig:stokesSphereHConvDrag} (a). 
Similar to previous examples, the FCFV shows significantly more accurate results on the velocity than in the pressure and the gradient of the velocity. The rate of convergence observed for the velocity corresponds to the higher than optimal (quadratic in this example), whereas in the pressure and the gradient of the velocity the expected linear rate is obtained.

The results in Figure~\ref{fig:stokesSphereHConvDrag} (b) show the convergence of the drag force as the number of degrees of freedom is increased. Using the sixth mesh, with approximately 2.5 million elements, an error on the drag force below 0.2\% is obtained, showing the potential of the proposed FCFV approach. The simulation took 10 minutes for the computation of all elemental matrices and 3 minutes for the assembly of the global system. The solution of the global system was performed using the biconjugate gradient method in a single processor and without pre-conditioner, taking less than 5 hours. The developed code is written in Matlab and the computation was performed in an Intel$^{\tiny{\textregistered}}$
Xeon$^{\tiny{\textregistered}}$ CPU $@$ 3.20GHz and 70GB main memory available.

\subsection{Stokes flow past a porous sphere}
\label{sc:stokesClusterSpheres} 

The last example, inspired by the results presented in~\cite{wittig2017drag}, considers the Stokes flow past a porous sphere that is formed by a cluster of solid spherical particles. This problem is of great interest in a variety of chemical engineering applications~\cite{neale1973creeping,wittig2017drag} (e.g. flow through catalysts) and natural processes (e.g. sedimentation).

The domain is defined as $\Omega = \left( [-5,10] \times [-5,5] \times [-5,5] \right)\setminus \bigcup_{i=1}^{N_s} \mathcal{B}_{\rho,\bm{x}_i}$, where $\mathcal{B}_{\rho,\bm{x}_i}$ denotes a ball of radius $\rho$ centred at $\bm{x}_i$ and $N_s$ is the total number of balls. 
In the present simulation, the arrangement of spherical particles is constructed so that no overlap between the spheres is possible. Given the number of spherical layers $n$, the radius of the spherical particles $\rho$, the radius of the spherical region $R$ containing all the particles and the minimum allowed gap between the particles $\delta$, the procedure devised to compute the centre of the spheres and the number of spheres to be generated within each layer is described next.

First, the radius of each spherical region that will be used to place the spheres of each layer is computed as $R_i = i(R-\rho)/n$, for $i=1,\ldots,n$. On each layer, a sphere is initially placed on the south pole of the spherical region with radius $R_i$ and the number of parallel arcs is computed as $a_i = \lfloor \pi R_i / d \rfloor$, where $d=2\rho + \delta$. The radius of the parallel arc is computed as $r_{i,j} = R_i\sin(jd/R_i)$ for $j=1,\ldots,a_i$. Finally, the centres of the spheres are equally-spaced along each parallel arc.

\begin{figure}[!tb]
	\centering
	\subfigure[]{\includegraphics[width=0.3\textwidth]{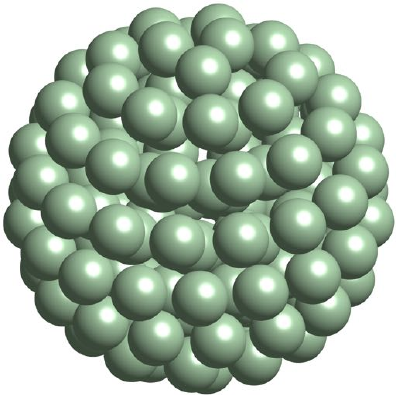}}
	\hspace{2cm}
	\subfigure[]{\includegraphics[width=0.3\textwidth]{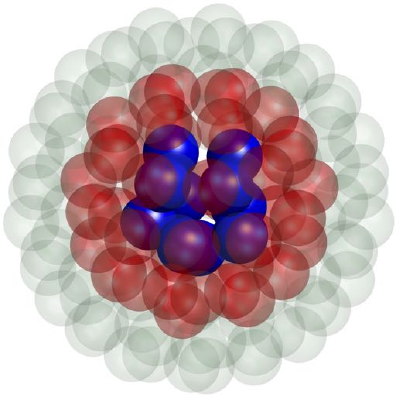}}
	\caption{Arrangement of spherical particles corresponding to $n=3$, $R=1$, $\rho=1/7$ and $\delta = 0.2\rho$. The colours in (b) are used to help the visualisation of the particles in each one of the three layers.}
	\label{fig:stokesPorousSpheresGeo}
\end{figure}

Figure~\ref{fig:stokesPorousSpheresGeo} (a) shows the arrangement of $N_s=126$ spherical particles corresponding to $n=3$, $R=1$, $\rho=1/7$ and $\delta = 0.2\rho$. Figure~\ref{fig:stokesPorousSpheresGeo} (b) shows the same distribution of particles using different colours to distinguish the particles in each one of the three spherical layers.
The arrangement of particles considered leads to a porosity of 0.6327, computed as $1 - N_s (\rho/R)^3$,
which is within the range of experiments considered in~\cite{wittig2017drag}.

A tetrahedral mesh with 9,646,810 elements and 38,587,240 nodes is utilised to compute the Stokes flow past the porous sphere formed by 126 spherical particles.  The mesh contains a total of 17,816,283 internal faces and 2,954,674 external faces, leading to a global problem with 63,095,659 degrees of freedom when Dirichlet boundary conditions are considered in the whole domain. 

\begin{figure}[!tb]
	\centering
	\subfigure[Velocity]{\includegraphics[width=0.45\textwidth]{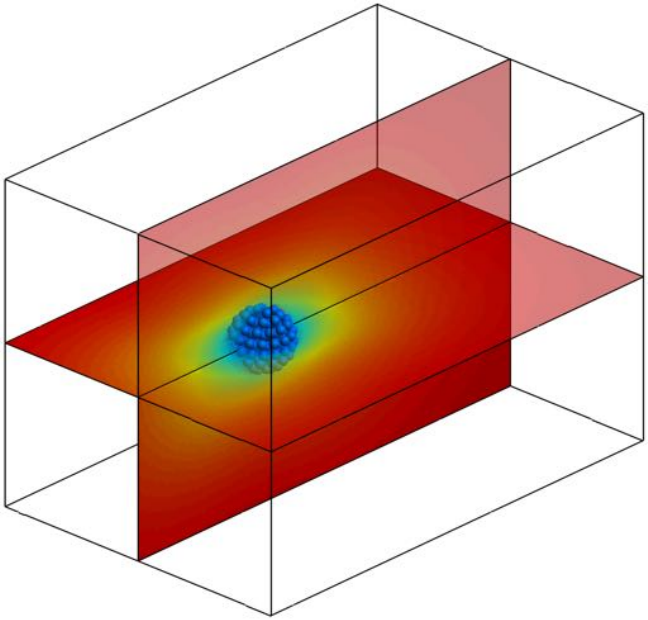}}
	\subfigure[Pressure]{\includegraphics[width=0.45\textwidth]{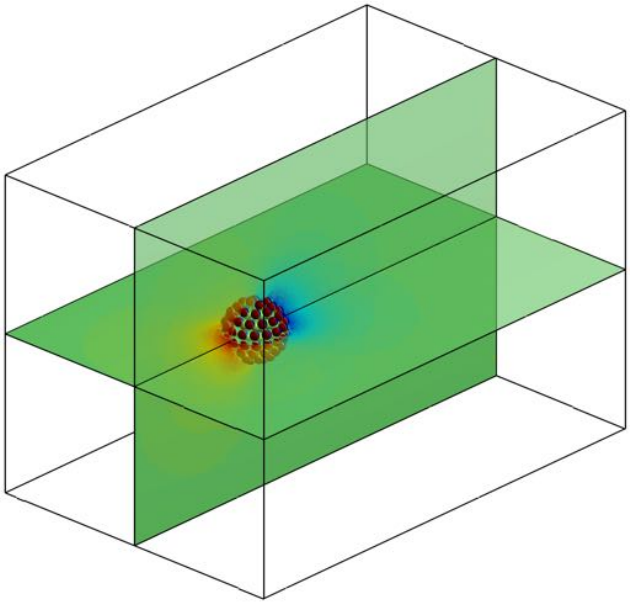}}
	\caption{Magnitude of the velocity and pressure distribution for the Stokes flow past a porous sphere.}
	\label{fig:stokesPorousSpheresSol}
\end{figure}

\begin{figure}[!tb]
	\centering
	\subfigure[Pressure]{\includegraphics[width=0.25\textwidth]{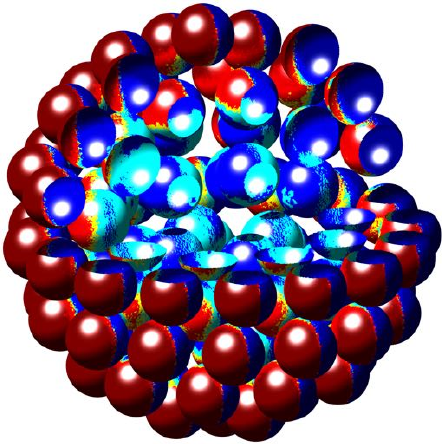}}
	\subfigure[Velocity and streamlines]{\includegraphics[width=0.74\textwidth]{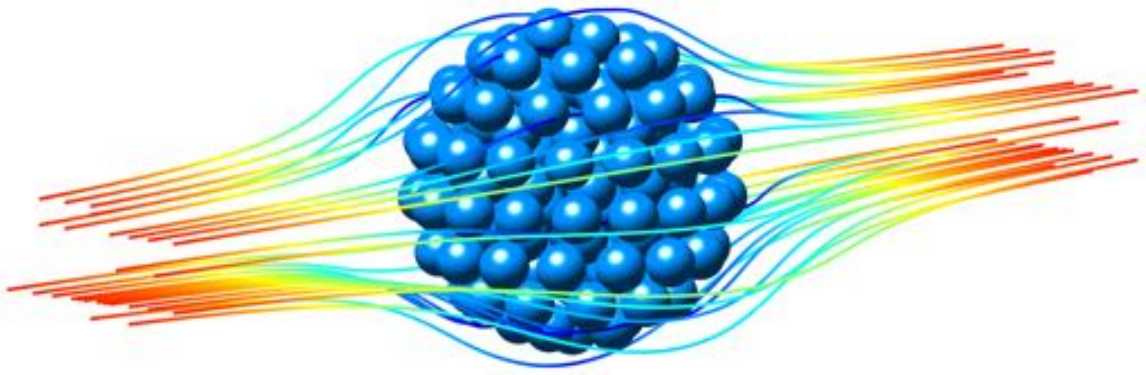}}
	\caption{(a) Pressure distribution over some spherical particles and (b) velocity and streamlines.}
	\label{fig:stokesPorousSpheresPress}
\end{figure}

The magnitude of the velocity and the pressure field are represented in Figure~\ref{fig:stokesPorousSpheresSol} over the surface of the spherical particles and two sections of the computational domain.
The results clearly show macroscopic behaviour of the flow around the spherical particles is very similar to the flow pattern obtained for the flow around a single sphere studied in Section~\ref{sc:stokesSphere}.
The pressure over the spherical particles in Figure~\ref{fig:stokesPorousSpheresSol} (b) has been amplified by a factor of $10^4$ to enable distinguish the high pressure over the first layer of spheres from the pressure over the inner layers. This phenomenon can be better observed in Figure~\ref{fig:stokesPorousSpheresPress} (a), where the pressure field over the spherical particles is depicted. A cut through the domain has been performed to enable the visualisation of some particles from the inner layers and therefore appreciate the high difference of pressure over the spheres in the outer  layer compared to the pressure over the spheres in the inner layers. Finally, Figure~\ref{fig:stokesPorousSpheresPress} (b) shows the velocity in the spherical particles together with some streamlines coloured according to the magnitude of the velocity.

\section{Concluding remarks}
\label{sc:Conclusion}

This papers proposes a new finite volume paradigm, called face-centred finite volume (FCFV), based on a hybridisable discontinuous Galerkin (HDG) method with constant degree of approximation. As any other HDG method, the FCFV technique requires the solution of a global system of equations whose size is equal to the total number of element faces. The solution and its gradient in each element are then recovered by solving a set of independent element-by-element problems. First order convergence on both the solution and its gradient is obtained without a reconstruction of the gradients. Therefore, contrary to other finite volume methodologies, the accuracy of the FCFV method is not compromised in the presence of highly stretched or distorted elements.

The application of the proposed method to scalar and vector second-order elliptic problems is considered, namely the Poisson and the Stokes equations. For the solution of Stokes flow problems, the FCFV method does not require the solution of a Poisson problem for computing the pressure, as required by segregated schemes, and does not require the use of different approximation spaces to satisfy the LBB condition, as required by other mixed finite element methods.

An exhaustive set of numerical studies has been presented to verify the optimal approximation properties of the method, to study the influence of the HDG stabilisation parameter, to analyse and compare the computational cost for different element types and to check the robustness of the method when using distorted and stretched meshes. The studies compromise two and three dimensional cases for both the Poisson and Stokes problems. The results show that, contrary to other FV methods, the accuracy of the FCFV method is not sensitive to mesh distortion and stretching. In addition, the FCFV method shows its better performance, accuracy and robustness using simplicial elements, facilitating its application to problems involving complex geometries in three dimensions. To illustrate the potential of the FCFV method, a set of more challenging three dimensional examples are used to illustrate the potential and efficiency of the methodology. 


\bibliographystyle{abbrv}
\bibliography{Ref-HDG}

\appendix
\section{FCFV method with Neumann local problems}\label{sc:Nlocal}

As described in~\cite{RS-SH:16}, a minor modification of a classical HDG formulation can be devised to obtain a smaller global problem. This modification consists of prescribing the Neumann boundary conditions in the local problem, rather than in the global problem, as done in Sections~\ref{sc:Poisson} and \ref{sc:Stokes}. This section presents the changes induced by this modification for the FCFV method applied to the Poisson problem. The derivation for the Stokes problem follows the same rationale.

\subsection{Strong form of the local and global problems}

First, the local (element-by-element) problem with both Dirichlet and Neumann boundary conditions is defined, namely
\begin{equation} \label{eq:Nlocal-strong}
	\left\{\begin{aligned}
		\bq_e + \grad u_e &=\bm{0}  &&\text{in $\Omega_e$, and for $e=1,\dotsc \numel$}\\
		\grad\cdot\bq_e &= s          &&\text{in $\Omega_e$, and for $e=1,\dotsc \numel$}\\		
		u_e &= u_D     &&\text{on $\partial\Omega_e\cap\Gamma_D$,}\\
		\bn_e\cdot\bq_e &= -t         &&\text{on $\partial\Omega_e\cap\Gamma_N$,}\\
		u_e &=\hu  &&\text{on $\partial\Omega_e\setminus\partial\Omega$,}
	\end{aligned} \right.
\end{equation}
for $e=1,\dotsc \numel$. In each element, $\Omega_e$, this problem produces an element-by-element solution $\bq_e$ and $u_e$ as a function of the unknown $\hu\in\eltwo(\Gamma)$. 

Second, a global problem is defined to determine $\hu$. It corresponds to the imposition of the  \emph{transmission condition}, 
\begin{equation} \label{eq:Ntransmission}
	\jump{\bn\cdot \bq} = 0  \text{ on $\Gamma$.}
\end{equation}

\subsection{Weak form of the local and global problems}

The discrete weak formulation of the previously introduced local problems is obtained by multiplying the problems by a test function in an appropriate discrete functional space and integrating by parts. For $e=1,\dotsc \numel$, seek $(\bq_e^h ,u_e^h)\in [\Vh(\Omega_e)]^{\nsd}\times\Vh(\Omega_e)$ such that, for all $(\bw ,v) \in [\Vh(\Omega_e)]^{\nsd}\times\Vh(\Omega_e)$
\begin{subequations}\label{eq:HDG-Nlocal}
  \begin{equation} \label{eq:Nweak-Nlocal1}
	(\grad\cdot\bw, u_e^h)_{\Omega_e} 
	- \langle \bn \cdot\bw, u_e^h \rangle_{\partial\Omega_e \cap \Gamma_N} 
	- (\bw,\bq_e^h)_{\Omega_e} \\
	= \langle \bn\cdot\bw, u_D \rangle_{\partial\Omega_e \cap \Gamma_D}
	+ \langle \bn\cdot\bw, \hu^h \rangle_{\partial\Omega_e\setminus\partial\Omega},
  \end{equation}	
  \begin{multline} \label{eq:Nweak-Nlocal2}
	\langle v,\tau_e\, u_e^h \rangle_{\partial\Omega_e\setminus\Gamma_N}
	+ (v, \grad\cdot\bq_e^h)_{\Omega_e} 
	- \langle v ,\bn_e\cdot\bq_e^h \rangle_{\partial\Omega_e\cap\Gamma_N} \\
	= (v,s)_{\Omega_e}
	+ \langle v,t \rangle_{\partial\Omega_e \cap \Gamma_N} 
	+ \langle v,\tau_e\, u_D \rangle_{\partial\Omega_e\cap\Gamma_D}  
	+ \langle v,\tau_e\,\hu^h \rangle_{\partial\Omega_e\setminus\partial\Omega}.	
  \end{multline}
\end{subequations}
It is worth noting that $\hu^h \in \VhHat(\Gamma)$ is not defined along $\Gamma_N$ and, consequently, $u_e^h$ is left along $\partial\Omega_e\cap\Gamma_N$. In addition, a new definition for the numerical traces of the normal fluxes has been introduced. They are defined element-by-element (i.e. for $e=1,\dotsc \numel$) as 
\begin{equation} \label{eq:N-EBENumFlux}
	\bn_e\cdot\widehat{\bq}_e^h := \begin{cases}
		\bn_e\cdot\bq_e^h +\tau_e (u_e^h- u_D    ) & \text{on $\partial\Omega_e\cap\Gamma_D$,} \\
		\bn_e\cdot\bq_e^h +\tau_e (u_e^h- \hu^h) & \text{on $\partial\Omega_e\cap\Gamma$,}  \\
		-t                                                       & \text{on $\partial\Omega_e\cap\Gamma_N$.} 
	\end{cases}
\end{equation}

For the global problem, continuity of the fluxes is now only imposed along the internal faces, see \eqref{eq:Ntransmission}. Hence, the global weak problem is: find $\hu^h \in \VhHat(\Gamma)$ for all $\hv \in \VhHat(\Gamma)$ such that
\begin{equation}\label{eq:Nglobal}
	\sum_{e=1}^{\numel}
	\langle \hv, \bigl[\bn_e\cdot\bq_e^h + \tau_e (u_e^h-\hu^h)\bigr] \rangle_{\partial\Omega_e\setminus\partial\Omega} =0 ,
\end{equation}
where the definition of the numerical normal flux in Equation~\eqref{eq:N-EBENumFlux} has already been used. 

\subsection{FCFV discretisation}\label{sc:FCFV_Neumann}

The discretisation of the local problem given by equations~\eqref{eq:HDG-Nlocal} with a degree of approximation $k=0$ in each element for both $\bq_e$ and $u_e$ and also a degree of approximation $k=0$ in each face/edge for $\hu$ leads to the following system of equations for the local problem, for each element $e=1,\dotsc ,\numel$,
\begin{equation}\label{eq:Nlocal-system}
	\begin{bmatrix} |\Omega_e| \mat{I}_{\nsd}    & \vect{w}_e \\  
		\vect{w}_e^T & \alpha_e   \end{bmatrix}
	\begin{Bmatrix} \qe \\ \ue \end{Bmatrix} =
	\begin{Bmatrix} \vect{z}_e \\ \beta_e \end{Bmatrix} +
	\sum_{j \in \Iset} |\Gamma_{e,j}|
	\begin{Bmatrix} \bm{n}_j \\ \tau_j \end{Bmatrix}
	\uHj,
\end{equation}
where
\begin{equation}
  \begin{aligned}
	\alpha_e &= \sum_{j\in\Mset} |\Gamma_{e,j}| \tau_j, \;&
	\beta_e &= s_e |\Omega_e|  +  \sum_{j\in\Dset} |\Gamma_{e,j}| \tau_j u_{D,j} + \sum_{j\in\Nset} |\Gamma_{e,j}| t_j, \\
	\vect{w}_e &= \sum_{j\in\Nset} |\Gamma_{e,j}| \bn_j , \;&
	\vect{z}_e &= \sum_{j\in\Dset} |\Gamma_{e,j}| \bn_j u_{D,j} ,
  \end{aligned}
\end{equation}
and $\Mset = \left\{j \in \Aset \; | \; \Gamma_{e,j} \cap \partial \Omega = \emptyset \right\}$ is the set of indices corresponding to the interior faces of element $\Omega_e$. 

It is worth noting that, contrary to the FCFV formulation presented in Section~\ref{sc:Poisson}, the formulation with Neumann local problems leads to a set of local problems coupling the degrees of freedom of the solution and its gradient. However, the particular structure of the matrix appearing in the local problem of Equation~\eqref{eq:Nlocal-system} can be exploited to obtain an explicit formula for its inverse, namely
\begin{equation}
\begin{bmatrix} |\Omega_e| \mat{I}_{\nsd}    & \vect{w}_e \\  
\vect{w}_e^T & \alpha_e  \end{bmatrix}^{-1} 
=
\frac{1}{|\Omega_e| \vartheta_e}
\begin{bmatrix}  \vartheta_e \mat{I}_{\nsd} - \vect{w}_e \otimes \vect{w}_e    & |\Omega_e| \vect{w}_e \\  
|\Omega_e| \vect{w}_e^T & -|\Omega_e|^2
\end{bmatrix},
\end{equation}
where $\vartheta_e = \| \vect{w}_e \|^2_2 - |\Omega_e| \alpha_e $.

Therefore, the following explicit expressions of $\bq_e$ and $u_e$ in terms of $\hu$ are obtained
\begin{subequations}\label{eq:HDG-Poisson-DlocalK0explicitNeumann}
 \begin{multline}
	\qe =  |\Omega_e|^{-1} \vect{z}_e  - |\Omega_e|^{-1} \vartheta_e^{-1} \left( \vect{w}_e \otimes \vect{w}_e \right) \cdot \vect{z}_e + \vartheta_e^{-1} \beta_e \vect{w}_e + \\
	\sum_{j \in \Iset}  |\Gamma_{e,j}| \Big( |\Omega_e|^{-1} \bn_j  - |\Omega_e|^{-1} \vartheta_e^{-1} \left( \vect{w}_e \otimes \vect{w}_e \right) \cdot \bn_j  + \vartheta_e^{-1} \tau_j \vect{w}_e \Big) \uHj ,
	\label{eq:Poisson-locQ-Neumann}
  \end{multline}
	\begin{equation}
	\ue = \vartheta_e^{-1} \vect{z}_e \cdot \vect{w}_e - |\Omega_e| \vartheta_e^{-1} \beta_e +  	\sum_{j \in \Iset} |\Gamma_{e,j}| \Big(  \vartheta_e^{-1} \vect{w}_e \cdot \vect{n}_j - |\Omega_e| \vartheta_e^{-1} \tau_j  \Big) \uHj.
	\label{eq:Poisson-locU-Neumann}
	\end{equation}
\end{subequations}

Similarly, the discretisation of the global problem given by Equation \eqref{eq:Nglobal} with a degree of approximation $k=0$ leads to
\begin{equation}\label{eq:HDG-Poisson-DglobalK0Neumann}
\sum_{e=1}^{\numel}\Bigl\{
|\Gamma_{e,i}| \bn_i \cdot \qe + |\Gamma_{e,i}| \tau_i \ue - |\Gamma_{e,i}| \tau_i \uHi \Bigr\}
= 0, \quad \text{for $i \in \Mset$}.
\end{equation}

By inserting \eqref{eq:Poisson-locQ-Neumann} and \eqref{eq:Poisson-locU-Neumann} into \eqref{eq:HDG-Poisson-DglobalK0Neumann}, the following system of equations containing only $\hu$ as an unknown, on the internal faces, is obtained:
\begin{equation}
\mat{\widetilde{K}}\vect{\hu}=\vect{\widetilde{f}},
\end{equation}
where the global matrix $\mat{\widetilde{K}}$ and right hand side vector $\vect{\widetilde{f}}$ are computed by assembling the contributions given by
\begin{subequations}\label{HDG-Poisson-globalSystemNeumann}
  \begin{multline}
    \mat{\widehat{K}}^e_{i,j} :=  |\Gamma_{e,i}| 
    \Big[ |\Omega_e|^{-1} |\Gamma_{e,j}| \bn_i \cdot \bn_j 
            + \vartheta_e^{-1} |\Gamma_{e,j}| ( \bn_i \cdot \vect{w}_e ) \big( \tau_j - |\Omega_e| ( \bn_j \cdot \vect{w}_e ) \big)  \\
            + \tau_i \vartheta_e^{-1} |\Gamma_{e,j}| ( \bn_j \cdot \vect{w}_e ) - |\Omega_e| \vartheta_e^{-1} \tau_j - \tau_i \delta_{ij}  \Big]	 , 
  \end{multline}
  \begin{multline}
    \vect{\widehat{f}}^e_i :=  |\Gamma_{e,i}| 
    \Big[ \vartheta_e^{-1} ( \bn_i \cdot \vect{w}_e ) \big( |\Omega_e| ( \vect{z}_e \cdot \vect{w}_e ) - \beta_e \big) 
            -|\Omega_e|^{-1} \bn_i \cdot \vect{z}_e \\
            - \vartheta_e^{-1} \tau_i ( \vect{z}_e \cdot \vect{w}_e ) + |\Omega_e| \vartheta_e^{-1} \tau_i \beta_e \Big].
  \end{multline}
\end{subequations}

It is important to note that the formulation with Neumann local problems described here induces a higher computational cost compared to the formulation with Dirichlet local problems described in Section~\ref{sc:Poisson}.

\end{document}